\documentclass[12pt]{amsart}
\usepackage{amsmath}
\usepackage{amsxtra}
\usepackage{amscd}
\usepackage{amsthm}
\usepackage{amsfonts}
\usepackage{amssymb}
\usepackage{eucal}
\usepackage{mathrsfs}

\input prepictex
\input pictex
\input postpictex

\newtheorem{thm}{Theorem}[section]
\newtheorem{lem}{Lemma}[section]

\newtheorem{prop}{Proposition}[section]

\theoremstyle{definition}

\theoremstyle{remark}
\newtheorem{rem}{Remark}[section]

\numberwithin{equation}{section}

\begin{document}

\newcommand{\thmref}[1]{Theorem~\ref{#1}}
\newcommand{\secref}[1]{Section~\ref{#1}}
\newcommand{\lemref}[1]{Lemma~\ref{#1}}
\newcommand{\propref}[1]{Proposition~\ref{#1}}
\newcommand{\corref}[1]{Corollary~\ref{#1}}
\newcommand{\remref}[1]{Remark~\ref{#1}}
\newcommand{\eqnref}[1]{(\ref{#1})}
\newcommand{\exref}[1]{Example~\ref{#1}}

\newcommand{\nc}{\newcommand}
\nc{\on}{\operatorname} \nc{\Z}{{\mathbb Z}} \nc{\C}{{\mathbb C}}
\nc{\R}{{\mathbb R}} \nc{\boldD}{{\mathbb D}} \nc{\oo}{{\mf O}}
\nc{\N}{{\mathbb N}} \nc{\bib}{\bibitem} \nc{\pa}{\partial}
 \nc{\rarr}{\rightarrow}
\nc{\larr}{\longrightarrow} \nc{\al}{\alpha} \nc{\ri}{\rangle}
\nc{\lef}{\langle} \nc{\W}{{\mc W}} \nc{\gam}{\ol{\gamma}}
\nc{\Q}{\ol{Q}} \nc{\q}{\widetilde{Q}} \nc{\la}{\lambda}
\nc{\ep}{\epsilon} \nc{\g}{{\mf g}} \nc{\h}{\mf h} \nc{\n}{\mf n}
\nc{\bb}{\mf b}  \nc{\G}{{\mf g}} \nc{\D}{\mc D} \nc{\Li}{{\mc L}}
\nc{\La}{\Lambda} \nc{\is}{{\mathbf i}} \nc{\V}{\mf V}
\nc{\bi}{\bibitem} \nc{\NS}{\mf N}
\nc{\dt}{\mathord{\hbox{${\frac{d}{d t}}$}}} \nc{\E}{\mc E}
\nc{\ba}{\tilde{\pa}} \nc{\half}{\frac{1}{2}}
\def\smapdown#1{\big\downarrow\rlap{$\vcenter{\hbox{$\scriptstyle#1$}}$}}
\nc{\mc}{\mathcal} \nc{\mf}{\mathfrak} \nc{\ol}{\fracline}
\nc{\el}{\ell} \nc{\etabf}{{\bf \eta}} \nc{\zetabf}{{\bf
\zeta}}\nc{\x}{{\bf x}} \nc{\xibf}{{\bf \xi}} \nc{\y}{{\bf y}}
\nc{\z}{{\bf z}} \nc{\parth}{\partial_{\theta}}
\nc{\cwo}{\C[w]^{(1)}} \nc{\cwe}{\C[w]^{(0)}}
\nc{\hf}{\frac{1}{2}} \nc{\hsdzero}{{}^0\widehat{\sd}}
\nc{\gltwo}{{\rm gl}_{\infty|\infty}}
\nc{\btwo}{{B}_{\infty|\infty}} \nc{\htwo}{{\h}_{\infty|\infty}}
\nc{\hglone}{\widehat{\rm gl}_{\infty}} \nc{\hgltwo}{\widehat{\rm
gl}_{\infty|\infty}} \nc{\hgltwof}{\widehat{\rm
gl}^f_{\infty|\infty}} \nc{\hbtwo}{\hat{B}_{\infty|\infty}}
\nc{\hhtwo}{\hat{\h}_{\infty|\infty}} \nc{\glone}{{\rm gl}_\infty}
\nc{\gl}{{\rm gl}} \nc{\pd}{\mc P} \nc{\hpd}{\widehat{\pd}}
\nc{\co}{\mc O} \nc{\Oe}{\co^{(0)}} \nc{\Oo}{\co^{(1)}}
\nc{\sdzero}{{}^0{\sd}} \nc{\hz}{\hf+\Z} \nc{\vac}{|0 \rangle}
\nc{\K}{\mf k} \nc{\bhf}{\bf\hf}

\def\mp#1#2#3{#1_{#2}^{#3}\mapsto\frac{\partial}{\partial#1_{#2}^{#3}},\quad
\frac{\partial}{\partial#1_{#2}^{#3}}\mapsto #1_{#2}^{#3}}
\def\cas#1#2#3#4#5#6#7#8{#1_#2^{#3,#4}=\left\{%
\begin{array}{ll}
    \frac{\partial}{\partial#5_{#2}^{#4}}, & #2#7 0, \\
    #6_{-#2}^{#4}, &#2#8 0 . \\
\end{array}%
\right.    }
\def\casa#1#2#3#4#5#6#7#8{#1_#2^{#3,#4}=\left\{%
\begin{array}{ll}
    -\frac{\partial}{\partial#5_{#2}^{#4}}, & #2#7 0, \\
    -#6_{-#2}^{#4}, &#2#8 0 . \\
\end{array}%
\right.    }
\def\casb#1#2#3#4#5#6#7#8{#1_#2^{#3,#4}=\left\{%
\begin{array}{ll}
    \frac{\partial}{\partial#5_{#2}^{#4}}, & #2#7 0, \\
    -#6_{-#2}^{#4}, &#2#8 0 . \\
\end{array}%
\right.    }

\nc{\fb}{{\mathfrak b}} \nc{\fg}{{\mathfrak g}}
\nc{\fh}{{\mathfrak h}}  \nc{\fk}{{\mathfrak k}}
\nc{\fl}{{\mathfrak l}} \nc{\fn}{{\mathfrak n}}
\nc{\fp}{{\mathfrak p}} \nc{\fu}{u} \nc{\cF}{{\mathcal F}}
\nc{\fsp}{{\mathfrak {sp}}} \nc{\fso}{{\mathfrak {so}}}
 \nc{\F}{{\mf F}^d}
 \nc{\Fhf}{{\mf F}^{d+\hf}}
 \nc{\Cxall}{\C[\x,\y,\etabf, \zetabf]}
 \nc{\Dxall}{\boldD[\x,\y,\etabf, \zetabf]}
 \nc{\cA}{\mc A} \nc{\cC}{\mc C} \nc{\cD}{\mc D}
 \nc{\hA}{\widehat{\mc A}} \nc{\hB}{\widehat{\mc B}} \nc{\hC}{\widehat{\mc C}} \nc{\hD}{\widehat{\mc D}}
 \nc{\sgn}{\rm sgn} \nc{\te}{{\tilde e}} \nc{\ccA}{{\mathscr A}^d}
 \nc{\ccAhf}{{\mathscr A}^{d+\hf}}
 \nc{\tE}{{\widetilde E}}
 \nc{\tM}{{\widetilde M}}
  \nc{\tm}{{\tilde m}}
 \nc{\tX}{{\widetilde X}}
\nc{\rank}{\mbox{rank}}
\def\xihalfge{\xi_{\frac{1}{2}}>\xi_{\frac{3}{2}}>\cdots>\xi_{r-\frac{1}{2}}}
\def\xionege{\xi_1>\xi_2>\cdots>\xi_r}
\def\xineghalfge{\xi_{s+\frac{1}{2}}>\xi_{s+\frac{3}{2}}>\cdots>\xi_{-\frac{1}{2}}}
\def\xinegzeroge{\xi_{s+1}>\xi_{s+2}>\cdots>\xi_{0}}
\def\xihalf{\xi_{\frac{1}{2}},\xi_{\frac{3}{2}},\cdots,\xi_{r-\frac{1}{2}}}
\def\xione{\xi_1,\xi_2,\cdots,\xi_r}
\def\xineghalf{\xi_{s+\frac{1}{2}},\xi_{s+\frac{3}{2}},\cdots,\xi_{-\frac{1}{2}}}
\def\xinegzero{\xi_{s+1},\xi_{s+2},\cdots,\xi_{0}}
\def\xiall{\xineghalf \mid \xinegzero\mid \xihalf\mid \xione}

\def\xihalfone{\xi_{\frac{1}{2}},\xi_{1},\cdots,\xi_{r-\frac{1}{2}},\xi_r}
\def\xineghalfzero{\xi_{s+1},\xi_{s+\frac{3}{2}},\cdots,\xi_{-\frac{1}{2}},\xi_0}
\setlength{\unitlength}{0.5cm}

\advance\headheight by 2pt

\title[Representations of Lie superalgebras of infinite rank]
{Quasi-finite representations, free field realizations,
and character formulae of Lie  superalgebras of infinite rank}
\thanks{{\bf 2000 Mathematics Subject Classification}:
17B65, 17B10}

\author[Ngau Lam]{Ngau Lam$^{1}$}
\thanks{$^{1}$Partially supported by NSC-grant 92-2115-M-006-016 of the
R.O.C.}
\address{Department of Mathematics, National Cheng Kung University, Tainan,
Taiwan 701} \email{nlam@mail.ncku.edu.tw}

\author[R.~B.~Zhang]{R.~B.~Zhang$^{2}$}
\thanks{$^{2}$Supported by the Australian Research Council.}
\address{School of Mathematics and Statistics, University of Sydney,
New South Wales 2006, Australia} \email{rzhang@maths.usyd.edu.au}

\begin{abstract} We classify the quasi-finite irreducible highest
weight modules over the infinite rank Lie  superalgebras
$\hgltwo$, $\hC$ and $\hD$, and determine the necessary and sufficient
conditions for quasi-finite irreducible highest
weight modules to be unitarizable with respect to natural
$\ast$-structures of the Lie superalgebras. The
unitarizable irreducible modules are constructed in terms of Fock
spaces of free quantum fields, and explicit formulae for their
formal characters are also obtained by investigating Howe dualities
between the infinite rank Lie  superalgebras and
classical Lie groups.

\bigskip
\noindent{\bf Key words}: Infinite dimensional Lie superalgebras,
quasi-finite representations, unitarizable representations, character
formulae.
\end{abstract}
\maketitle


\section{Introduction}
Supersymmetry permeated many areas of mathematics in the last decade,
producing deep results such as the Seiberg-Witten theory and
mirror symmetry. In all applications, supersymmetry manifests itself as
concrete representations of the relevant Lie superalgebras \cite{K}. Thus it is of
central importance to develop the representation theory of Lie superalgebras
in order to use supersymmetry as a tool to
address problems in other areas.

In this paper we investigate the representation theory of the Lie
superalgebra $\hgltwo$ and its $osp$-type Lie sub superalgebras.
These Lie superalgebras constitute a class of $\hf\Z$-graded infinite rank Lie
superalgebras arising from central extensions of Lie superalgebras
of complex matrices of infinite size. The Lie superalgebra $\hgltwo$
and its $osp$-type Lie sub superalgebras featured very prominently
in the study \cite{CW3} of the super $W_{1+\infty}$ algebra, i.e.,
the central extension of the superalgebra of differential operators
on the super circle, which plays a fundamental role in
conformal filed theory and the theory of superstrings.
Also, it was demonstrated in \cite{CL} that the representation theory
of the infinite rank Lie superalgebras is intimately related to
that of affine Kac-Moody superalgebras arising from central extensions
of the loop algebras of finite dimensional simple Lie superalgebras.
In this paper we shall focus on $\hgltwo$ and its subalgebras $\hC$ and $\hD$
(see \secref{sect superalgebras} for their definitions).
Aspects of a $\hB$ type subalgebra of $\hgltwo$ were studied in \cite{CW3}.

Recall that the infinite dimensional Lie algebra $\hglone$ and its
various subalgebras were extensively studied in \cite{KR, KWY, W1, W2} in
relation to the $W_{1+\infty}$ algebra. In particular, the notion of
quasi-finite modules \cite{KR} over infinite dimensional graded Lie
(super)algebras were introduced.
Such modules are close to finite dimensional representations of finite
dimensional Lie (super)algebras in spirit. In our context,
a $\hf\Z$-graded module $M=\oplus_{j\in\hf\Z} M_j$ over a
$\hf\Z$-graded Lie superalgebra will be called quasi-finite if all
its homogeneous subspaces $M_j$ are finite dimensional.
One of our results in the present paper is the classification of all the
quasi-finite irreducible highest weight modules over $\hgltwo$,
$\hC$ and $\hD$.

It is well known that the energy of a quantum system is always bounded below. Also,
the space of the physical states of the quantum system always admits a
positive definite contravariant Hermitian form
as required by the probabilistic interpretation of quantum theory.
Therefore, the representations of Lie superalgebras,  which are
potentially useful in quantum physics, are the unitarizable highest weight (or
lowest weight) representations. Another result of the present paper
is the classification of the unitarizable quasi-finite irreducible highest weight
modules over $\hgltwo$, $\hC$ and $\hD$ with respect to some natural
$\C$-conjugate linear anti-involutions of these Lie superalgebras.

We analyse the unitarizable irreducible quasi-finite highest weight modules in some detail.
The main results obtained are the following.
We first realize these irreducible representations on
Fock spaces of free quantum fields. We then prove generalized Howe dualities
between the infinite rank Lie  superalgebras and certain classical Lie groups.
This way we are able to set up one to one correspondences between the unitarizable
irreducible quasi-finite highest weight modules of the infinite rank Lie
superalgebras and the finite dimensional irreducible representations of the associated
classical groups. Finally we derive explicit formulae for the formal characters
of the unitarizable quasi-finite irreducible highest weight modules over the Lie
superalgebras $\hC$ and $\hD$. (We recall that the formal
characters of the unitarizable irreducible modules over $\hgltwo$
were obtained in \cite{CL}.)

The method used here for the construction of the character formulae is
a generalization of that developed in \cite{CL, CZ, CLZ}, which relies
in an essential way on Howe dualities \cite{H1, H2}. Howe dualities for Lie superalgebras
were known in the original paper of Howe \cite{H1}, and also in \cite{S1, S2},
and were further investigated in \cite{W1, W2, N, OP, CW1, CW2}. Recent investigations on
Howe dualities led to a thorough understanding of the Segal-Shale-Weil
representations of Lie superalgebras \cite{CZ, CL}, in particular,
the construction of character formulae for them.  In \cite{H}and \cite{W1, W2}
Howe dualities were established respectively in the contexts of
affine Kac-Moody algebras and infinite rank Lie algebras.
The Howe dualities obtained in the present paper are generalizations
of those studied by Wang in \cite{W1, W2}.

The arrangement of the paper is as follows. \secref{Preliminaries}
provides some background material on generalized partitions and unitarizable
modules over Lie superalgebras. The material will be used throughout
the paper.  \secref{sect superalgebras} examines central extensions of the
Lie superalgebra $\gltwo$ of infinite matrices and its $osp$-type subalgebras,
and gives the definitions of the Lie superalgebra $\hgltwo$, and its
subalgebras $\hA$,  $\hC$ and $\hD$. The remaining three sections
constitute the main body of the paper. In
\secref{sect_quasi-finiteness}, we classify the quasi-finite
irreducible highest weight modules over the Lie superalgebra $\hgltwo$,
and its subalgebras $\hA$,  $\hC$ and $\hD$. In \secref{sect
unitarizable}, we classify the unitarizable quasi-finite
irreducible highest weight modules over these Lie
superalgebras with respect to specific $\ast$-structures, and
construct Fock space realizations of the unitarizable irreducible
modules. Generalized Howe dualities between the infinite rank
Lie  superalgebras and classical Lie groups will also be
established in this section, which are used in \secref{sect character} to
derive character formulae for the unitarizable quasi-finite
irreducible highest weight modules over $\hC$ and $\hD$.

\section{Preliminaries}\label{Preliminaries} We work on the field
$\C$ of complex numbers throughout the paper. For any vector space
$V$, we shall denote its dual space by $V^*$.

\subsection{Shifted Frobenius notation for generalized partitions}
\label{Frobenius} By a {\em partition $\la$ of length $d$} we mean
a non-increasing finite sequence of non-negative integers
$(\la_1,\cdots,\la_d )$ and shall use $l(\la)$ to denote the length of $\la$.
We will let $\la'$ denote the {\em transpose} of the partition
$\la$. We define the {\em rank} of a partition
$\la=(\la_1,\cdots,\la_d )$, denoted by $\rank(\la)$, to be the
largest integer $i$, for which $\la_i\geq i$. Note that
$\rank(\la)=\rank(\la')\le d$. For example, if $\la=(4,3,1,0,0)$,
then $l(\la)=5$, $\la'=(3,2,2,1)$, and $\rank(\la)=2$. By a {\em
generalized partition of length $d$}, we shall mean a
non-increasing finite sequence of integers $(\la_1,\cdots,\la_d )$
and the length of $\la$ is also denoted by $l(\la)$. A generalized
partition of  $\la=(\la_1,\dots,\la_d)$ is called a generalized
partition of non-positive integers if $\la_i\leq 0$ for all $i$.
Corresponding to each generalized partition
$\la=(\la_1,\dots,\la_d)$, we will define
$\la^*:=(-\la_d,\dots,-\la_1)$. Then $\la^*$ is also a generalized
partition. In particular, if
$\la=(\la_1,\dots,\la_d)$ is a generalized partition of non-positive integers,
then $\la^*=(-\la_d,\dots,-\la_1)$ is a partition. In this case, we
define the {\em rank}, $\rank(\la )$, of $\la$ by $\rank(\la ):=-\rank(\la^* )$.
We also set $\la'_{j}:=-(\la^*)'_{-j}$ for all $j\in\{ -1,-2,\cdots,\la_d\}$.

Each generalized partition $\la=(\la_1,\cdots,\la_d)$ of length
$d$ can be uniquely expressed as $\la=\la^+ + \la^-$, with
\begin{align*}\la^+:=(\hbox{max}\{\la_1,0\},\cdots,\hbox{max}\{\la_d,0\}),
\\
\la^-:=(\hbox{min}\{\la_1,0\},\cdots,\hbox{min}\{\la_d,0\}).
\end{align*}
Note that $\la^+$ is a partition of length $d$ , while $\la^-$ is
a generalized partition of non-positive integers of length $d$.
Furthermore,
\begin{eqnarray}\label{depth}
\mbox{depth of $\la^+$}\, + \ \mbox{depth of $(\la^-)^*$}&\le & d,
\end{eqnarray}
where the depth of a partition is the number of positive integers
in it.( Note that the depth of a partition $\la$ equals $\la'_1$.)

Now we will define the shifted Frobenius notation for generalized
partitions (see \cite M) which is very useful for describing the
highest weights of unitarizable irreducible quasi-finite
modules over $\hgltwo$. Given a partition
$\la=(\la_1,\dots,\la_d)$ of length $d$ and $\rank(\la)=r>0$, we let
$\xi_i:=\la_{i+\frac{1}{2}}-i+\frac{1}{2}$ and $\xi_j:=\la_j'-j$,
for $i\in \frac{1}{2}+\Z_+$ with $\frac{1}{2}\leq i
\leq r-\frac{1}{2} $,  and $j\in \N$ with $1\leq j \leq r$. We have
\begin{equation}\label{xipos}
\xihalfge > 0,  \quad \xionege \geq 0.
\end{equation}
The {\em shifted Frobenius notation for the partition $\la$}
is given by
$$
  \Lambda(\la):=(\xihalf\mid \xione).
$$
Clearly, we have
\begin{equation}\label{d}
\xi_1+\mbox{min}\{\xi_\hf,1\}\le d,
\end{equation}
When $\la=(0,0,\cdots,0)$, we define
$\Lambda(0,0,\cdots,0):=(0,0)$. Note that \eqnref{xipos} and
\eqnref{d} implies $r\leq d$.

Conversely, if two finite sequences $\xihalf$ and $\xione$
of non-negative integers of length $d$ satisfy \eqnref{xipos} and
\eqnref{d}, we may regard them as the shifted Frobenius notation
of a unique partition of length $d$, which we will denote by
$F(\xihalf\mid \xione)$. We put $F(0,0):=(0,0,\cdots,0)$. Thus we
have a one-to-one correspondence between the set of all partitions of
length $d$ and the set of all pairs of finite sequences of
non-negative integers, $\xihalf$ and $\xione$ satisfying
\eqnref{xipos} and \eqnref{d}, and also the requirement that
\begin{equation}\label{xigeq}
\begin{array}{ll}
\qquad\qquad \xi_{r-\hf}=0 \ \mbox{if and only if }\ r=1\
\mbox{and}\ \xi_\hf=\xi_1=0.
\end{array}
\end{equation}

Similarly, given any non-zero generalized partition
$\la=(\la_1,\dots,\la_d)$ of non-positive integers with $\rank
(\la)=s$, the {\em shifted Frobenius notation for the generalized
partition $\la$ of non-positive integers}, also denoted by
$\Lambda(\la)$, is defined by
$$
\Lambda(\la):=(\xineghalf\mid\xinegzero)
$$
where $\xi_i:=\la_{i-1}'-i$ and $\xi_j:=\la_{d+\hf
+j}-j+\frac{1}{2}$ for all $i\in\{0,-1,-2,\cdots,s+1\}$ and
$j\in\{-\frac{1}{2},-1\frac{1}{2},\cdots,s+\frac{1}{2}\}$. We also
define $\Lambda(0,0,\cdots,0):=(0,0)$. Similarly we have a
one-to-one correspondence between the set of all generalized
partitions of non-positive integers of length $d$ and the set of
all pairs of finite sequences of non-positive integers,
$\xineghalf$ and $\xinegzero$,  satisfying the following
conditions
\begin{equation}\label{xineggeq}
\begin{array}{ll}
0\geq \xineghalfge,  \quad \quad 0\geq\xinegzeroge ,\\
\xi_{s+1}=0 \ \mbox{if and only if }
s=-1\
 \mbox{and}\  \xi_{-\hf}=\xi_0=0,
\end{array}
\end{equation}
and
\begin{equation}\label{d-neg}
\xi_0\le d.
\end{equation}

Now we define the shifted Frobenius notation for the generalized
partitions as follows.  For a nonzero generalized partition $\la$
of length $d$, the {\em shifted Frobenius notation for the
generalized partition $\la$}, also denoted by $\Lambda(\la)$, is
defined by
$$
\Lambda(\la):=(\Lambda(\la^-)|\Lambda(\la^+)).
$$
Similarly we have a one-to-one correspondence between the set of
all generalized partitions of length $d$ and the set of all quartets
of finite sequences of integers, $\xineghalf$; $\xinegzero$;
$\xihalf$ and $\xione$ satisfying \eqref{xipos}, \eqnref{xigeq},
\eqnref{xineggeq} and
\begin{equation}\label{length}
\min\{\xi_\hf,1\}+\xi_1-\xi_0\leq d.
\end{equation}
We will denote by
$$
F(\xineghalf\mid \xinegzero\mid \xihalf\mid\xione),
$$
the unique generalized partition
corresponding to the quartet of finite sequences of integers $\xineghalf$;\
$\xinegzero$;\  $\xihalf$ and $\xione$  satisfying \eqref{xipos}, \eqnref{xigeq},
\eqnref{xineggeq} and \eqnref{length}.

\subsection{Unitarizable modules}
Let us recall some basic facts about $\ast$-superalgebras and
their unitarizable modules.  A {\em $\ast$-superalgebra} is an
associative superalgebra $A$ together with an anti-linear
anti-involution $\omega: A\rightarrow A$. Here we should emphasize
that for any $a, b\in A$, we have $\omega(a b)=
\omega(b)\omega(a)$, where no sign factors are involved. A
$\ast$-superalgebra homomorphism $f: (A, \omega) \rightarrow (A',
\omega')$ is a superalgebra homomorphism obeying $f\circ \omega =
\omega'\circ f$. Let $(A, \omega)$ be a $\ast$-superalgebra, and
let $M$ be a $\Z_2$-graded $A$-module. A Hermitian form $\langle\
\cdot\ | \ \cdot \ \rangle$ on $M$ is said to be {\em
contravariant} if $\langle a v | v'\rangle = \langle v | \omega(a)
v'\rangle$, for all $a\in A$, $v, v'\in M$.
An $A$-module $M$ is called {\em unitarizable} if $M$
admits a positive definite contravariant Hermitian form.

Let $\g$ be a Lie superalgebra together with an anti-linear
anti-involution $\omega$ (i.e. $\omega$ is an anti-linear map
satisfying $\omega([x,y])=[\omega(y),\omega(x)]$ for all $x,y\in
\g$. In this case, we also call $\omega$ a $\ast$-structure of $\g$.).
Let $M$ be a $\g$-module. A Hermitian form $\langle\ \cdot\
| \ \cdot \ \rangle$ on $M$ is said to be {\em contravariant} if
$\langle x v | v'\rangle = \langle v | \omega(x) v'\rangle$, for
all $x\in \g$, $v, v'\in M$. When the Hermitian form
$\langle\ \cdot\ | \ \cdot \ \rangle$ is positive definite,
we define $\|u\| :=\sqrt{\langle u|u\rangle} $ for all $u\in M$.
A $\g$-module $M$ is called {\em unitarizable} if
$M$ admits a positive definite contravariant Hermitian form. The
anti-linear anti-involution $\omega$ can be naturally extended to
an anti-linear anti-involution, also denoted by
$\omega$, on the universal enveloping algebra
${\mc U}(\g)$ of $\g$, making $({\mc U}(\g), \omega)$  a
$\ast$-superalgebra. Moreover, a $\g$-module $M$ is unitarizable
if and only if it is a unitarizable ${\mc U}(\g)$-module.
(Throughout the paper, ${\mc U}( \mf s)$ stands for the universal
enveloping algebra of the Lie superalgebra $\mathfrak s$.)

\section{Lie  superalgebras of infinite rank} \label{sect superalgebras}
We present here the infinite rank
Lie superalgebras to be studied in this paper.

Consider the infinite-dimensional complex superspace
$\C^{\infty|\infty}$ with a basis $\{ e_j \mid  j\in \Z\}$ for the
even subspace, and a basis $\{ e_r \mid r\in \hf+\Z\}$ for the odd
subspace. We introduce a $\hf\Z$-gradation on $\C^{\infty|\infty}$  by
setting the degree of $e_p$ equal
to $-p$ for all $p\in \hf \Z$. For any $p, q\in\hf\Z$, let
let $e_{p q}$ be the
endomorphism of $\C^{\infty|\infty}$ defined by
$e_{p q}(e_r)=\delta_{q r} e_p$. Then
$T$ is a homogeneous endomorphism on $\C^{\infty|\infty}$ of degree $p$ if
and only if $T=\sum_{j\in \hf \Z }a_j e_{j-p,j}$, where $ \ a_j\in
\C$. Denote by $(M_{\infty|\infty})_p$ the set of all endomorphisms
of $\C^{\infty|\infty}$ of degree $p$, and let
$M_{\infty|\infty}:=\oplus_{p\in\hf\Z}(M_{\infty|\infty})_p$. Then
$M_{\infty|\infty}$ is a $\hf\Z$-graded associative superalgebra,
which also acquires a Lie superalgebra structure with the usual Lie
super-bracket
\begin{equation}
[A, B]:= A B - (-1)^{4 deg(A) deg(B)} B A,
\end{equation}
where $deg(A)$ and $deg(B)$ are the degrees of $A$ and $B$
respectively. We shall denote this Lie superalgebra by
$\gltwo:=\oplus_{p\in\hf\Z}(\gltwo)_p$. Note
that the subspace $\gltwo^f$ generated by $\{e_{ij}|i,j\in\half\Z\}$ is
a subalgebra of $\gltwo$. By arranging the basis elements of
$\C^{\infty|\infty}$ in strictly increasing order,  any endomorphism
of $\C^{\infty|\infty}$ may be written as an infinite-sized square
matrix with coefficients in $\C$. Thus
\begin{equation*}
\gltwo:=\{(a_{ij}), i,j\in\hf\Z|\ a_{ij}=0\ {\rm for}\ |j-i|>>0\}.
\end{equation*}

The Lie superalgebra $\gltwo$ has a central extension by a
non-trivial two co-cycle. Let $J=\sum_{r\le 0} e_{r r}$. Define
\begin{equation}\label{centralext}
\alpha(A,B):={\rm Str}([J,A]B),\quad A,B\in\gltwo,
\end{equation}
where $\rm Str$ stands for the supertrace defined for a matrix
$D=(d_{ij})\in\gltwo$
by ${\rm Str}(D)=\sum_{r\in\hf\Z}(-1)^{2r}d_{rr}$, provided that
the infinite sum is not divergent.
Then $\alpha(A,B)$ is well behaved for all $A,B\in\gltwo$, and
indeed defines a two co-cycle.

We denote by $\hgltwo$ the central extension of $\gltwo$ by the
even central element $C$ associated with this two co-cycle.
By setting the degree of $C$ equal to $0$, the Lie superalgebra
$\hgltwo$ acquires a $\half\Z$-gradation from that of
$\gltwo$. Let
\begin{equation*}
\hgltwof:=\{(a_{ij})\in \hgltwo | \ \mbox{finitely many of the}\
a_{ij} \ \mbox{are non-zero}\ \}\oplus\C C.
\end{equation*} It is easy to see that
$\hgltwof$ is a $\hf \Z$-graded subalgebra of $\hgltwo$.

Let us now introduce the Lie sub superalgebra $\cA$ of $\gltwo$
defined by
\begin{equation*}
\cA:=\{(a_{ij})\in \gltwo | \  a_{ij}=0 \ \mbox{if}\ i=0 \
\mbox{or}\ j=0 \}.
\end{equation*}
It also admits a central extension by an even central element $C$
associated with the two co-cycle \eqnref{centralext}. We shall
denote the central extension of $\cA$ by $\hA$. Then the Lie superalgebra
$\hA$ also acquires a $\half\Z$-gradation from that of
$\hgltwo$ by declaring $C$ to have degree $0$.

An alternative way to describe the Lie superalgebra $\cA$ is as
follows. Consider the infinite-dimensional complex superspace
$\C^{\infty|\infty}$ with even basis $\{ e_j \mid  j\in \Z^*\}$
and odd basis $\{ e_r \mid r\in \hf+\Z\}$. Then $\cA$ is the Lie
superalgebra of graded endomorphisms of $\C^{\infty|\infty}$.

Let us now construct a Lie sub superalgebra $\cC$ of $\cA$.
Introduce a non-degenerate skew-supersymmetric bilinear form
$(\cdot|\cdot)$ on $\C^{\infty|\infty}$ defined by
\begin{eqnarray*}
&(e_i|e_j)=-(e_j|e_i)={\sgn(i)}\delta_{i,-j},\qquad i,j\in\Z^*;\\
&(e_r|e_s)=(e_s|e_r)=\delta_{r,-s},\qquad r,s\in\hf+\Z; \\
&(e_i|e_r)=(e_r|e_i)=0,\qquad i\in\Z^*, r\in\hf+\Z;
\end{eqnarray*}
where $\sgn(i):=+1$ if $i\in \hf\N$ and $\sgn(i):=-1$ if $i\in
-\hf\N$. We define the Lie superalgebra
$\cC =\cC_0 \oplus \cC_1$ to be the $\hf\Z$-graded Lie sub superalgebra of
$\cA$ preserving this form, i.e.
$$\cC_\epsilon=\{A\in\cA_{\epsilon}|(Av|w)=-(-1)^{\epsilon|v|}(v|Aw)\},
 \;\;\epsilon =0,1,
$$
where $|v|$ denotes the parity of $v\in \C^{\infty|\infty}$, namely,
$|v|=0$ (respectively $1$) if $v$ belongs to the even (respectively
odd) homogeneous subspace of $\C^{\infty|\infty}$.
Then $\cC$ is a Lie superalgebra of type $SPO$. It is easy to see that
the subspace $\cC^f$ spanned by the following elements is
a subalgebra of $\cC$ ($i,j\in\Z^*$, $r,s\in\hf+\Z$):
\begin{align*}
&\te_{i,j}:=-\te_{-j,-i}:=e_{i,j}-e_{-j,-i}, \quad  ij>0 \ (\rm{i.e.}, \ i, j>0 \ \rm{or} \ i, j <0);\\
&\te_{i,j}:=\te_{-j,-i}:=e_{i,j}+e_{-j,-i}, \quad  ij<0 \ (\rm{i.e.}, \ i, -j>0 \ \rm{or} \ i, -j <0) ;\\
&\te_{r,s}:=-\te_{-s,-r}:=e_{r,s}-e_{-s,-r};\\
&\te_{i,r}:=\te_{-r,-i}:=e_{i,r}+e_{-r,-i}, \quad i>0;\\
&\te_{i,r}:=-\te_{-r,-i}:=e_{i,r}-e_{-r,-i}, \quad i<0.
\end{align*}
Note that $\cC^f_0$ is a direct sum of an infinite dimensional
symplectic Lie algebra and an infinite dimensional orthogonal Lie
algebra. Let $\hC$ denote the central extension of $\cC$ by an even
central element $C$ associated with the two-cocycle \eqnref{centralext}.
By setting the degree of $C$ to zero, $\hC$ becomes a $\half\Z$-graded Lie superalgebra,
with the gradation compatible with that of $\hA$.

Consider the non-degenerate supersymmetric bilinear form
$(\cdot|\cdot)$ on $\C^{\infty|\infty}$ defined by
\begin{eqnarray*}
&(e_i|e_j)=(e_j|e_i)=\delta_{i,-j},\qquad i,j\in\Z^*;\\
&(e_r|e_s)=-(e_s|e_r)={\sgn(r)}\delta_{r,-s},\qquad r,s\in\hf+\Z; \\
&(e_i|e_r)=(e_r|e_i)=0,\qquad i\in\Z^*, r\in\hf+\Z.\\
\end{eqnarray*}
We define the Lie superalgebra $\cD =\cD_0
\oplus \cD_1$ to be the subalgebra of $\cA$ preserving this
form, i.e.
$$\cD_\epsilon=\{A\in\cA_{\epsilon}|(Av|w)=-(-1)^{\epsilon|v|}(v|Aw)\},
 \;\;\epsilon =0,1.
$$
This is a Lie superalgebra of type $OSP$. It is easy to see that
the subspace $\cD^f$ spanned by the following elements is a
subalgebra of $\cD$ ($i,j\in\Z^*$, $r,s\in\hf+\Z$):
\begin{align*}
&\te_{i,j}:=-\te_{-j,-i}:=e_{i,j}-e_{-j,-i};\\
&\te_{r,s}:=-\te_{-s,-r}:=e_{r,s}-e_{-s,-r},\quad rs>0 \ (\rm{i.e.}, \ r, s>0 \ \rm{or} \ r, s <0);\\
&\te_{r,s}:=\te_{-s,-r}:=e_{r,s}+e_{-s,-r},\quad  rs<0 \ (\rm{i.e.}, \ r, -s>0 \ \rm{or} \ r, -s <0);\\
&\te_{i,r}:=\te_{-r,-i}:=e_{i,r}+e_{-r,-i}, \quad r>0; \\
&\te_{i,r}:=-\te_{-r,-i}:=e_{i,r}-e_{-r,-i}, \quad r<0.
\end{align*}
Note that $\cD^f_0$ is a direct sum of an infinite dimensional
symplectic Lie algebra and an infinite dimensional orthogonal Lie
algebra. The Lie superalgebra $\cD$ has a central extension by an even
central element $C$ associated with the two-cocycle given in \eqnref{centralext}.
We shall denote this central extension by $\hD$. The Lie
superalgebra $\hD$ also has a $\half\Z$-gradation compatible with that
of $\hA$ with $C$ being of degree zero.

\begin{rem}\label{embeddings}
Both $\hC$ and $\hD$ are $\hf\Z$-graded Lie sub superalgebras of $\hA$.
Thus the triangular decomposition $\hA=\hA_+\oplus\hA_0\oplus\hA_-$ of
$\hA$ leads to natural triangular decompositions of $\hC$ and $\hD$:
\begin{equation*}
\hC=\hC_+\oplus\hC_0\oplus\hC_-, \quad \hD=\hD_+\oplus\hD_0\oplus\hD_-,
\end{equation*}
where $\hC_\varpi=\hC\cap\hA_\varpi$, \ $\hD_\varpi=\hD\cap\hA$, for
$\varpi$ being $+$, $-$ and $0$. This will be of considerable
importance when we develop the representation theory of these Lie
superalgebras.
\end{rem}
\begin{rem} For $\hat\g$ being $\hA$, $\hC$ or $\hD$, we shall use the notation
$\hat\g^f$ to denote the $\hf\Z$-graded subalgebra $\hat\g\cap\hgltwo^f$.
\end{rem}

\section{Criterion for quasi-finiteness of modules}
\label{sect_quasi-finiteness}

In this section we give a complete classification of all the
quasi-finite irreducible highest weight modules over the Lie
superalgebras discussed in \secref{sect superalgebras}.

\subsection{Quasi-finite modules}
Let $\g = \oplus_{j \in \hf \Z} {\g}_j$ (possibly $\dim{\g_j}=\infty$)
be a $\hf \Z$-graded Lie superalgebra,
with the even subspace $\oplus_{j \in \Z} {\g}_j$,  and odd subspace
$\oplus_{j \in \hf + \Z} {\g}_j$. We assume that $\g_0$ is abelian.
We have the triangular decomposition
\begin{equation*}
\g=\g_{-}\oplus \g_0\oplus \g_+,\quad  \rm{with} \  \ \g_{\pm} =
\oplus_{r>0}\g_{\pm r}.
\end{equation*}
A $\g$-module $M =
\oplus_{j \in \hf \Z} M_j $ is {\em graded} if ${\g}_i M_j
\subseteq M_{i+j}$. A vector $v\in M$ is called {\em homogeneous}
of degree $j$ if $v\in M_j$ for some $j\in \hf\Z$. Following the
terminology of Kac and Radul \cite{KR}, we shall call $M$ {\em
quasi-finite} if $ \dim M_j < \infty$ for all $j\in\hf\Z$.

A $\g$-module $M$ is called a {\em highest weight module} with highest
weight $\xi \in {\g}_0^*$ if there is a
nonzero vector $v_{\xi}\in M$ satisfying the following conditions:
\begin{enumerate}
\renewcommand{\labelenumi}{(\roman{enumi})}
\item $hv_{\xi}=\xi(h)v_{\xi}$, for all $h\in \g_0$,

\item $\g_+ v_\xi = 0\/$,

\item ${\mc U} ( \g_- ) v_\xi = M$.
\end{enumerate}
Then $v_\xi$ is called a {\em highest weight vector} of $M$.
Note that by declaring the highest weight vector of the highest
weight module $M$ to be of degree zero, the module $M$ is
naturally $\half\Z$-graded. More precisely, we have $
M=\oplus_{r\in\half\Z_+}M_{-r} $ and $M_0=\C v_\xi$.
A homogeneous nonzero vector $v$ in the highest weight module $M$
is said to be {\em singular} if $\g_+v=0$. A highest weight module
is irreducible if and only if the space of singular vectors is
$1$-dimensional.

Associated with every $\xi \in {\g}_0^*$, there is a {\em Verma module} $$M(\g,
\xi):={\mc U}( \g)\otimes_{{\mc U}( \g_0\oplus\g_+)}\C v_\xi,$$
where $\C v_\xi$ is regarded as a ${\mc U}(\g_0\oplus\g_+)$-module
by setting $hv_\xi=\xi(h)v_\xi$ for all $h\in \g_0$ and
$\g_+v_\xi=0$. Note that $M(\g, \xi)$ is a highest weight module
and for every highest weight module $M$ of highest weight $\xi$,
there is a natural epimorphism $\varphi$ from $M(\g, \xi)$ onto
$M$ determined by $\varphi(v_\xi)=u_\xi$, where $u_\xi$ is a
highest weight vector in $M$. Thus every highest weight module $M$
of highest weight $\xi$ is a quotient of $M ( \g, \xi)$. Moreover,
$M(\g, \xi)$ contains a unique maximal proper submodule which is
also graded. Hence, for any $\xi \in {\g}_0^*$, there is the
unique irreducible highest weight module, denoted by $L ( \g, \xi
)$, which is isomorphic to the quotient of $M ( \g, \xi)$ by the
maximal proper graded submodule.

We recall the following criterion for quasi-finite highest weight
modules.

\begin{prop}\label{quasi-finite}\cite{CW3}
Let $\g=\oplus_{j\in\hf \Z}\g_j$ be a $\hf \Z$-graded Lie
superalgebra such that $\g_0$ is abelian.
Let $M = \oplus_{j \in \hf \Z_+} M_{-j} $ be a highest weight $\g$-module
with highest weight $\xi\in\g_0^*$.  For any non-zero highest weight
vector $v_\xi$ in $M$, the subspace $\g_jv_\xi$ is finite-dimensional for
all $j$ if and only if $M$ is quasi-finite.
\end{prop}

 Let $M$ be a $\g$-module. For any $\la\in \g_0^*$, set
\begin{equation*}
M_\la=\{v\in M\ \mid\ hv=\la (h)v,\ \hbox {for all}\  h \in
\g_0\}.
\end{equation*}
When $M_\la\not= 0$, $\la$ is called a {\em weight} of $M$, and
$M_\la$ is called the {\rm weight space} of weight $\la$. We let
$P(M)$ denote the set of all weights of $M$. A graded $\g$-module
$M = \oplus_{j \in \hf \Z} M_j $ is called {\em
$\g_0$-diagonalizable} if $M$ satisfies the following conditions:
\begin{enumerate}
\renewcommand{\labelenumi}{(\roman{enumi})}
\item $M_\la$ is finite dimensional,

\item $M_j=\oplus_{\la\in P(M)}(M_\la\cap M_j)$, for all $j \in
\hf \Z$.
\end{enumerate}

For any $\la\in \g_0^*$, we  also set
\begin{equation*}
\g_\la=\{ x\in {\g}\ \mid\ [h,x]=\la (h)x,\ \hbox {for all}\  h
\in \g_0\}.
\end{equation*}

As all the infinite rank Lie superalgebras in
\secref{sect superalgebras} are $\hf\Z$-graded, the
representation theoretical notions discussed above are all valid
for them.  The following easy lemma is also useful for the purpose
of studying their representation theory.

\begin{lem}\label{pi}
Let $\pi$ be any transcendental real number over the field of
rational numbers. For any integers $j_1<j_2<\cdots<j_n$, we let
$v_i:=(\pi^{ij_1},\pi^{ij_2},\cdots,\pi^{ij_n})$, for
$i=1,2,\cdots,n$. Then $v_1,v_2,\cdots,v_n$ are linearly
independent in $\C^n$.
\end{lem}
\begin{proof}
Let
\begin{equation*}
f(x)={\rm det}\begin{pmatrix}
x^{j_1}&x^{j_2}&\cdots,&x^{j_n}\\
x^{2j_1}&x^{2j_2}&\cdots,&x^{2j_n}\\
\vdots&\vdots&\cdots&\vdots \\
x^{nj_1}&x^{nj_2}&\cdots,&x^{nj_n}
\end{pmatrix}.
\end{equation*}
Then $f(x)$ is a nonzero Laurent polynomial with integral
coefficients. Therefore $f(\pi)\not= 0$ and this implies the
lemma.
\end{proof}

\subsection{Quasi-finite $\hgltwo$-modules}
For any $k\in \hf \Z$ and $N\in\hf\Z_+$, we let
\begin{equation*}
(\hgltwo)_{k,N}:=\{ x\in \hgltwo\ \mid \ x=\sum_{|j|\ge N, \atop
j\in \hf \Z }a_je_{j-k,j}, \ a_j\in \C \ \}.
\end{equation*}
The following lemma can be confirmed by a straightforward computation.

\begin{lem}\label{growth}
Given any fixed positive integer or half integer $N$, we have
\begin{equation*}
[(\hgltwo)_{p}, (\hgltwo)_{-k, k+N}] \subset
(\hgltwo)_{-(k-p),(k-p)+N},
\end{equation*}
for all $k,p \in \hf \Z_+ $ with $p\le k$.
\end{lem}

\begin{prop}\label{Nvanish}
Let $M = \oplus_{j \in \hf \Z_+} M_{-j} $ be a highest weight
$\hgltwo$-module and $v_0$ a non-zero highest weight vector.
If $(\hgltwo)_{-r}v_0$ is finite dimensional for a fixed number
$r\in \hf\N$, then for every $p \in \hf \Z_+ $ with $p < r$,
there exists $N\in \N$ such that
\begin{equation*}
(\hgltwo)_{-p,N}v_0=0.
\end{equation*}
In particular, $(\hgltwo)_{-p}v_0$ is finite dimensional for all
$p\le r$.
\end{prop}

\begin{proof}
Fixing a transcendental real number $\pi$, we let
$w_i=\sum_{j\in \hf \Z }\pi^{2ij}e_{j+r,j}$, for each $i\in \N$, which belong to
$(\hgltwo)_{-r}$.  For any
$x=\sum_{l=1}^k\alpha_{i_l}w_{i_l}$, where
$\alpha_{i_1},\alpha_{i_2},\cdots,\alpha_{i_k}$ are nonzero
complex numbers with $i_1<i_2<\cdots<i_k$, we can re-write
it as $x=\sum_{j\in \hf \Z }\beta_{j}e_{j+r,j}$.
Then by applying \lemref{pi} we easily show that
the $\beta_j$ are nonzero except for finitely many $j$. Thus there always exists
some positive integer $N$ with $N>r$ such that $\beta_j\not=0$ for all
$j$ with $|j|\ge N$. Since $(\hgltwo)_{-r}v_0$ is finite dimensional,
we can always find nonzero complex
numbers $\alpha_{i_1},\alpha_{i_2},\cdots,\alpha_{i_k}$ so that
$x=\sum_{l=1}^k\alpha_{i_l}w_{i_l}$ satisfies $x v_0=0$. We fix
such an $x$.

We shall prove the proposition by contradiction. Assume that there exists $p \in
\hf \Z_+ $ with $p < r$ such that $(\hgltwo)_{-p,q}v_0\not= 0$,
for all $q\in \N$. Then we can find $y:=\sum_{|j|\ge N
}a_je_{j+p,j}\in (\hgltwo)_{-p,N}$ such that $y v_0\not=0$. We
claim that corresponding to each such $y$,
there exits a $u=\sum_{j\in \hf \Z  }b_j e_{j+p,j+r}\in (\hgltwo)_{r-p}$
such that
\begin{equation}
[u,x]=y.\label{contradict}
\end{equation}
Indeed if we choose an element $u$ with the coefficients $b_j$, $-N-2r< j< N+r$, given by
\begin{equation*}
b_j:=\left \{ \begin{array}{ll} 0,& \mbox{if \ \ $-N-r<j<N$;} \\
\\
\frac{a_j}{\beta_j}, & \mbox{if \ \ $N\le j < N+r$;}\\
\\
\frac{-(-1)^{4r(r-p)}a_{j+r}}{\beta_{j+p}}, & \mbox{if \ \
$-N-2r< j \le -N-r$};\\
\end{array}
\right.
\end{equation*}
and the $b_j$ for $j\le -N-2r$ or $j \ge N+r$ given recursively by
\begin{equation*}
b_j:=\left \{ \begin{array}{ll}
{{a_j+(-1)^{4r(r-p)}\beta_{j+p-r}b_{j-r}}\over \beta_j},
& \mbox{if \ \  $j \ge N+r$;}\\
\\
{{(-1)^{4r(r-p)}(b_{j+r}\beta_{j+r}-a_{j+r})}\over {\beta_{j+p}}},
& \mbox{if \ \ $j\le -N-2r$},
\end{array}
\right.
 \end{equation*}
then \eqref{contradict} holds true as can be shown by a direct computation.
However, equation \eqref{contradict} leads to the obvious contradiction
$y v_0=[u, x]v_0=0$. This completes the proof.
\end{proof}

The following theorem is an obvious consequence of
\propref{quasi-finite} and \propref{Nvanish}.

\begin{thm}\label{criterion-gl-N}
Let $M = \oplus_{j \in \hf \Z_+} M_{-j} $ be a highest weight
$\hgltwo$-module with highest weight $\xi$ and $v_0$ a
non-zero highest weight vector. Then M is quasi-finite if and
only if for every $r \in \hf \Z_+ $, there exists $N\in \N$ such
that
\begin{equation*}
(\hgltwo)_{-r,N}v_0=0.
\end{equation*}
In this case,  $M={\mc U}(\hgltwo^f)v_0$
and hence is $(\hgltwo)_0$-diagonalizable.
\end{thm}

Denote by $\Lambda_0, \omega_s$ ($s\in\hf\Z$), the fundamental
weights of $\hgltwo$, where $\Lambda_0, \omega_s\in(\hgltwo)_0^*$,
are defined by
\begin{align*}
\omega_s(\sum_{r\in\hf\Z}a_re_{rr}+dC)=a_s, \quad
\Lambda_0(\sum_{r\in\hf\Z}a_re_{rr}+dC)=d,
\end{align*}
for all $\sum_{r\in\hf\Z}a_re_{rr}+dC \in (\hgltwo)_0.$

\begin{thm}\label{criterion-gl}
Let $M$ be an irreducible highest weight $\hgltwo$-module and
$v_\xi$ a non-zero highest weight vector. Then the following
are equivalent:
\begin{enumerate}
\renewcommand{\labelenumi}{(\roman{enumi})}
\item $M$ is quasi-finite,

\item $(\hgltwo)_{-\hf}v_\xi$ is finite dimensional,

\item there exists $N\in \N$ such that
$\xi=\sum_{|j|\le N, \atop j\in\hf\Z} \xi_j\omega_j+d\La_0$,
where $\xi_j$, $d\in \C$.
\end{enumerate}
\end{thm}

\begin{proof}
It clearly follows from \propref{Nvanish} that (i) implies (ii) and
(ii) implies (iii). Now we show that (iii) implies
(i). Assume that $\xi=\sum_{|j|\le N_0, \atop j\in\hf\Z}
\xi_j\omega_j+d\La_0$. Then by \thmref{criterion-gl-N}, it is
sufficient to show that for all $r\in \hf\Z_+$, there exists $N\in
\N$ such that $(\hgltwo)_{-r,N}v_\xi= 0$. We shall prove it by
induction. It is obviously true for $r=0$ and we assume that it
is also true for all $p\in\hf\N$ with $0\le p<r$.  Choose $N_p\in \N$ such that
$(\hgltwo)_{-p,N_p}v_\xi=0$. Let
$N={\max}\{N_0+1,N_{1/2},\cdots,N_{r-1/2}\}$. For all $p\in\hf\N$
with $p> r$, it is clear that
$(\hgltwo)_{p}(\hgltwo)_{-r,r+N}v_\xi \subseteq
(\hgltwo)_{p-r}v_\xi=0$. By \lemref{growth}, we also have
\begin{align*}
(\hgltwo)_{p} (\hgltwo)_{-r,r+N}v_\xi& \subseteq [(\hgltwo)_{p},
(\hgltwo)_{-r,r+N}]v_\xi \\
&\subseteq (\hgltwo)_{-(r-p),(r-p)+N}v_\xi \\
&=0,
\end{align*}
for all $p\in\hf\N$ with $0< p\le r$. Thus
$(\hgltwo)_+(\hgltwo)_{-r,r+N}v_\xi =0$. Similarly using
\lemref{growth} again,
\begin{align*}
&(\hgltwo)_{0} (\hgltwo)_{-r,r+N}v_\xi\\
\subseteq & [(\hgltwo)_{0},
(\hgltwo)_{-r,r+N}]v_\xi+ (\hgltwo)_{-r,r+N}(\hgltwo)_{0}v_\xi \\
\subseteq & (\hgltwo)_{-r,r+N}v_\xi. \\
\end{align*}
Therefore ${\mc U}((\hgltwo)_-)(\hgltwo)_{-r,r+N}v_\xi$ is a
proper submodule of $M$. Thus we have $(\hgltwo)_{-r,r+N}v_\xi=0$
since $M$ is irreducible.
\end{proof}

We will let $\hglone$ denote the the $\Z$-graded subalgebra of
$\hgltwo$ defined by
\begin{equation*}
\hglone:=\{(a_{ij})\in \hgltwo | \ a_{ij}=0 \ \mbox{for $i=\hf$ or
$j=\hf$}\ \}\oplus\C C.
\end{equation*}
Then $\hglone$ is a Lie algebra with the natural $\Z$-gradation
induced from $\hgltwo$. Therefore, the notions of highest weight
modules, quasi-finite highest weight modules, etc. can also be
defined for the Lie algebra $\hglone$. We also let $\Lambda_0$,  $\omega_i$,
$i\in\Z$, denote the fundamental weights of $\hglone$. That is,
$\Lambda_0, \omega_i\in(\hglone)_0^*$ with
$\omega_i(\sum_{j\in\Z}a_je_{jj}+dC)=a_i$ and
$\Lambda_0(\sum_{j\in\Z}a_je_{jj}+dC)=d$, for all
$\sum_{j\in\Z}a_je_{jj}+dC \in \g_0$.

Using arguments analogous to those in the proof of \thmref{criterion-gl}
we can prove the following theorem.

\begin{thm}
Let $M$ be an irreducible highest weight $\hglone$-module and
$v_\xi$ a  non-zero highest weight vector. Then the following
are equivalent:
\begin{enumerate}
\renewcommand{\labelenumi}{(\roman{enumi})}
\item $M$ is quasi-finite,

\item $(\hglone)_{-1}v_\xi$ is finite dimensional,

\item there exists $N\in \N$ such that $\xi=\sum_{|j|\le N}
\xi_j\omega_j+d\La_0$, where $\xi_j$, $d\in \C$.
\end{enumerate}
\end{thm}

\subsection{Quasi-finite $\hA$-modules}
Results proved in the last subsection all generalize to the Lie
superalgebra $\hA$. We shall summarize them here, but omit all the proofs,
as they are the same as in the case of $\hgltwo$.

For any $k\in \hf \Z$ and $N\in \hf \Z_+$, we let
\begin{equation*}
\hA_{k,N}:=\{ x\in \hA\ \mid \ x=\sum_{|j|\ge N, \atop j\in \hf
\Z^* }a_je_{j-k,j}, \ a_j\in \C \ \}.
\end{equation*}

\begin{prop}\label{hA-Nvanish}
Let $M = \oplus_{j \in \hf \Z_+} M_{-j} $ be a highest weight
$\hA$-module and $v_0$ be a non-zero highest weight vector in $M$. If
$\hA_{-r}v_0$ is finite dimensional for a fixed number $r\in
\hf\N$, then for every $p \in \hf \Z_+ $ with $p < r$,  there
exists $N\in \N$ such that
\begin{equation*}
\hA_{-p,N}v_0=0.
\end{equation*}
In particular, $\hA_{-p}v_0$ is finite dimensional for all $p\le
r$.
\end{prop}

The following theorem is an obvious consequence of
\propref{quasi-finite} and \propref{hA-Nvanish}.

\begin{thm}\label{criterion-hA-N}
Let $M = \oplus_{j \in \hf \Z_+} M_{-j} $ be a highest weight
$\hA$-module and $v_0$ be a non-zero highest weight vector in $M$. Then
M is quasi-finite if and only if for every $r \in \hf \Z_+ $, there
exists $N\in \N$ such that
\begin{equation*} \hA_{-r,N}v_0=0. \end{equation*}
In this case, $M={\mc U}(\hA^f)v_0$ and is $\hA_0$-diagonalizable.
\end{thm}

We let $\Lambda_0$, $\omega_s\in \hA_0^*$ ($s\in\hf\Z^*$,
$\Z^*:=\Z\backslash\{0\}$),
denote the fundamental weights of $\hA$, which are  defined
by $\omega_s(\sum_{r\in\hf\Z^*}a_r e_{rr}+dC)=a_s$, \
$\Lambda_0(\sum_{r\in\hf\Z^*}a_re_{rr}+dC)=d$, for all
$\sum_{r\in\hf\Z^*}a_re_{rr}+dC \in \hA_0$.

\begin{thm}\label{criterion-hA}
Let $M$ be an irreducible highest weight $\hA$-module and $v_\xi$
a non-zero highest weight vector in $M$. Then the following are
equivalent:
\begin{enumerate}
\renewcommand{\labelenumi}{(\roman{enumi})}
\item $M$ is quasi-finite,

\item $\hA_{-\hf}v_\xi$ is finite dimensional,

\item there exists $N\in \N$ such that
$\xi=\sum_{|j|\le N, \atop j\in\hf\Z^*} \xi_j\omega_j+d\La_0$,
where $\xi_j$, $d\in \C$.
\end{enumerate}
\end{thm}

\subsection{Quasi-finite $\hC$-modules}
For any $k\in \hf \Z$ and $N\in \hf\Z_+$, we let
\begin{equation*}
\hC_{k,N}:=\{ x\in \hC\ \mid \ x=\sum_{j\ge N, \atop j\in \hf\N
}a_j\te_{j-k,j}, \ a_j\in \C \ \}.
\end{equation*}

We have the following lemma, which can be confirmed by a direct
computation.

\begin{lem}\label{growth-hC}
Given any fixed positive integer or half integer $N$, we have
\begin{equation*}
[\hC_{p}, \hC_{-k,k+N}] \subset \hC_{-(k-p),(k-p)+N},
\end{equation*}
for all $k,p \in \hf \Z_+ $ with $p\le k$.
\end{lem}

\begin{prop}\label{Nvanish-hC}
Let $M = \oplus_{j \in \hf \Z_+} M_{-j} $ be a highest weight
$\hC$-module and $v_0$ a non-zero highest weight vector in $M$. If
$\hC_{-r}v_0$ is finite dimensional for a fixed number $r\in
\hf\N$, then for every $p \in \hf \Z_+ $ with $p < r$, there
exists $N\in \N$ such that
\begin{equation*}
\hC_{-p,N}v_0=0.
\end{equation*}
In particular, $\hC_{-p}v_0$ is finite dimensional for all $p\le
r$.
\end{prop}

\begin{proof}
The proof is quite similar to \propref{Nvanish}. Choosing a fixed
transcendental real number $\pi$, we let $w_i=\sum_{j\in \hf \N
}\pi^{2ij}\te_{j+r,j}$ be an element in $\hC_{-r}$ for each $i\in
\N$. We put $x=\sum_{l=1}^k\alpha_{i_l}w_{i_l}$, where
$\alpha_{i_1},\alpha_{i_2},\cdots,\alpha_{i_k}$ are nonzero
complex numbers. By \lemref{pi}, we have $x=\sum_{j\in \hf \N
}\beta_{j}\te_{j+r,j}$ such that $\beta_j$ are nonzero except
finitely many $j$. Choose a positive integer $N$ with $N>r$ such
that $\beta_j\not=0$ for all $j$ with $j\ge N$.

We shall prove by contradiction. Assume that there exists $p \in
\hf \Z_+ $ with $p < r$ such that $\hC_{-p,q}v_0\not= 0$, for all
$q\in \N$. Therefore, we can find $y:=\sum_{j\ge N
}a_j\te_{j+p,j}\in \hC_{-p,N}$ such that $yv_0\not=0$. Let
$u=\sum_{j\in \hf \N }b_j\te_{j+p,j+r}\in \hC_{r-p}$, where $b_j$
is defined as follows:
 \begin{equation*}
b_j:=\left \{ \begin{array}{ll} 0,&
\mbox{if \ \ $0<j<N$;} \\
\\
{a_j\over \beta_j}, & \mbox{if \ \ $N\le j < N+r$;}\\
\end{array}
\right.
 \end{equation*}
and for $j\ge N+r$, $b_j$ is defined by the following recursive
relations:
\begin{equation*}
b_j:= {{a_j+(-1)^{4r(r-p)}\beta_{j+p-r}b_{j-r}}\over \beta_j}.
 \end{equation*}
Direct computations show that $[u,x]=y$,  and hence
$uxv_0=yv_0\not=0$. On the other hand, we can find nonzero complex
numbers $\alpha_{i_1},\alpha_{i_2},\cdots,\alpha_{i_k}$ so that
$xv_0=\sum_{l=1}^k\alpha_{i_l}w_{i_l}v_0=0$ since $\hC_{-r}v_0$ is
finite dimensional, which contradicts $uxv_0\not=0$.
\end{proof}

The following theorem is an immediate consequence of the
\propref{quasi-finite} and \propref{Nvanish-hC}.

\begin{thm}\label{criterion-hC-N}
Let $M = \oplus_{j \in \hf \Z_+} M_{-j} $ be a highest weight
$\hC$-module with highest weight $\xi$ and $v_0$ a non-zero
highest weight vector in $M$. Then M is quasi-finite if and only if for
every $r \in \hf \Z_+ $, there exists $N\in \N$ such that
\begin{equation*}
\hC_{-r,N}v_0=0.
\end{equation*}
In this case, $M={\mc U}(\hC^f)v_0$ and is $\hC_0$-diagonalizable.
\end{thm}

Let $\Lambda_0$, $\omega_s$, $s\in\hf\N$, denote the fundamental weights of
$\hC$, that is, $\Lambda_0$, $\omega_s\in(\hC_0)^*$ defined by
$\omega_s(\sum_{r\in\hf\N}a_r\te_{rr}+dC)=a_s$,\
$\Lambda_0(\sum_{r\in\hf\N}a_r\te_{rr}+dC)=d$, for all
$\sum_{r\in\hf\N}a_r\te_{rr}+dC \in \hC_0$. We have the following
theorem, the proof of which will be omitted here, since it is
similar to the proof of \thmref{criterion-gl}.

\begin{thm}\label{criterion-hC}
Let $M$ be an irreducible highest weight $\hC$-module with a non-zero
highest weight vector $v_\xi$. Then the following are
equivalent:
\begin{enumerate}
\renewcommand{\labelenumi}{(\roman{enumi})}
\item $M$ is quasi-finite,

\item $\hC_{-\hf}v_\xi$ is finite dimensional,

\item there exists $N\in \N$ such that $\xi=\sum_{j\le N, \atop
j\in\hf\N} \xi_j\omega_j+d\La_0$, where $\xi_j$, $d\in \C$.
\end{enumerate}
\end{thm}

\subsection{Quasi-finite $\hD$-modules}
All theorems and propositions proved in the last subsection can be adapted to
the Lie superalgebra $\hD$. We summarize the results here, but omit their proofs
as they are much the same as in the case of $\hC$.

For any $k\in \hf \Z$ and $N\in\hf \Z_+$, we let
\begin{equation*}
\hD_{k,N}:=\{ x\in \hD\ \mid \ x=\sum_{j\ge N, \atop j\in \hf\N
}a_j\te_{j-k,j}, \ a_j\in \C \ \}.
\end{equation*}
The following lemma can be proven by a direct computation.

\begin{lem}\label{growth-hD}
Given any fixed positive integer or half integer $N$, we have
\begin{equation*}
[\hD_{p}, \hD_{-k,k+N}] \subset \hD_{-(k-p),(k-p)+N},
\end{equation*}
for all $k,p \in \hf \Z_+ $ with $p\le k$.
\end{lem}

\begin{prop}\label{Nvanish-hD}
Let $M = \oplus_{j \in \hf \Z_+} M_{-j} $ be a highest weight
$\hD$-module and $v_0$ a non-zero highest weight vector in $M$. If
$\hD_{-r}v_0$ is finite dimensional for a fixed number $r\in
\hf\N$, then for every $p \in \hf \Z_+ $ with $p < r$, there
exists $N\in \N$ such that
\begin{equation*}
\hD_{-p,N}v_0=0.
\end{equation*}
In particular, $\hD_{-p}v_0$ is finite dimensional for all $p\le
r$.
\end{prop}

The following theorem is an obvious consequence of the
\propref{quasi-finite} and \propref{Nvanish-hD}.

\begin{thm}\label{criterion-hD-N}
Let $M = \oplus_{j \in \hf \Z_+} M_{-j} $ be a highest weight
$\hD$-module with highest weight $\xi$ and $v_0$ a non-zero
highest weight vector. Then M is quasi-finite if and only if for
every $r \in \hf \Z_+ $, there exists $N\in \N$ such that
\begin{equation*}
\hD_{-r,N}v_0=0.
\end{equation*}
In this case, $M={\mc U}(\hD^f) v_0$ and  is $\hD_0$-diagonalizable.
\end{thm}

Let $\Lambda_0$, $\omega_s$, $s\in\hf\N$, denote the fundamental weights of
$\hD$, that is, $\Lambda_0$, $\omega_s\in(\hD_0)^*$ defined by
$\omega_s(\sum_{r\in\hf\N}a_r\te_{rr}+dC)=a_s$,\
$\Lambda_0(\sum_{r\in\hf\N}a_r\te_{rr}+dC)=d$, for all
$\sum_{r\in\hf\N}a_r\te_{rr}+dC \in \hD_0$.

\begin{thm}\label{criterion-hD}
Let $M$ be an irreducible highest weight $\hD$-module and $v_\xi$
a non-zero highest weight vector in $M$. Then the following are
equivalent:
\begin{enumerate}
\renewcommand{\labelenumi}{(\roman{enumi})}
\item $M$ is quasi-finite,

\item $\hD_{-\hf}v_\xi$ is finite dimensional,

\item there exists $N\in \N$ such that $\xi=\sum_{j\le N, \atop
j\in\hf\N} \xi_j\omega_j+d\La_0$, where $\xi_j$, $d\in \C$.
\end{enumerate}
\end{thm}

\begin{rem}
All results in this section can be restricted to the non-super case, leading to
descriptions and classifications of the irreducible quasi-finite highest weight
modules over $\hglone$ and its $\Z$-graded subalgebras.
\end{rem}

\section{Unitarizable representations and their free field realizations}
\label{sect unitarizable}

In this section we study in detail a particularly nice class of
modules over the infinite rank Lie  superalgebras
$\hgltwo$, $\hA$, $\hC$ and $\hD$, namely,  the
quasi-finite irreducible highest weight modules, which are
unitarizable with respect to natural choices of $\ast$-structures
for these Lie superalgebras. Several results are obtained here,
including (1). the classification of the unitarizable quasi-finite
irreducible highest weight modules over the Lie superalgebras
$\hgltwo$, $\hA$, $\hC$ and $\hD$; (2). realizations of these
irreducible modules on Fock spaces; and (3).
generalized Howe dualities between these infinite rank Lie
superalgebras and classical Lie algebras.   The generalized Howe
dualities will provide the principal tool for constructing
character formulae for the unitarizable quasi-finite irreducible
highest weight modules in the next section.

\subsection{Unitarizable $\hgltwo$-modules and their Fock
space realizations}\label{unitarizable-gl}
\subsubsection{Free field realization of  $\hgltwo$ and $(\hgltwo, gl_d)$-duality}
Let $gl_d$ denote the Lie algebra of all complex $d\times d$
matrices. Let $\{ e^1,\ldots,e^d\}$ be the standard basis for
$\C^d$. Denote by $E_{ij}$ the elementary matrix with $1$ in the
$i$-th row and $j$-th column and $0$ elsewhere. Then
$\fh_d=\sum_{i=1}^d\C E_{ii}$ is a Cartan subalgebra, while
$\fb_d=\sum_{1\le i\le j\le d}\C E_{ij}$ is the standard Borel
subalgebra containing $\h_d$. The weight of $e^i$ is denoted by
$\tilde{\epsilon}_i$ for $1\le i\le d$.

We regard a finite sequence $\la=(\la_1,\cdots,\la_d )$ of complex
numbers as an element of the dual vector space $\fh^*_d$ of
$\fh_d$ defined by $\la(E_{ii})=\la_i$, for $i=1,\cdots,d$. Denote
by $V^\la_d$ the irreducible $gl_d$-module with highest weight
$\la$ relative to the standard Borel subalgebra $\fb_d$.

Consider $d$ pairs of free fermionic fields $\psi^{\pm,i}(z)$ and
$d$ pairs of free symplectic bosonic fields $\gamma^{\pm,i}(z)$,
$i=1,\cdots,d$,  with the following mode expansions:
\begin{align*}
\psi^{+,i}(z)&=\sum_{n\in\Z}\psi^{+,i}_n z^{-n-1},\quad\quad\
\psi^{-,i}(z)=\sum_{n\in\Z}\psi^{-,i}_nz^{-n},\\
\gamma^{+,i}(z)&=\sum_{r\in\frac{1}{2}+\Z}\gamma^{+,i}_rz^{-r-1/2},\quad
\gamma^{-,i}(z)=\sum_{r\in\frac{1}{2}+\Z}\gamma^{-,i}_rz^{-r-1/2},
\end{align*}
where the operators $\psi^{+,i}_n$ and $\gamma^{+,i}_r$ satisfy
the usual (anti-)commutation relations with the non-trivial ones
being given by
\begin{align*}
\psi^{+,i}_m \psi^{-,j}_n + \psi^{-,j}_n \psi^{+,i}_m&=\delta_{ij}\delta_{m+n,0}, \\
\gamma^{+,i}_r \gamma^{-,j}_s - \gamma^{-,j}_s \gamma^{+,i}_r
&=\delta_{ij}\delta_{r+s,0}.
\end{align*}
Denote by $\ccA$ the associative superalgebra generated by these
operators. $\ccA$ admits a $\ast$-structure $\omega:
\ccA\rightarrow \ccA$ defined by
\begin{align}
\omega(\psi_m^{+,k})=\psi_{-m}^{-,k},&\quad
\omega(\psi_m^{-,k})=\psi_{-m}^{+,k}, \label{conjugate-psi}\\
\omega(\gamma_r^{+,k})=\left\{%
\begin{array}{ll}
    \gamma_{-r}^{-,k}, &\text{if } r>0, \\
    -\gamma_{-r}^{-,k}, & \text{if } r<0. \\
\end{array}%
\right.,&\quad
 \omega(\gamma_r^{-,k})=\left\{%
\begin{array}{ll}
    -\gamma_{-r}^{+,k}, &\text{if } r>0, \\
    \gamma_{-r}^{+,k}, & \text{if } r<0, \\
\end{array}%
\right.\label{conjugate-gamma}
\end{align}
for all $m\in \Z$, $r\in \frac{1}{2}+\Z$, $k=1,2,\cdots, d$.  It
can be easily shown that $\omega$ indeed defines an anti-linear
anti-involution on $\ccA$.

We shall take the operators $\psi^{+,i}_n$, $\psi^{-,i}_m$,
$\gamma^{\pm,i}_r$, with $i=1,2\cdots,d$, $n\ge 0$, $m>0$, $r>0$,
as annihilation operators, and the rest as creation operators. Let
$\F$ be the Fock space of $\ccA$ generated by the vacuum vector
$\vac$, which is nullified by the annihilation operators, i.e.,
\begin{align*}
&\psi^{+,i}_n\vac=\psi^{-,i}_m\vac=\gamma^{\pm,i}_r\vac=0, \\
&\rm{for \ all} \  i=1,2, \cdots,d;  \ \ n\ge 0; \ \ m>0;   \ \ r>0.
\end{align*}
Let $\langle\cdot|\cdot\rangle$ be the contravariant Hermitian from
on the Fock space $\F$ defined with respect to the anti-linear
anti-involution $\omega$ given in \eqref{conjugate-psi} and
\eqref{conjugate-gamma}. We normalize the form on the vacuum vector $\vac$ so that
$\langle 0 | 0 \rangle=1$.
\begin{lem}
The Fock space $\F$ equipped with the contravariant Hermitian from
$\langle\cdot|\cdot\rangle$ is a unitarizable $\ccA$-module.
\end{lem}
\begin{proof} In the `particle number basis' for the Fock space, the lemma can be
confirmed by a straightforward calculation, details of which
are omitted here. However, see Remark \ref{remark-symplectic} below.
\end{proof}
\begin{rem} \label{remark-symplectic} If we remove all the minus signs from
\eqref{conjugate-gamma},
we still obtain a $\ast$-structure for $\ccA$. In fact this is the
usual conjugation rule for symplectic bosons adopted in the
physics literature. However, it is quite well known that the Fock
space of symplectic bosons is not unitarizable with respect to the
usual conjugation rule.
\end{rem}

The Lie superalgebra $\hgltwo$ can be realized on the Fock
space $\F$ as follows:
 \begin{align}
C&=d, \nonumber \\
e_{ij}&=\sum_{p=1}^d:\psi^{+,p}_{-i}\psi^{-,p}_{j}:, &
e_{rs}&=-\sum_{p=1}^d:\gamma^{+,p}_{-r}\gamma^{-,p}_{s}:, \label{realization}\\
e_{is}&=\sum_{p=1}^d:\psi^{+,p}_{-i}\gamma^{-,p}_{s}:, &
e_{rj}&=-\sum_{p=1}^d:\gamma^{+,p}_{-r}\psi^{-,p}_{j}:, \nonumber
\end{align}
where $i,j\in\Z$ and $r,s\in\frac{1}{2}+\Z$. There also exists an
action of $gl_d$ on $\F$, which is given by the formula
\begin{equation*}
E_{ij}=\sum_{n\in\Z}:\psi^{+,i}_{-n}\psi^{-,j}_{n}:-
\sum_{r\in1/2+\Z}:\gamma^{+,i}_{-r}\gamma^{-,j}_{r}:.
\end{equation*}
Here and further $::$ denotes the normal ordering of operators.
That is, if $A$ and $B$ are two operators, then $:AB:=AB$, if $B$
is an annihilation operator, while $:AB:=(-1)^{|A||B|}BA$,
otherwise, where $|X|$ denotes the parity of $X$.

The following result is the $(\hgltwo, \ gl_d)$-duality of
\cite{CW3}.
\begin{thm} \cite{CW3}\label{duality}
The Lie superalgebra $\hgltwo$ and $gl_d$ form a dual pair on $\F$
in the sense of Howe.  Furthermore we have the following
(multiplicity-free) decomposition of $\F$ with respect to their
joint action
\begin{equation*}
\F\cong\sum_{\la}L(\hgltwo,\Lambda(\la))\otimes V_d^{\la},
\end{equation*}
where the summation is over all generalized partitions of length
$d$.
\end{thm}
Here the notation $\Lambda(\la)$ requires some explanation. For a
generalized partition $\la=(\la_1,\la_2,\cdots,\la_d)$ of length
$d$, we have the shifted Frobenius notation (see
\secref{Frobenius}) $$\La(\la)=(\xiall).$$
We identify $\La(\la)$ with an element of
the dual space $(\hgltwo)_0^*$ of $(\hgltwo)_0$
defined by
\begin{equation}
\La(\la):=\sum_{{{s+\hf}\le j\le r}\atop
j\in\hf\Z}\xi_j\omega_j+d\La_0.
\end{equation}

\subsubsection{Unitarizable $\hgltwo$-modules}
Let us consider the anti-linear anti-involution on $\hgltwo^f$
defined by
\begin{equation*}
C\mapsto C, \quad e_{p q} \mapsto (-1)^{[p]+[q]} e_{q p}, \ \ \forall
p, q,
\end{equation*}
where $[r]=1$ if $r$ is a negative half integer, and $[r]=0$ otherwise.
This map naturally
extends to an anti-linear anti-involution on $\hgltwo$ with
\begin{equation*}
C\mapsto C, \quad \sum_{p\in \frac{1}{2}\Z}a_pe_{p-k,p} \mapsto\sum_{p\in
\frac{1}{2}\Z} (-1)^{[p]+[p-k]}\overline{a}_pe_{p,p-k},
\end{equation*}
for all $\sum_{p\in \hf\Z}a_pe_{p-k,p}\in (\hgltwo)_k $ and for
all $k\in \hf\Z$. Here $\overline{a}$ denotes the
complex conjugate of the complex number $a$. Abusing the notation,
we shall also denote this map by $\omega$.

The realization \eqref{realization} of $\hgltwo^f$ in $\ccA$ defines an
associative superalgebra homomorphism $\Phi: {\mc U}(\hgltwo^f)\rightarrow
\ccA$. By using \eqref{conjugate-psi} and \eqref{conjugate-gamma} we
can show by a direct computation that $\omega\Phi(x) = \Phi(\omega(x))$
for all $x\in{\mc U}(\hgltwo^f)$. Thus $\Phi$ is a $\ast$-superalgebra
homomorphism.

Since the Fock space $\F$ equipped with the Hermitian form
$\langle\ \cdot | \cdot\ \rangle$ is a unitarizable
$\ccA$-module, it naturally restricts to a unitarizable module over
$\hgltwo^f$. Hence $\F$ forms a unitarizable $\hgltwo$-module as can be easily
seen by examining the action of $\hgltwo$ on $\F$.
By \thmref{duality}, for each generalized
partition $\lambda$ of length $d$, the irreducible
$\hgltwo$-module $L(\hgltwo,\Lambda(\lambda))$ is unitarizable. We shall show that
such modules exhaust all the unitarizable irreducible quasi-finite highest weight
$\hgltwo$-modules with respect to the $\ast$-structure $\omega$.

\begin{thm}\label{unitary}
Let $M$ be an irreducible quasi-finite highest weight
$\hgltwo$-module with highest weight $\xi$. Then $M$ is
unitarizable with respect to the $\ast$-structure $\omega$
if and only if $\xi=\Lambda(\lambda)$ for some
generalized partition $\lambda$. In other words, $M$ is
unitarizable if and only if
$$\xi=\sum_{{{s+\hf}\le j\le r}\atop
j\in\hf\Z}\xi_j\omega_j+d\La_0$$ such that $d\in \Z_+$,
$-s,r\in\N$ and $\xi_j\in\Z$ for all $j$ satisfying the following
conditions:
\begin{enumerate}
\renewcommand{\labelenumi}{(\roman{enumi})}
\item $\xihalfge \geq 0$, $\xionege \geq 0$,
and $\xi_{r-\hf}=0$ if only if $ r=1$ and $ \xi_\hf=\xi_1=0$,

\item $0\geq\xineghalfge $, $0\geq\xinegzeroge $, and
$\xi_{s+1}=0$ if only if $s=-1$ and $\xi_{-\hf}=\xi_0=0$,

\item $\min\{\xi_\hf,1\}+\xi_1-\xi_0\leq d$.
\end{enumerate}
\end{thm}

\begin{proof} We already know that for any
generalized partition $\lambda$, the $\hgltwo$-module
$L(\hgltwo,\Lambda(\lambda))$ is unitarizable.
Now we show that if $M$ is a unitarizable irreducible quasi-finite highest weight
$\hgltwo$-module with the highest weight $\xi$, then
$\xi=\Lambda(\lambda)$ for some generalized partition $\lambda$.
Let $\langle \ \cdot |  \cdot \  \rangle$ be a positive definite
contravariant Hermitian form on $M$ and $v_\xi$ a highest weight
vector of $M$ such that $ \langle v_\xi | v_\xi\rangle =1$. We put
$\xi(e_{i i})=\xi_i$ for all $i\in \hf\Z$. By
\thmref{criterion-gl}, there exists $N\in \N$ such that
$\xi=\sum_{|j|\le N, \atop j\in\hf\Z} \xi_j\omega_j+d\La_0$, where
$\xi_j$, $d\in \C$. For each $i\in \N$, $\{
e_{i,i}-e_{i+1,i+1},e_{i,i+1},e_{i+1,i}\}$ forms a standard basis
for the Lie algebra $sl(2,\C)$ and $\omega(e_{i,i+1})=e_{i+1,i}$.
Unitarizability of $M$ with respect to this subalgebra requires (see, e.g.,
Theorem 11.7 in \cite{K2})
$\xi_i-\xi_{i+1}=\xi(e_{i,i}-e_{i+1,i+1})\in\Z_+$. Since $M$ is a
quasi-finite highest weight $\hgltwo$-module, we have $\xi_i\in
\Z_+$ for all $i\in \N$ and
\begin{equation*}
\xi_1\ge\xi_2\ge\cdots\ge\xi_n\ge\xi_{n+1}=\xi_{n+2}=\cdots=0\qquad
\mbox{ for some } n\in \N.
\end{equation*}
Similarly, we have $\xi_i\in \Z_+$ for all $i\in \hf+\Z_+$ and
\begin{equation*}
\xi_{\frac{1}{2}}\ge\xi_{\frac{3}{2}}\ge\cdots\ge\xi_{r-\frac{1}{2}}
\ge\xi_{r+\frac{1}{2}}=\xi_{r+\frac{3}{2}}=\cdots=0\qquad \mbox{
for some } r\in \N.
\end{equation*}
Now we are going to show that for each $i\in \hf\N$, either
\begin{equation}\label{>}
\xi_i>\xi_{i+1} \qquad \mbox{or}\qquad \xi_i=\xi_{i+\hf}=0.
\end{equation}
It is
sufficient to show that $\xi_i=\xi_{i+1}$ implies
$\xi_i=\xi_{i+\hf}=0$. Assume that $\xi_i=\xi_{i+1}$. Then
\begin{equation*}
\begin{array}{ll}
\langle e_{i+1,i}v_\xi\mid e_{i+1,i}v_\xi\rangle
 &=\langle v_\xi\mid \omega(e_{i+1,i})e_{i+1,i}v_\xi\rangle\\
 &=\langle v_\xi\mid e_{i,i+1}e_{i+1,i}v_\xi\rangle\\
 &=\xi_i-\xi_{i+1}\\
 &=0.
\end{array}
\end{equation*}
Therefore, we have $e_{i+1,i}v_\xi=0$. On the other hand, we have
\begin{equation*}
\begin{array}{ll}
 \langle e_{i+\hf,i+1}e_{i+1,i}v_\xi\mid e_{i+\hf,i+1}e_{i+1,i}v_\xi\rangle
 &=\langle e_{i+\hf,i}v_\xi\mid e_{i+\hf,i}v_\xi\rangle\\
 &=\langle v_\xi\mid \omega(e_{i+\hf,i})e_{i+\hf,i}v_\xi\rangle\\
 &=\langle v_\xi\mid e_{i,i+\hf}e_{i+\hf,i}v_\xi\rangle\\
 &=\langle v_\xi\mid e_{i,i}+e_{i+\hf,i+\hf}v_\xi\rangle\\
 &=\xi_i+\xi_{i+\hf}.
\end{array}
\end{equation*}
Thus $\xi_i+\xi_{i+\hf}=0$ since $e_{i+1,i}v_\xi=0$. Since
$\xi_i\ge 0$ and $\xi_{i+\hf}\ge 0$, we have
$\xi_i=\xi_{i+\hf}=0$. By \eqnref{>}, there is $r\in \N$ such that
\begin{equation*}
\begin{array}{ll}
\xi_1>\xi_2>\cdots>\xi_r\ge\xi_{r+1}=\xi_{r+2}=\cdots=0,\\
\\
\xi_{\frac{1}{2}}>\xi_{\frac{3}{2}}>\cdots>\xi_{r-\frac{1}{2}}
>\xi_{r+\frac{1}{2}}=\xi_{r+\frac{3}{2}}=\cdots=0,
\end{array}
\end{equation*}
or $ \xi_\hf=\xi_1=0$.

By using similar argument as above (note that
$\omega(e_{i-\hf,i})=-e_{i,i-\hf}$ for all $-i\in \Z_+$), we have
$-\xi_i\in \Z_+$ for all $-i\in \hf\Z_+$, and there is $s\in \N$
such that
\begin{equation*}
\begin{array}{ll}
0=\cdots=\xi_{s-\frac{3}{2}}=\xi_{s-\frac{1}{2}}
\ge\xineghalfge,\\
\\
0=\cdots=\xi_{s-1}=\xi_{s} > \xinegzeroge,
\end{array}
\end{equation*}
or $\xi_{-\hf}=\xi_0=0$.

Now we choose a large positive integer $n$ such that
$\xi(e_{n,n})=\xi(e_{-n,-n})=0$. Consider the subalgebra
$sl(2,\C)$ spanned by $\{
e_{-n,-n}-e_{n,n}+C,e_{-n,n},e_{n,-n}\}$. Note that
$\omega(e_{-n,n})=e_{n,-n}$. Using the standard trick on
unitarizable modules again, we have
$d=\xi(C)=\xi(e_{-n,-n}-e_{n,n}+C)\in\Z_+$. Finally, we need to
show $\min\{\xi_\hf,1\}+\xi_1-\xi_0\leq d$. Since $\langle
e_{1,0}v_\xi\mid e_{1,0}v_\xi\rangle\geq 0$ and
\begin{equation*}
\begin{array}{ll}
 \langle e_{1,0}v_\xi\mid e_{1,0}v_\xi\rangle
 &=\langle v_\xi\mid \omega(e_{1,0})e_{1,0}v_\xi\rangle\\
 &=\langle v_\xi\mid (e_{0,0}-e_{1,1}+C)v_\xi\rangle\\
 &=\xi_0-\xi_{1}+d,
\end{array}
\end{equation*}
we have $\xi_0-\xi_{1}+d\geq 0$. If $\xi_0-\xi_{1}+d > 0$, the
proof the theorem is completed. Otherwise, we have $\xi_0-\xi_1+d=
0$ and $e_{1,0}v_\xi=0$. Therefore we have
\begin{equation}\label{dxi}
d+\xi_0=\xi_1\geq 0.
\end{equation}
On the other hand, by using \eqnref{dxi} and $e_{1,0}v_\xi=0$, we
have
\begin{equation*}
\begin{array}{ll}
0&= \langle e_{\hf,1}e_{1,0}v_\xi\mid
e_{\hf,1}e_{1,0}v_\xi\rangle\\
 &=\langle e_{\hf,0}v_\xi\mid e_{\hf,0}v_\xi\rangle\\
 &=\langle v_\xi\mid \omega(e_{\hf,0})e_{\hf,0}v_\xi\rangle\\
 &=\langle v_\xi\mid (e_{0,0}+e_{\hf,\hf}+C)v_\xi\rangle\\
 &=\xi_0+\xi_{\hf}+d\\
 &=\xi_\hf+\xi_{1}.\\
\end{array}
\end{equation*}
Since $\xi_\hf\geq 0$ and $\xi_1\geq 0$, we have
$\xi_\hf=\xi_{1}=0$ and $\min\{\xi_\hf,1\}+\xi_1-\xi_0=-\xi_0\leq
d$ by \eqnref{dxi} again. This completes the proof of the theorem.
\end{proof}

Recall that in \cite{W1} Wang showed that there is a
$\hglone\times gl_d$ duality on the subspace of the Fock
space $\F$ generated by the fermionic operators.
By modifying the arguments in the proof of \thmref{unitary} we can show
that

\begin{thm}\label{unitary-hglone}
Let $M$ be an irreducible quasi-finite highest weight
$\hglone$-module with the highest weight $\xi$. Then $M$ is
unitarizable if and only if
$$\xi=\sum_{{{s}\le j\le r}}\xi_j\omega_j+d\La_0$$ such that $d\in \Z_+$, $r\in\N$,
$-s\in\Z_+$ and $\xi_j\in\Z$ for all $j$ satisfying the following
conditions:
\begin{enumerate}
\renewcommand{\labelenumi}{(\roman{enumi})}
\item $\xionege \geq 0$,

\item $0 \geq \xi_{s}>\xi_{s+1}>\cdots>\xi_{0} $,

\item $\xi_1-\xi_0\leq d$.
\end{enumerate}
\end{thm}

\subsection{Unitarizable $\hA$-modules and their Fock space
realizations}\label{unitarizable-hA}
\subsubsection{Free field realization of $\hA$ and $(\hA, gl_d)$-duality}
Consider $d$ pairs of free fermions $\widetilde{\psi}^{\pm,i}(z)$
and $d$ pairs of free symplectic bosons $\gamma^{\pm,i}(z)$ (
$i=1,\cdots,d$)
\begin{align*}
\widetilde{\psi}^{+,i}(z)&=\sum_{n\in\Z^{^*}}\psi^{+,i}_nz^{-n-1},\quad\quad\
\widetilde{\psi}^{-,i}(z)=\sum_{n\in\Z^{^*}}\psi^{-,i}_nz^{-n},\\
\gamma^{+,i}(z)&=\sum_{r\in\frac{1}{2}+\Z}\gamma^{+,i}_rz^{-r-1/2},\quad
\gamma^{-,i}(z)=\sum_{r\in\frac{1}{2}+\Z}\gamma^{-,i}_rz^{-r-1/2}
\end{align*}
with the non-trivial anti-commutation relations
$[\psi^{+,i}_m,\psi^{-,j}_n]=\delta_{ij}\delta_{m+n,0}$ and
commutation relations
$[\gamma^{+,i}_r,\gamma^{-,j}_s]=\delta_{ij}\delta_{r+s,0}$.
We take $\psi^{\pm, i}_m$, $\gamma^{\pm, i}_r$, $m, r>0$, as
annihilation operators, and $\psi^{\pm, i}_m$, $\gamma^{\pm, i}_r$, $m,
r<0$ as creation operators.
Let $\F_0$ denote the corresponding Fock space generated by the
vacuum vector $\vac$ with
$\psi^{\pm,i}_m\vac=\gamma^{\pm,i}_r\vac=0$, for $i=1,2\cdots,d$,
$m>0$ and $r>0$.

We have an action of $\hA$ with central charge $d$ on
$\F_0$ given by ($i,j\in\Z^*$ and $r,s\in\frac{1}{2}+\Z$)
\begin{align}
e_{ij}&=\sum_{p=1}^d:\psi^{+,p}_{-i}\psi^{-,p}_{j}:,\nonumber &
e_{rs}&=-\sum_{p=1}^d:\gamma^{+,p}_{-r}\gamma^{-,p}_{s}:,\nonumber\\
e_{is}&=\sum_{p=1}^d:\psi^{+,p}_{-i}\gamma^{-,p}_{s}:,\nonumber &
e_{rj}&=-\sum_{p=1}^d:\gamma^{+,p}_{-r}\psi^{-,p}_{j}:.\nonumber
\end{align}

There is also an action of $gl_d$ on $\F_0$, which is given by the formula
\begin{equation}\label{Eij}
E_{ij}=\sum_{n\in\Z^{^*}}:\psi^{+,i}_{-n}\psi^{-,j}_{n}:-
\sum_{r\in1/2+\Z}:\gamma^{+,i}_{-r}\gamma^{-,j}_{r}:.
\end{equation}

For $j\in\N$ we define the matrices $X^{-j}$ as follows:
\begin{eqnarray*}
 X^{-1}=&
\begin{pmatrix}
\gamma_{-\hf}^{-,d}&\gamma_{-\hf}^{-,d-1}&\cdots
&\gamma_{-\hf}^{-,1}\\ \psi_{-1}^{-,d}&\psi_{-1}^{-,d-1}&\cdots
&\psi_{-1}^{-,1}\\ \vdots&\vdots&\cdots &\vdots\\
\psi_{-1}^{-,d}&\psi_{-1}^{-,d-1}&\cdots &\psi_{-1}^{-,1}\\
\end{pmatrix},\allowdisplaybreaks\\
 X^{-2}=&
\begin{pmatrix}
\gamma_{-\hf}^{-,d}&\gamma_{-\hf}^{-,d-1}&\cdots
&\gamma_{-\hf}^{-,1}\\
\gamma_{-\frac{3}{2}}^{-,d}&\gamma_{-\frac{3}{2}}^{-,d-1}&\cdots
&\gamma_{-\frac{3}{2}}^{-,1}\\
\psi_{-2}^{-,d}&\psi_{-2}^{-,d-1}&\cdots &\psi_{-2}^{-,1}\\
\vdots&\vdots&\cdots &\vdots\\
\psi_{-2}^{-,d}&\psi_{-2}^{-,d-1}&\cdots &\psi_{-2}^{-,1}\\
\end{pmatrix},\allowdisplaybreaks\\
&\vdots\\&\vdots\\
 X^{-k} \equiv X^{-d}=&
\begin{pmatrix}
\gamma_{-\hf}^{-,d}&\gamma_{-\hf}^{-,d-1}&\cdots
&\gamma_{-\hf}^{-,1}\\
\gamma_{-\frac{3}{2}}^{-,d}&\gamma_{-\frac{3}{2}}^{-,d-1}&\cdots
&\gamma_{-\frac{3}{2}}^{-,1}\\ \vdots&\vdots&\cdots &\vdots\\
\gamma_{-d+\hf}^{-,d}&\gamma_{-d+\hf}^{-,d-1}&\cdots
&\gamma_{-d+\hf}^{-,1}\\
\end{pmatrix}, \quad k \ge d.\allowdisplaybreaks\\
\end{eqnarray*}
The matrices $ X^j$, for $j\in\N$, are defined similarly. Namely,
$ X^j$ is obtained from $ X^{-j}$ by replacing $\psi^{-,k}_i$ by
$\psi^{+,d-k+1}_i$ and $\gamma^{-,k}_r$ by $\gamma^{+,d-k+1}_r$.
For $0\le r\le d$ and $i\in\Z^*$, we let $X^i_{r}$ denote the
first $r\times r$ minor of the matrix $ X^i$. (We use the
definition given in \cite{CW1} for the determinant of
a matrix with Grassmannian entries.)

Given a generalized partition $\la=(\la_1,\la_2,\cdots,\la_d)$ of
length $d$, we have the shifted Frobenius notations
$\La(\la^+)=({\xi^+_{\frac{1}{2}},\xi^+_{1\frac{1}{2}},\cdots,\xi^+_{r-\frac{1}{2}}}
\mid {\xi^+_1,\xi^+_2,\cdots,\xi^+_r})$ and
$\La(\la^-)=({\xi^-_{\frac{1}{2}},\xi^-_{1\frac{1}{2}},\cdots,\xi^-_{s-\frac{1}{2}}}
\mid {\xi^-_1,\xi^-_2,\cdots,\xi^-_s})$ for the partitions $\la^+$
and $\la^-$, respectively (see \secref{Frobenius}). We let
$\La^{\hA}(\la)$ be an element of the dual space $(\hA_0)^*$ of $\hA_0$ defined by
\begin{equation}
\La^{\hA}(\la):=\sum_{{j\le r}\atop
j\in\hf\N}\xi^+_j\omega_j-\sum_{j\ge {-s}\atop
j\in-\hf\N}\xi^-_{-j}\omega_j+d\La_0.
\end{equation}

Let $\ccA_0$ be the subalgebra of $\ccA$ generated by those
$\psi^{+,i}_{m}$, $\psi^{-,i}_{m}$, $\gamma^{+,i}_{r}$ and
$\gamma^{-,i}_{r}$, $m\in\Z^*$, $r\in \frac{1}{2}+\Z$,
$i=1,2,\cdots,d$. Then $gl_d$ acts on $\ccA_0$ by the usual commutator.
This action lifts to an action of $GL_d$ on $\ccA_0$ by conjugation. The
$GL_d$ invariants of the associative algebra $\ccA_0$ is generated by $\hA^f$, hence
the $gl_d$-action on $\F_0$ commutes with the $\hA$-action.
Therefore we have the following result, which is analogous to the
$(\hgltwo, gl_d)$-duality of \cite{CW3}.

\begin{thm} \label{duality-hA}
The Lie superalgebra $\hA$ and $gl_d$ form a dual pair on $\F_0$
in the sense of Howe.  In particular,  we have the following
(multiplicity-free) decomposition of $\F_0$ with respect to their
joint action
\begin{equation*}
\F_0\cong\sum_{\la}L(\hA,\Lambda^{\hA}(\la))\otimes V_d^{\la},
\end{equation*}
where the summation is over all generalized partitions of length
$d$. Furthermore, the joint highest weight vector of the
$\la$-component is given by
\begin{equation*}
 X^{\la_{d}}_{{(\la^{-*})}'_{-\la_{d}}}\cdots
X^{-1}_{{(\la^{-*})}'_{1}} \cdot  X^1_{{(\la^+)}_1'} \cdot
X^2_{{(\la^+)}_2'}\cdots X^{\la_1}_{{(\la^+)}_{\la_1}'}\vac.
\end{equation*}
\end{thm}

\subsubsection{Unitarizable $\hA$-modules}
The restriction of the anti-linear anti-involution $\omega$ on the
Lie superalgebra $\hgltwo$ to $\hA$ gives an anti-linear
anti-involution on $\hA$,  which will also denoted by $\omega$. We have
$\omega(C)=C$, and
$$\omega(\sum_{p\in \hf\Z^{^*}}a_pe_{p-k,p}) =\sum_{p\in
\hf\Z^{^*}} (-1)^{[p]+[p-k]}\overline{a}_pe_{p,p-k},
$$
for all $\sum_{p\in \hf\Z^{^*}}a_pe_{p-k,p}\in \hA_k $ and for all
$k\in \hf\Z$. It is clear that the Fock space $\F_0$ is
a subspace of the Fock space $\F$ which is defined in the last
subsection. Since $\F$ is a unitarizable
$\hgltwo$-module with respect to the Hermitian form $\langle\
\cdot|\cdot\ \rangle$, the Fock space $\F_0$ is a unitarizable
$\hA$-module with respect to the anti-linear anti-involution
$\omega$ on $\hA$. By \thmref{duality-hA}, for each generalized
partition $\lambda$ of length $d$, the irreducible $\hA$-module
$L(\hA,\Lambda^{\hA}(\lambda))$ is unitarizable. In fact they are
all the unitarizable irreducible $\hA$-modules with respect to $\omega$.
We have the following theorem.

\begin{thm}\label{unitary-hA}
Let $M$ be an irreducible quasi-finite highest weight $\hA$-module
with highest weight $\xi$. Then $M$ is unitarizable if and
only if $\xi=\Lambda^{\hA}(\lambda)$ for some generalized
partition $\lambda$. In other words, $M$ is unitarizable if and
only if
\begin{equation*}
\La^{\hA}(\la):=\sum_{{j\le r}\atop
j\in\hf\N}\xi^+_j\omega_j-\sum_{j\ge {-s}\atop
j\in-\hf\N}\xi^-_{-j}\omega_j+d\La_0.
\end{equation*}
such that $d\in \Z_+$, $s,r\in\N$ and $\xi_j\in\Z$ for all $j$
satisfying the following conditions:
\begin{enumerate}
\renewcommand{\labelenumi}{(\roman{enumi})}
\item $\xi^+_{\hf}>\xi^+_{\frac{3}{2}}>\cdots>\xi^+_{r-\hf}
\geq 0$, $\xi^+_{1}>\xi^+_{2}>\cdots>\xi^+_{r} \geq 0$, and
$\xi^+_{r-\hf}=0$ if and only if $ r=1$ and $ \xi^+_\hf=\xi^+_1=0$,

\item $\xi^-_{\hf}>\xi^-_{\frac{3}{2}}>\cdots>\xi^-_{s-\hf}
\geq 0$, $\xi^-_{1}>\xi^-_{2}>\cdots>\xi^-_{s} \geq 0$, and
$\xi^-_{s-\hf}=0$ if and only if $s=1$ and $\xi^-_{\hf}=\xi^-_{1}=0$,

\item $\min\{\xi^+_\hf,1\}+\min\{\xi^-_{\hf},1\}+\xi^+_1+\xi^-_1 \leq d$.
\end{enumerate}
\end{thm}

\begin{proof} We already know that $L(\hA,\Lambda^{\hA}(\lambda))$ are unitarizable
irreducible quasi-finite highest weight $\hA$-modules for any
generalized partition $\lambda$. Now we are going to show that if
$M$ is a unitarizable irreducible quasi-finite highest weight
$\hA$-module with the highest weight $\xi$, then
$\xi=\Lambda^{\hA}(\lambda)$ for some generalized partition
$\lambda$. Let $\langle \ \cdot | \cdot \ \rangle$ be a positive
definite contravariant Hermitian form on $M$ and $v_\xi$ a highest
weight vector of $M$ such that $ \langle v_\xi | v_\xi\rangle =1$.
We put $\xi(e_{i,i})=\xi_i$ for all $i\in \hf\Z$. By
\thmref{criterion-hA}, there exists $r,s\in \N$ such that
$\xi=\sum_{{j\le r}\atop j\in\hf\N}\xi^+_j\omega_j-\sum_{j\ge
{-s}\atop j\in-\hf\N}\xi^-_{-j}\omega_j+d\La_0$, where $\xi^+_j$,
$\xi^-_j$, $d\in \C$. By using similar argument as those in the proof
of \thmref{unitary}, we have $\xi^+_i$, $\xi^-_j$, $d\in \Z_+$ for
$i=\hf,1,\cdots,r-\hf,r$, $j=\hf,1,\cdots,s-\hf,s$, and
$$\xi^+_{\hf}>\xi^+_{1\hf}>\cdots>\xi^+_{r-\hf} \geq 0,\qquad
\xi^+_{1}>\xi^+_{2}>\cdots>\xi^+_{r} \geq 0.$$
and
$$
\xi^-_{\hf}>\xi^-_{1\hf}>\cdots>\xi^-_{s-\hf} \geq 0, \qquad
\xi^-_{1}>\xi^-_{2}>\cdots>\xi^-_{s} \geq 0.
$$
Moreover, $\xi^+_{r-\hf}=0$ if and only if $ r=1$ and $
\xi^+_\hf=\xi^+_1=0$, and $\xi^-_{s-\hf}=0$ if and only if $s=1$ and
$\xi^-_{\hf}=\xi^-_{1}=0$.

Finally, we need to show
$\min\{\xi^+_\hf,1\}+\min\{\xi^-_{\hf},1\}+\xi^+_1+\xi^-_1 \leq
d$. Since $\langle e_{1,0}v_\xi\mid e_{1,0}v_\xi\rangle\geq 0$ and
\begin{equation*}
\begin{array}{ll}
 \langle e_{1,-1}v_\xi\mid e_{1,-1}v_\xi\rangle
 &=\langle v_\xi\mid \omega(e_{1,-1})e_{1,-1}v_\xi\rangle\\
 &=\langle v_\xi\mid (e_{-1,-1}-e_{1,1}+C)v_\xi\rangle\\
 &=-\xi^-_{1}-\xi^+_{1}+d,
\end{array}
\end{equation*}
we have $-\xi^-_{1}-\xi^+_{1}+d\geq 0$. If
$-\xi^-_{1}-\xi^+_{1}+d= 0$, then we have $e_{1,-1}v_\xi=0$.
Therefore we have
\begin{equation}\label{plusminus1}
d-\xi^-_{1}=\xi^+_1\geq 0.
\end{equation}
On the other hand, by using \eqnref{plusminus1} and
$e_{1,-1}v_\xi=0$, we have
\begin{equation*}
\begin{array}{ll}
0&= \langle e_{\hf,1}e_{1,-1}v_\xi\mid
e_{\hf,1}e_{1,-1}v_\xi\rangle\\
 &=\langle e_{\hf,-1}v_\xi\mid e_{\hf,-1}v_\xi\rangle\\
 &=\langle v_\xi\mid \omega(e_{\hf,-1})e_{\hf,-1}v_\xi\rangle\\
 &=\langle v_\xi\mid (e_{-1,-1}+e_{\hf,\hf}+C)v_\xi\rangle\\
 &=-\xi^-_1+\xi^+_{\hf}+d\\
 &=\xi^+_\hf+\xi^+_{1}.\\
\end{array}
\end{equation*}
Thus we have $\xi^+_\hf=\xi^+_{1}=0$. Similarly, we also have
\begin{equation*}
\begin{array}{ll}
0= \langle e_{-1,-\hf}e_{1,-1}v_\xi\mid
e_{-1,-\hf}e_{1,-1}v_\xi\rangle
 &=\langle e_{1,-\hf}v_\xi\mid e_{1,-\hf}v_\xi\rangle\\
 &=\langle v_\xi\mid \omega(e_{1,-\hf})e_{1,-\hf}v_\xi\rangle\\
 &=\langle v_\xi\mid (e_{-\hf,-\hf}+e_{1,1}-C)v_\xi\rangle\\
 &=-\xi^-_{\hf}+\xi^+_{1}-d\\
 &=-\xi^-_\hf-\xi^-_{1}.\\
 \end{array}
\end{equation*}
Thus we have $\xi^-_\hf=\xi^-_{1}=0$. Hence
$\xi^+_\hf=\xi^+_{1}=\xi^-_\hf=\xi^-_{1}=0$ and
$\min\{\xi^+_\hf,1\}+\min\{\xi^-_{\hf},1\}+\xi^+_1+\xi^-_1 =0\leq
d$ when $-\xi^-_{1}-\xi^+_{1}+d= 0$.

Now we assume that $-\xi^-_{1}-\xi^+_{1}+d>0$. If $\xi^+_{\hf}=0$,
then $\min\{\xi^+_\hf,1\}+\min\{\xi^-_{\hf},1\}+\xi^+_1+\xi^-_1
\le 1+\xi^+_1+\xi^-_1 \leq d$ and the proof of the theorem is
completed. Otherwise, we may assume that $\xi^+_{\hf}>0$. Also we
may assume that $-\xi^-_{1}-\xi^+_{1}+d=1$, otherwise the proof is
also completed. Now we consider
\begin{equation*}
\begin{array}{ll}
 \langle e_{1,-1}e_{1,\hf}v_\xi\mid e_{1,-1}e_{1,\hf}v_\xi\rangle
 &=\langle e_{1,\hf}v_\xi\mid e_{-1,1}e_{1,-1}e_{1,\hf}v_\xi\rangle\\
 &=\langle e_{1,\hf}v_\xi\mid (e_{-1,-1}-e_{1,1}+C)e_{1,\hf}v_\xi\rangle\\
 &=(-\xi^-_{1}-\xi^+_{1}-1+d)\langle e_{1,\hf}v_\xi\mid e_{1,\hf}v_\xi\rangle\\
 &=(-\xi^-_{1}-\xi^+_{1}-1+d)(\xi^+_{\hf}+\xi^+_{1})\\
 &=0.
\end{array}
\end{equation*}
Thus we have $\langle e_{-\hf,1}e_{1,-1}e_{1,\hf}v_\xi\mid
 e_{-\hf,1}e_{1,-1}e_{1,\hf}v_\xi\rangle=0$ and
\begin{equation*}
\begin{array}{ll}
 \langle e_{-\hf,1}e_{1,-1}e_{1,\hf}v_\xi\mid
 e_{-\hf,1}e_{1,-1}e_{1,\hf}v_\xi\rangle
 &=\langle e_{-\hf,-1}e_{1,\hf}v_\xi\mid
 e_{-\hf,-1}e_{1,\hf}v_\xi\rangle\\
 &=\langle e_{1,\hf}v_\xi\mid -e_{-1,-\hf}e_{-\hf,-1}e_{1,\hf}v_\xi\rangle\\
 &=\langle e_{1,\hf}v_\xi\mid (e_{-1,-1}+e_{-\hf,-\hf})e_{1,\hf}v_\xi\rangle\\
 &=(-\xi^-_{1}-\xi^-_{\hf})\langle e_{1,\hf}v_\xi\mid e_{1,\hf}v_\xi\rangle\\
 &=(-\xi^-_{1}-\xi^-_{\hf})(\xi^+_{\hf}+\xi^+_{1}).\\
\end{array}
\end{equation*}
Therefore, we have $\xi^-_{\hf}=0$ and
$\min\{\xi^+_\hf,1\}+\min\{\xi^-_{\hf},1\}+\xi^+_1+\xi^-_1 \le
1+\xi^+_1+\xi^-_1 \leq d$. This completes the proof of the
theorem.
\end{proof}

\subsection{Unitarizable $\hC$-modules and
their Fock space realizations}\label{unitarizable-hC}
\subsubsection{Free field realization of $\hC$ and $(\hC, Sp(2d))$-duality}
Let us first recall some facts about the complex
symplectic group $Sp(2d)$ (see, e.g., \cite{BT, FH, GW}).
Consider the non-degenerate skew-symmetric
bilinear form on $\C^{2d}$ given by the $2d\times
2d$ matrix
\begin{equation*}
J_{2d}=\begin{pmatrix}
0&I_{d}\\
-I_{d}&0
\end{pmatrix},
\end{equation*}
where $I_d$ is the $d\times d$ identity matrix. The symplectic
group $Sp(2d)$ is the subgroup of $GL(2d)$ which consists of those
$A\in GL(2d)$ with $A^tJ_{2d}A=J_{2d}$, where $A^t$ is the
transpose of the matrix $A$. The Lie algebra of $Sp(2d)$ is
$\fsp(2d)$ which consists of those $A\in gl(2d)$ with
$A^tJ_{2d}+J_{2d}A=0$. Denote by $e_{ij}$ the elementary matrix
with $1$ in the $i$-th row and $j$-th column and $0$ elsewhere.
Then $\h:=\sum_{1\le i\le d}\C (e_{ii}-e_{d+i,d+i})$ is a Cartan
subalgebra, while $\bb:=\sum_{i\le j \le d}\C
(e_{ij}-e_{j+d,i+d})+ \sum_{i\le j}\C (e_{i,j+d}+e_{j+d,i})$ is
the standard Borel subalgebra containing $\h$. Let
$h_{i}=e_{ii}-e_{d+i,d+i}$.

We write an element $\lambda\in\h^*$ as $\lambda= (\lambda_1,
\lambda_2, \cdots, \lambda_d)$ where $\lambda_i=\lambda(h_i)$
for $i=1, 2, \cdots, d$.
Let $V^\la_{\fsp(2d)}$ denote the irreducible
$\fsp(2d)$-module with highest weight $\la\in\h^*$ defined
with respect to the standard Borel subalgebra. Then $V^\la_{\fsp(2d)}$ is
finite-dimensional if and only if
$\la_1\ge\la_2\ge\cdots\ge\la_d$ and $\la_i\in\Z_+$ for
$i=1,\cdots,d$. Furthermore every such representation lifts to a
unique irreducible representation of $Sp(2d)$, which is denoted by
$V^\la_{Sp(2d)}$ and so we obtain an obvious parametrisation of
$Sp(2d)$-highest weights in terms of Young diagrams $\la$ with
$l(\la)\le d$.

We let $\epsilon_i\in\h^*$ so that
$\epsilon_i(h_{j})=\delta_{ij}$. We put $z_i=e^{\epsilon_i}$ when
dealing with characters of $Sp(2d)$.

Introduce the following operators on the Fock space $\F_0$:
\begin{align}
\label{sp+} E^{sp+}_{ij}&=\sum_{
n\in\N}:\psi^{+,i}_{-n}\psi^{+,j}_{n}: -\sum_{n\in
-\N}:\psi^{+,i}_{-n}\psi^{+,j}_{n}:+
\sum_{r\in1/2+\Z}:\gamma^{+,i}_{-r}\gamma^{+,j}_{r}:,\\
\label{sp-} E^{sp-}_{ij}&=\sum_{
n\in\N}:\psi^{-,i}_{-n}\psi^{-,j}_{n}: -\sum_{n\in
-\N}:\psi^{-,i}_{-n}\psi^{-,j}_{n}:-
\sum_{r\in1/2+\Z}:\gamma^{-,i}_{-r}\gamma^{-,j}_{r}:,
\end{align}
where $1\le i,j\le d$. It is clear that \eqnref{sp+} and
\eqnref{sp-} together with \eqnref{Eij} form a basis for the Lie
algebra $\fsp (2d)$. The action of the Lie algebra $\fsp(2d)$ on
the Fock space $\F_0$ can be lifted to an action of Lie group
$Sp(2d)$. Moreover $\F_0$ is a direct sum of finite dimensional
irreducible $Sp(2d)$-modules.

On the other hand, $Sp(2d)$ acts on $\ccA_0$ by conjugation.
It is not hard to see that the $Sp(2d)$-invariants of
the associative algebra $\ccA_0$ is generated by the following
combinations of the elements of \eqref{realization}:
\begin{align}
\label{C1}&C; &\te_{r,s}=e_{r,s}-e_{-s,-r};\\
\label{C2}&\te_{i,j}=e_{i,j}-e_{-j,-i},\ \  ij>0;
&\te_{i,j}=e_{i,j}+e_{-j,-i},\ \  ij<0;\\
\label{C3}&\te_{i,r}=\te_{-r,-i}=e_{i,r}+e_{-r,-i}, \  i>0;
&\te_{i,r}=-\te_{-r,-i}=e_{i,r}-e_{-r,-i},\ i<0,
\end{align}
where $i,j\in \Z^*$ and $r,s\in \hf+\Z$. Note that \eqnref{C1},
\eqnref{C2} and \eqnref{C3} form the Lie superalgebra $\hC^f$ and
hence $Sp(2d)$ commutes with $\hC$. By a result of Howe
\cite{H1}, we have the following multiplicity-free decomposition
\begin{equation}\label{Sp-duality}
\F_0\cong\sum_{\la}W^\la\otimes V_{Sp(2d)}^{\la}
\end{equation}
of $\F_0$ with respect to the joint action of $\hC$ and $Sp(2d)$,
where the summation is over a subset of all partitions of length
$d$. Here $W^\la$ denotes an irreducible module over $\hC$.

For any given partition $\la=(\la_1,\la_2,\cdots,\la_d)$ of length $d$,
we have the shifted Frobenius notation
$\La(\la)=({\xi_{\frac{1}{2}},\xi_{1\frac{1}{2}},\cdots,\xi_{r-\frac{1}{2}}}
\mid {\xi_1,\xi_2,\cdots,\xi_r})$ (see
\secref{Frobenius}). Let $\La^{\hC}(\la)$ be an element of the
dual space $\hC_0^*$ of $\hC_0$ defined by
\begin{equation}
\La^{\hC}(\la):=\sum_{{j\le r}\atop
j\in\hf\N}\xi_j\omega_j+d\La_0.
\end{equation}

Let $\la$ be a partition of length $d$. Then the following vector
\begin{equation}\label{Sp-hwv}
X^1_{\la_1'} \cdot  X^2_{\la_2'}\cdots
X^{\la_1}_{{\la}_{\la_1}'}\vac
\end{equation}
is a highest weight vector of the Lie superalgebra $\hC$, which is
defined with respect to the Borel subalgebra $\hC_0\oplus\hC_+$ of $\hC$ (see Remark
\ref{embeddings}) and has weight $\La^{\hC}(\la)$. This follows from the fact
that the vector is actually a highest weight vector of $\hA$. It is easy to see that
the vector given by \eqnref{Sp-hwv} is also annihilated by \eqnref{sp+}.
Thus it is a joint highest weight vector of $Sp(2d)$ and $\hC$. Note that
the weights of the vectors of the form \eqnref{Sp-hwv} with respect to the
Lie group $Sp(2d)$ are exactly the weights associated with all
partitions of length $d$. By the decomposition
of \eqnref{Sp-duality}, we have the following theorem.

\begin{thm} \label{duality-hC}
The Lie superalgebra $\hC$ and $Sp(2d)$ form a dual pair on $\F_0$
in the sense of Howe.  Furthermore we have the following
(multiplicity-free) decomposition of $\F_0$ with respect to their
joint action
\begin{equation*}
\F_0\cong\sum_{\la}L(\hC,\Lambda^{\hC}(\la))\otimes
V_{Sp(2d)}^{\la},
\end{equation*}
where the summation is over all partitions of length $d$.
The joint highest weight vector of the
$\la$-component is given by
\begin{equation*}
X^1_{\la_1'} \cdot  X^2_{\la_2'}\cdots
X^{\la_1}_{{\la}_{\la_1}'}\vac.
\end{equation*}
\end{thm}

\subsubsection{Unitarizable $\hC$-modules}
The restriction of the anti-linear anti-involution $\omega$ on the
Lie superalgebra $\hA$ to $\hC$ gives rise to an anti-linear
anti-involution on $\hC$, which will also be denoted by $\omega$. It satisfies
$\omega(C)=C$ and
$$\omega(\sum_{p\in \hf\Z^{^*}}a_p\te_{p-k,p}) =\sum_{p\in
\hf\Z^{^*}} (-1)^{[p]+[p-k]}\overline{a}_p\te_{p,p-k},
$$
for all $\sum_{p\in \hf\Z^{^*}}a_p\te_{p-k,p}\in (\hC)_k $ and for
all $k\in \hf\Z$. Since the Fock space $\F_0$ is a unitarizable
$\hA$-module with the positive definite contravariant Hermitian form $\langle\
\cdot|\cdot\ \rangle$, it is also a unitarizable
$\hC$-module with respect to the anti-linear anti-involution
$\omega$. By \thmref{duality-hC}, the irreducible $\hC$-module
$L(\hC,\Lambda^{\hC}(\lambda))$ is unitarizable with respect to
$\omega$ for every partition $\lambda$ of length $d$.
In fact these modules exhaust all the irreducible
quasi-finite highest weight $\hC$-modules, which are
unitarizable with respect to $\omega$. We have the following theorem.

\begin{thm}\label{unitary-hC}
Let $M$ be an irreducible quasi-finite highest weight $\hC$-module
with highest weight $\xi$. Then $M$ is unitarizable if and
only if $\xi=\Lambda^{\hC}(\lambda)$ for some partition $\lambda$
of length $d$. In other words, $M$ is unitarizable if and only if
\begin{equation*}
\La^{\hC}(\la):=\sum_{{j\le r}\atop j\in\hf\N}\xi_j\omega_j+d\La_0
\end{equation*}
such that $d\in \Z_+$, $r\in\N$ and $\xi_j\in\Z$ for all $j$
satisfying the following conditions:
\begin{enumerate}
\renewcommand{\labelenumi}{(\roman{enumi})}
\item $\xi_{\hf}>\xi_{1\hf}>\cdots>\xi_{r-\hf} \geq 0$,
$\xi_{1}>\xi_{2}>\cdots>\xi_{r} \geq 0$, and $\xi_{r-\hf}=0$ if
only if $ r=1$ and $ \xi_\hf=\xi_1=0$,

\item $\min\{\xi_\hf,1\}+\xi_1 \leq d$.
\end{enumerate}
\end{thm}

\begin{proof} By the argument above,
$L(\hC,\Lambda^{\hC}(\lambda))$ are unitarizable
irreducible quasi-finite highest weight $\hC$-modules for any
partition $\lambda$. Now we are going to show that if $M$ is a
unitarizable irreducible quasi-finite highest weight $\hC$-module
with the highest weight $\xi$, then $\xi=\Lambda^{\hC}(\lambda)$
for some partition $\lambda$. Let $\langle \ \cdot | \cdot \
\rangle$ be a positive definite contravariant Hermitian form on
$M$ and $v_\xi$ a highest weight vector of $M$ such that $ \langle
v_\xi | v_\xi\rangle =1$. We put $\xi(\te_{i,i})=\xi_i$ for all
$i\in \hf\N$. By \thmref{criterion-hC}, there exists $r,s\in \N$
such that $\xi=\sum_{{j\le r}\atop
j\in\hf\N}\xi_j\omega_j+d\La_0$, where $\xi_j$, $d\in \C$. Using similar
arguments as in the proof of \thmref{unitary}, we can show that
$\xi_i\in \Z_+$ for $i=\hf,1,\cdots,r-\hf,r$, and
$$\xi_{\hf}>\xi_{1\hf}>\cdots>\xi_{r-\hf} \geq 0,\qquad
\xi_{1}>\xi_{2}>\cdots>\xi_{r} \geq 0.$$ Moreover, we also have
$\xi_{r-\hf}=0$ if only if $ r=1$ and $ \xi_\hf=\xi_1=0$.

Now we choose a large positive integer $n$ such that
$\xi(\te_{n,n})=0$. Consider the subalgebra $sl(2,\C)$ with
standard basis $\{ -\te_{n,n}+C,\ \hf\te_{-n,n},\
\hf\te_{n,-n}\}$. Note that $\omega(\hf\te_{-n,n})=\hf\te_{n,-n}$,
$\omega(\hf\te_{n,-n})=\hf\te_{-n,n}$ and
$\omega(-\te_{n,n}+C)=-\te_{n,n}+C$. A standard result on
unitarizable $sl(2, \C)$-modules (see, e.g., \cite{K2})
leads to $d=\xi(C)=\xi(-\te_{n,n}+C)\in\Z_+$.

Finally, we need to show $\min\{\xi_\hf,1\}+\xi_1\leq d$. Since
$\langle \te_{1,-1}v_\xi\mid \te_{1,-1}v_\xi\rangle\geq 0$ and
\begin{equation*}
\begin{array}{ll}
 \langle \hf\te_{1,-1}v_\xi\mid \hf\te_{1,-1}v_\xi\rangle
 &=\langle v_\xi\mid \hf\omega(\hf\te_{1,-1})\te_{1,-1}v_\xi\rangle\\
 &=\langle v_\xi\mid (-\te_{1,1}+C)v_\xi\rangle\\
 &=d-\xi_{1},
\end{array}
\end{equation*}
we have $d\geq \xi_{1}$. If $ d > \xi_{1}$, the proof the theorem
is completed. Otherwise, we assume $d-\xi_1= 0$ and hence
$\te_{1,-1}v_\xi=0$.

On the other hand, by using $\te_{1,-1}v_\xi=0$ and $d-\xi_1= 0$,
we have
\begin{equation*}
\begin{array}{ll}
0 &= \langle \te_{\hf,1}\te_{1,-1}v_\xi\mid
\te_{\hf,1}\te_{1,-1}v_\xi\rangle\\
 &=\langle \te_{\hf,-1}v_\xi\mid \te_{\hf,-1}v_\xi\rangle\\
 &=\langle v_\xi\mid \omega(\te_{\hf,-1})\te_{\hf,-1}v_\xi\rangle\\
 &=\langle v_\xi\mid (2C-\te_{1,1}+\te_{\hf,\hf})v_\xi\rangle\\
 &=2d-\xi_{1}+\xi_{\hf}\\
 &=d+\xi_\hf.\\
 \end{array}
\end{equation*}
Hence we have $d=\xi_\hf=\xi_{1}=0$ and
$\min\{\xi_\hf,1\}+\xi_1\leq d$. This completes the proof of the
theorem.
\end{proof}

\subsection{Unitarizable $\hD$-modules and their
Fock space realizations}\label{unitarizable-hD}
In this subsection we will construct two types of free field
realizations of $\hD$ which are respectively associated with the
$(\hD, O(2d))$ and  $(\hD, O(2d+1))$-dualities.

\subsubsection{Basic facts on representations of
the complex orthogonal group}
We start by recalling some facts about finite dimensional
representations of the complex orthogonal group $O(k)$ (se, e.g., \cite{BT, FH, GW}).
We first consider the case when $k=2d$ is even. Consider the
non-degenerate symmetric bilinear form on
$\C^{2d}$ given by the $2d\times 2d$ matrix
\begin{equation*}
K_{2d}=\begin{pmatrix}
0&I_{d}\\
I_{d}&0
\end{pmatrix},
\end{equation*}
where $I_d$ is the $d\times d$ identity matrix. The orthogonal
group $O(2d)$ is the subgroup of $GL(2d)$ which consists of those
$A\in GL(2d)$ with $A^tK_{2d}A=K_{2d}$, where $A^t$ is the
transpose of the matrix $A$. The Lie algebra of $O(2d)$ is
$\fso(2d)$ which consists of those $A\in gl(2d)$ with
$A^tK_{2d}+K_{2d}A=0$. Denote by $e_{ij}$ the elementary matrix
with $1$ in the $i$-th row and $j$-th column and $0$ elsewhere.
Let $h_{i}:=e_{ii}-e_{d+i,d+i}$,
$E^{so+}_{ij}:=e_{i,j+d}-e_{j,i+d}$ and
$E^{so-}_{ij}:=e_{i+d,j}-e_{j+d,i}$ for $1\le i,j\le d$. Then
$\h:=\sum_{1\le i\le d}\C h_i$ is a Cartan subalgebra, while
$\bb:= \sum_{1\le i\le j\le d}\C (e_{i,j}-e_{j+d,i+d})+\sum_{1\le
i, j\le d} \C E^{so+}_{ij}$ is the standard Borel subalgebra
containing $\h$.

Write an element $\lambda\in\h^*$ as $\lambda=(\la_1, \la_2, \cdots,\la_d
)$, where $\la_i=\la(h_i)$, for $i=1, 2, \cdots, d$.
Let $V^\la_{\fso(2d)}$ denote the irreducible highest weight $\fso(2d)$-module
with highest weight $\la\in\h^*$ defined
with respect to the standard Borel subalgebra. Then
$V^\la_{\fso(2d)}$ is finite dimensional if and only if
$\la_1\ge\la_2\ge\cdots\ge |\la_d|$ with either $\la_i\in\Z$ or
else $\la_i\in\hf + \Z$  for all $i=1,\cdots,d$. Furthermore,
the $\fso(2d)$-module $V^\la_{\fso(2d)}$ lifts to an $SO(2d)$-module
if and only if $\la_1\ge\la_2\ge\cdots\ge |\la_d|$ with $\la_i\in\Z$
for all $i=1,\cdots,d$.

Let $V$ be a finite-dimensional irreducible $O(2d)$-module.  When
regarded as an $\fso(2d)$-module we have the following
possibilities:
\begin{itemize}
\item[(i)] $V$ is a direct sum of two irreducible
$\fso(2d)$-modules with integral highest weights $(\la_1,\la_2,\cdots,\la_d)$
and $(\la_1,\la_2,\cdots,\la_{k-1},-\la_d)$ respectively, where
$\la_k>0$.
 \item[(ii)] $V$ is an irreducible $\fso(2d)$-module with integral
 highest weight of the form $(\la_1,\la_2,\cdots,\la_{d-1},0)$.
 \end{itemize}
In the first case, that is when $V$ is the direct sum of the two
irreducible $\fso(2d)$-modules,  we denote $V$ by
$W_{O(2d)}^{\widetilde\la}$, where we let
${\widetilde\la}=(\la_1,\la_2,\cdots,\la_{d-1},\la_d>0)$.  In the
second case there are two possible choices of $V$, which we denote
by $W_{O(2d)}^{\widetilde\la}$ and
$W_{O(2d)}^{\widetilde\la}\otimes {\rm det}$, respectively.
Recalling that $O(2d)$ is a semidirect product of $SO(2d)$ and
$\Z_2$. Thus the $O(2d)$-modules $W_{O(2d)}^{\widetilde\la}$ and
$W_{O(2d)}^{\widetilde\la}\otimes {\rm det}$ restrict to isomorphic
$SO(2d)$-modules.  However as $O(2d)$-modules they differ by the
determinant representation so that we may distinguish these two
modules as follows: consider the element $\tau\in O(2d)-SO(2d)$
that switches the basis vector $e^{d}$ with $e^{2d}$ and leaves
all other basis vectors of $\C^{2d}$ invariant. We declare
$W_{O(2d)}^{\widetilde\la}$ to be the $O(2d)$-module on which
$\tau$ transforms an $SO(d)$-highest weight vector trivially.
Note that $\tau$ transforms an $SO(2d)$-highest weight vector in
the $O(2d)$-module $W_{O(2d)}^{\widetilde\la}\otimes{\rm det}$ by
$-1$.

We may associate Young diagrams $\la$ of length $2d$ to these
integral highest weights of $O(2d)$ as follows (cf.~\cite{H2}).  For
$\la_1\ge\la_2\ge\cdots\ge\la_d>0$ with $\la_i\in\Z_+$ for all $i$,
we have an obvious Young diagram of length $2d$ by putting
$\la:=(\la_1,\la_2,\cdots,\la_{d},0,\cdots,0)$. When $\la_d=0$, we
associate to the highest weight of $W_{O(2d)}^{\widetilde\la}$ the
usual Young diagram of length $2d$ by putting
$\la:=(\la_1,\la_2,\cdots,\la_{d},0,\cdots,0)$. We put
$V^\la_{O(2d)}:=W_{O(2d)}^{\widetilde\la}$. To the highest weight
of $W^{\tilde\la}_{O(2d)}\otimes {\rm det}$ we associate the Young diagram
$\bar\la$ obtained from $\la$ by replacing its first
column by a column of length $2d-\la'_1$.
Let $V^{\bar\la}_{O(2d)}:=W_{O(2d)}^{\widetilde\la}\otimes {\rm det}$.

Hereafter, we shall adopt the following convention.  Given any
partition $\la$ of length $2 d$ satisfying the condition
$\la'_1+\la'_2\le 2 d$, we denote by $\bar\la$ the partition
obtained from $\la$ by replacing its first
column by a column of length $2d-\la'_1$. There is a one to one
correspondence between the finite dimensional irreducible
$O(2d)$-representations and the partitions $\la$ of length $2 d$
satisfying the condition $\la'_1+\la'_2\le 2 d$.

Next consider the case when $k=2d+1$ is odd. Take the
non-degenerate symmetric bilinear form on $\C^{2d+1}$ given
by the $(2d+1)\times (2d+1)$ matrix
\begin{equation*}
K_{2d+1}=\begin{pmatrix}
0&0&I_{d}\\
0&1&0\\
I_{d}&0&0
\end{pmatrix},
\end{equation*}
where $I_d$ is the $d\times d$ identity matrix. The orthogonal
group $O(2d+1)$ is the subgroup of $GL(2d+1)$ which consists of
those $A\in GL(2d+1)$ with $A^tK_{2d+1}A=K_{2d+1}$, where $A^t$ is
the transpose of the matrix $A$. The Lie algebra of $O(2d+1)$ is
$\fso(2d+1)$ which consists of those $A\in gl(2d+1)$ with
$A^tK_{2d+1}+K_{2d+1}A=0$. Denote by $e_{ij}$ the elementary
matrix with $1$ in the $i$-th row and $j$-th column and $0$
elsewhere. Let $h_{i}:=e_{ii}-e_{d+1+i,d+1+i}$,
$E^{so+}_{i}:=e_{i,d+1}-e_{d+1,d+1+i}$,
$E^{so-}_{j}:=e_{d+1,j}-e_{j+d+1,d+1}$,
$E^{so+}_{ij}:=e_{i,j+d+1}-e_{j,i+d+1}$ and
$E^{so-}_{ij}:=e_{i+d+1,j}-e_{j+d+1,i}$ for $1\le i,j\le d$. Then
$\h:=\sum_{1\le i\le d}\C h_i$ is a Cartan subalgebra, while
$\bb:=\sum_{1\le i\le j\le d}\C
(e_{i,j}-e_{j+d+1,i+d+1})+\sum_{1\le i\le d} \C
E^{so+}_{i}+\sum_{1\le i, j\le d} \C E^{so+}_{ij}$ is the standard
Borel subalgebra containing $\h$.

Write an element $\lambda\in\h^*$ as $\lambda=(\la_1, \la_2, \cdots,\la_d
)$, where $\la_i=\la(h_i)$, for $i=1, 2, \cdots, d$.
Let $V^\la_{\fso(2d+1)}$ denote the irreducible highest weight $\fso(2d+1)$-module
with highest weight $\la\in\h^*$ defined
with respect to the standard Borel subalgebra. Then
$V^\la_{\fso(2d+1)}$ is finite dimensional if and only if
$\la_1\ge\la_2\ge\cdots\ge \la_d$ with either $\la_i\in\Z_+$ or
$\la_i\in\hf + \Z_+$ for all $i=1,\cdots,d$.
Furthermore $V^\la_{\fso(2d+1)}$
lifts to a representation of $SO(2d+1)$ if and only if
$\la_i\in\Z_+$.

Recall that $O(2d+1)$ is a direct product of $SO(2d+1)$ and
$\Z_2$. Thus any finite-dimensional irreducible representation of
$O(2d+1)$, when regarded as an $SO(2d+1)$-module, remains
irreducible. Conversely an irreducible representation of
$SO(2d+1)$ gives rise to two non-isomorphic $O(2d+1)$-modules that
differ from each other by the determinant representation ${\rm
det}$. We let $W^{\widetilde\la}_{O(2d+1)}$ stand for the
irreducible $O(2d+1)$-module corresponding to
${\widetilde\la}=(\la_1\ge\la_2\ge\cdots\ge\la_k\ge 0)$ on which
the element $-I_{2d+1}$ transforms trivially, so that
$\{W^{\widetilde\la}_{O(2d+1)},W^{\widetilde\la}_{O(2d+1)}\otimes{\rm
det}\}$ with ${\widetilde\la}$ ranging over all partitions with
length $d$ as above is a complete set of finite-dimensional
non-isomorphic irreducible $O(2d+1)$-modules, where $I_{2d+1}$ is
the identity $(2d+1)\times(2d+1)$ matrix.

Similarly as before we may associate Young diagrams of length
$2d+1$ to these $O(2d+1)$-highest weights.  For the highest weight
${\widetilde\la}=(\la_1\ge\la_2\cdots\ge\la_d\ge 0)$ of
$W^{\widetilde\la}_{O(2d+1)}$, we have an obvious Young diagram
$\la=(\la_1,\la_2\cdots,\la_d,0,\cdots,0)$ of length $2d+1$.
Let $V^{\la}_{O(2d+1)}:=W^{\widetilde\la}_{O(2d+1)}$.

To the highest weight of $W^{\widetilde\la}_{O(2d+1)}\otimes{\rm
det}$ we associate the Young diagram $\bar\la$ obtained from
$\la$ by replacing its first column by a column
of length $2d+1-\la'_1$. In this case, we let
$V^{\bar\la}_{O(2d+1)}:=W^{\widetilde\la}_{O(2d+1)}\otimes{\rm det}$.

Hereafter, we adopt the convention that for any given any partition
$\la$ of length $2 d+1$ satisfying
$\la'_1+\la'_2\le 2 d+1$, we denote by $\bar\la$ the partition of
length $2 d+1$ obtained from $\la$ by replacing its first column by a column
of length $2d+1-\la'_1$. There is a one to one correspondence between
the finite dimensional irreducible representations of $O(2d +1)$ and
the partitions $\la$ of length $2 d+1$ satisfying
$\la'_1+\la'_2\le 2 d+1$.

Let $\epsilon_i\in\h^*$ so that $\epsilon_i(h_j)=\delta_{ij}$.  We
put $z_i=e^{\epsilon_i}$ when dealing with characters of $O(2d)$
and $O(2d+1)$.

\subsubsection{Free field realization of $\hD$ and
$(\hD, O(2d))$-duality}
Let us consider the realization of $\hD$ on the Fock space $\F_0$ related to the
$(\hD, O(2d))$-duality. Introduce the following operators on $\F_0$:
\begin{align}
\label{so+} E^{so+}_{ij}&=\sum_{
n\in\Z^*}:\psi^{+,i}_{-n}\psi^{+,j}_{n}: +
\sum_{r\in1/2+\Z_+}:\gamma^{+,i}_{-r}\gamma^{+,j}_{r}:
-\sum_{r\in-1/2-\Z_+}:\gamma^{+,i}_{-r}\gamma^{+,j}_{r}:,\\
\label{so-} E^{so-}_{ij}&=\sum_{
n\in\Z^*}:\psi^{-,i}_{-n}\psi^{-,j}_{n}: -
\sum_{r\in1/2+\Z_+}:\gamma^{-,i}_{-r}\gamma^{-,j}_{r}:+
\sum_{r\in-1/2-\Z_+}:\gamma^{-,i}_{-r}\gamma^{-,j}_{r}:,
\end{align}
where $1\le i,j\le d$. It is easy to see that \eqnref{so+} and
\eqnref{so-} together with \eqnref{Eij} form a basis for the Lie
algebra $\fso (2d)$. The action of the Lie algebra $\fso(2d)$ on
the Fock space $\F_0$ can be lifted to an action of Lie group
$SO(2d)$ and extend to the action of the Lie group $O(2d)$.
Moreover $\F_0$ is a direct sum of finite dimensional irreducible modules
over $O(2d)$.

On the other hand, $O(2d)$  acts on $\ccA_0$ by
conjugation. It is not hard to see that the $O(2d)$-invariants in the
associative algebra $\ccA_0$ is generated by the following combinations
of the elements of \eqnref{realization}:
\begin{align}
\label{D1}&C; &\te_{i,j}=e_{i,j}-e_{-j,-i};\\
\label{D2}&\te_{r,s}=e_{r,s}-e_{-s,-r},\  rs>0;
&\te_{r,s}=e_{r,s}+e_{-s,-r},\  rs<0;\\
\label{D3}&\te_{i,r}=\te_{-r,-i}=e_{i,r}+e_{-r,-i}, \ r>0;
&\te_{i,r}=-\te_{-r,-i}=e_{i,r}-e_{-r,-i},\ r<0,
\end{align}
where $i,j\in \Z^*$ and $r,s\in \hf+\Z$. Note that \eqnref{D1},
\eqnref{D2} and \eqnref{D3} form the Lie superalgebra $\hD^f$ in $\ccA_0$.
Therefore the $\hD$-action on $\F_0$ commutes with the $O(2d)$-action.
Following the general reasoning of
\cite{H2}, we have the following
multiplicity-free decomposition of $\F_0$ with respect to the
joint action of $\hD$ and $O(2d)$:
\begin{equation}\label{evenduality}
\F_0\cong\sum_{\la}W^\la\otimes V_{O(2d)}^{\la},
\end{equation}
where the summation is over a subset of all partitions of length
with $\la'_1+\la'_2\le 2d$, and $W^\la$ is a certain irreducible
module over $\hD$.

For each $j\in\{1,2,\cdots,d\}$, we define the $d\times d$ matrix
$\tX^{j}$ as follows:
\begin{eqnarray*}
 \tX^{j}:=&
\begin{pmatrix}
 \gamma_{-\hf}^{+,1}&\gamma_{-\hf}^{+,2}&\cdots&\gamma_{-\hf}^{+,d-1}
 &\gamma_{-\hf}^{-,d}\\
 \gamma_{-\frac{3}{2}}^{+,1}&\gamma_{-\frac{3}{2}}^{+,2}
 &\cdots&\gamma_{-\frac{3}{2}}^{+,d-1}&\gamma_{-\frac{3}{2}}^{-,d}\\
 \vdots&\vdots&\cdots&\vdots &\vdots\\
 \gamma_{-\frac{2j-1}{2}}^{+,1}&\gamma_{-\frac{2j-1}{2}}^{+,2}
 &\cdots&\gamma_{-\frac{2j-1}{2}}^{+,d-1}&\gamma_{-\frac{2j-1}{2}}^{-,d}\\
 \psi_{-j}^{+,1}&\psi_{-j}^{+,2}&\cdots&\psi_{-j}^{+,d-1}
 &\psi_{-j}^{-,d}\\
 \vdots&\vdots&\cdots&\vdots &\vdots\\
 \psi_{-j}^{+,1}&\psi_{-j}^{+,2}&\cdots&\psi_{-j}^{+,d-1}
 &\psi_{-j}^{-,d}\\
\end{pmatrix}.
\end{eqnarray*}
For any integer $j\ge d$, $\tX^{j}:=\tX^{d}$. Note that the matrix
$ \tX^j$ is obtained from $ X^{j}$ by replacing its last column by
$(\gamma_{-\hf}^{-,d},\gamma_{-\frac{3}{2}}^{-,d},\cdots,
\gamma_{-\frac{2j-1}{2}}^{-,d}\psi_{-j}^{-,d}\cdots\psi_{-j}^{-,d})$.
For $0\le r\le d$ and $i > 0$, we let $\tX^i_{r}$ denote the first
$r\times r$ minor of the matrix $ \tX^i$.

We define the $2d\times 2d$ matrix $\Gamma$ as follows:
\begin{eqnarray*}
 \Gamma :=
\begin{pmatrix}
 \gamma_{-\hf}^{+,1}&\gamma_{-\hf}^{+,2}&\cdots&\gamma_{-\hf}^{+,d}
 &\gamma_{-\hf}^{-,d}&\gamma_{-\hf}^{-,d-1}&\cdots&\gamma_{-\hf}^{-,1}\\
 \psi_{-1}^{+,1}&\psi_{-1}^{+,2}&\cdots&\psi_{-1}^{+,d}
 &\psi_{-1}^{-,d}&\psi_{-1}^{-,d-1}&\cdots&\psi_{-1}^{-,1}\\
 \psi_{-1}^{+,1}&\psi_{-1}^{+,2}&\cdots&\psi_{-1}^{+,d}
 &\psi_{-1}^{-,d}&\psi_{-1}^{-,d-1}&\cdots&\psi_{-1}^{-,1}\\
 \vdots&\vdots&\cdots&\vdots&\vdots&\vdots&\cdots &\vdots\\
 \psi_{-1}^{+,1}&\psi_{-1}^{+,2}&\cdots&\psi_{-1}^{+,d}
 &\psi_{-1}^{-,d}&\psi_{-1}^{-,d-1}&\cdots&\psi_{-1}^{-,1}\\
\end{pmatrix}.
\end{eqnarray*}
For $r > d$, we let $\Gamma_r$ denote the first $r\times r$ minor
of the matrix $ \Gamma$.

Given a partition $\la=(\la_1,\la_2,\cdots,\la_{2d})$ of length
$2d$ with $\la'_1+\la'_2\le 2d$, we have the shifted Frobenius
notation
$\La(\la)=({\xi_{\frac{1}{2}},\xi_{1\frac{1}{2}},\cdots,\xi_{r-\frac{1}{2}}}
\mid {\xi_1,\xi_2,\cdots,\xi_r})$ for the partition $\la$ (see
\secref{Frobenius}). Let $\La^{\hD}(\la)$ an element of the
dual space $\hD_0^*$ of the vector space $\hD_0$ defined by
\begin{equation*}
\La^{\hD}(\la):=\sum_{{j\le r}\atop
j\in\hf\N}\xi_j\omega_j+d\La_0.
\end{equation*}

By using similar arguments as in \cite{CW3} (see also
\cite{W1}), we can explicitly construct the joint highest weight vectors of $\hD$ and
$SO(2d)$ appearing in \eqnref{evenduality}. The following theorem is an
easy consequence of the decomposition \eqnref{evenduality} and the description of the
joint highest weight vectors.

\begin{thm} \label{duality-hD}
The Lie superalgebra $\hD$ and $O(2d)$ form a dual pair on $\F_0$
in the sense of Howe, namely, we have the following
multiplicity-free decomposition of $\F_0$ with respect to their
joint action
\begin{equation*}
\F\cong\sum_{\la}L(\hD,\Lambda^{\hD}(\la))\otimes V_{O(2d)}^{\la},
\end{equation*}
where the summation is over all partitions of length $2d$ with
$\la'_1+\la'_2\le 2d$. Furthermore the $\hD\times\fso(2d)$ joint highest weight
vector of the $\la$-component can be described in the following way.
\begin{enumerate}
\renewcommand{\labelenumi}{(\roman{enumi})}
\item When $\la'_1< d$,
\begin{equation*}
X^1_{\la_1'} \cdot  X^2_{\la_2'}\cdots
X^{\la_1}_{{\la}_{\la_1}'}\vac
\end{equation*}
is a joint highest weight vector of $\hD\times\fso(2d)$ with the
joint highest weight $\La^{\hD}(\la)+\sum_{i=1}^d\la_i\epsilon_i\in
\hD^*_0\oplus\h^*$.
\item When $\la'_1= d$,
\begin{equation*}
X^1_{\la_1'} \cdot  X^2_{\la_2'}\cdots
X^{\la_1}_{{\la}_{\la_1}'}\vac
\end{equation*}
and
\begin{equation*}
\tX^1_{\la_1'} \cdot
\tX^2_{\la_2'}\cdots\tX^{\la_1}_{{\la}_{\la_1}'}\vac
\end{equation*}
are joint highest weight vectors of $\hD\times\fso(2d)$ with the
joint highest weights $\sum_{i=1}^d\la_i\epsilon_i+\La^{\hD}(\la)$
and
$\sum_{i=1}^{d-1}\la_i\epsilon_i-\la_d\epsilon_d+\La^{\hD}(\la)$,
respectively.
\item When  $\la'_1 > d$,
\begin{equation*}
\Gamma_{\la_1'} \cdot  X^2_{\la_2'}\cdots
X^{\la_1}_{{\la}_{\la_1}'}\vac
\end{equation*}
is a joint highest weight vector of $\hD\times\fso(2d)$ with the
joint highest weight
$\sum_{i=1}^{2d-\la'}\la_i\epsilon_i+\La^{\hD}(\la)$.
\end{enumerate}
\end{thm}

\begin{rem}\label{F0} By examining the embedding of $\hD$ in $\hA$
(also recalling Remark \ref{embeddings}), we can see that the
restriction of the anti-linear anti-involution $\omega$ on the
Lie superalgebra $\hA$ to $\hD$ gives an anti-linear
anti-involution also denoted by $\omega$ on $\hD$ such that
$\omega(C)=C$ and
$$\omega(\sum_{p\in \hf\Z^{^*}}a_p\te_{p-k,p}) =\sum_{p\in
\hf\Z^{^*}} (-1)^{[p]+[p-k]}\overline{a}_p\te_{p,p-k},
$$
for all $\sum_{p\in \hf\Z^{^*}}a_p\te_{p-k,p}\in \hD_k $,
$k\in \hf\Z$. Since the Fock space $\F_0$ is a unitarizable
$\hA$-module with the positive contravariant Hermitian form $\langle\
\cdot|\cdot\ \rangle$, it is also a unitarizable
$\hD$-module with respect to the anti-linear anti-involution
$\omega$. By \thmref{duality-hD}, for every partition $\lambda$ of
length $2d$ satisfying $\la_1'+\la_2'\le 2d$, the irreducible $\hD$-module
$L(\hD,\Lambda^{\hD}(\lambda))$ is unitarizable. \end{rem}

\subsubsection{Free field realization of $\hD$ and
$(\hD, O(2d+1))$-duality. }
Now we turn to the free field realization of $\hD$ associated with the
$(\hD, O(2d+1))$-duality. Introduce a free fermionic field
$\phi(z):=\sum_{n\in\Z^{^*}}\phi_nz^{-n-1}$ and a free bosonic
field $\chi(z):=\sum_{r\in\frac{1}{2}+\Z}\chi_rz^{-r-1/2}$ with
the non-trivial anti-commutation relations
$[\phi_m,\phi_n]=\delta_{ij}\delta_{m+n,0}$ and commutation
relations $[\chi_r,\chi_s]=\delta_{ij}\delta_{r+s,0}$ for $r>0$. We
shall denote by $\ccAhf_0$ the associative superalgebra generated by
the modes of all the quantum fields
$\widetilde{\psi}^{\pm,i}(z)$, $\gamma^{\pm,i}(z)$, $i=1,\cdots,d$,
$\phi(z)$, and $\chi(z)$. Let $\Fhf_0$ denote the Fock space of
the quantum fields generated by the vacuum
vector $\vac$, where
$\psi^{\pm,i}_m\vac=\gamma^{\pm,i}_r\vac=\phi_m\vac=\chi_r\vac=0$,
for $i=1,2\cdots,d$, $m>0$ and $r>0$.

Introduce an anti-linear anti-involution $\omega$ on $\ccAhf_0$ in the
following way.  It is defined by \eqref{conjugate-psi} and \eqref{conjugate-gamma}
on all the $\widetilde{\psi}^{\pm,i}_j$ and $\gamma^{\pm,i}_j$, and
\begin{align}
\omega(\phi_i) = \phi_{-i}, \quad \mbox{for all} \ i;  &&
\omega(\chi_r)=\chi_{-r}, \quad \mbox{for all}\ r. \label{conjugate-phichi}
\end{align}
The Fock space $\Fhf_0$ admits a positive definite contravariant
Hermitian form with respect to this $\ast$-structure of $\ccAhf_0$. As
usual, we shall normalize the form on the vacuum vector $\vac$
so that $\langle 0|0 \rangle=1$.

We have an action of $\hD$ of central charge $d+\hf$ on $\Fhf_0$
given by ($i,j\in\Z^*$ and $r,s\in\frac{1}{2}+\Z$)
\begin{align*}
 \te_{ij}&:=\sum_{p=1}^d:\psi^{+,p}_{-i}\psi^{-,p}_{j}:
 -\sum_{p=1}^d:\psi^{+,p}_{j}\psi^{-,p}_{-i}:+:\phi_{-i}\phi_j:;\\
 \te_{rs}&:=-\te_{-s,-r}:=-\sum_{p=1}^d:\gamma^{+,p}_{-r}\gamma^{-,p}_{s}:
 +\sum_{p=1}^d:\gamma^{+,p}_{s}\gamma^{-,p}_{-r}:+:\chi_{-r}\chi_s:,\; rs>0;\\
  \te_{rs}&:=-\sum_{p=1}^d:\gamma^{+,p}_{-r}\gamma^{-,p}_{s}:
 -\sum_{p=1}^d:\gamma^{+,p}_{s}\gamma^{-,p}_{-r}:+:\chi_{-r}\chi_s:,\, r<0,s>0;\\
   \te_{rs}&:=-\sum_{p=1}^d:\gamma^{+,p}_{-r}\gamma^{-,p}_{s}:
 -\sum_{p=1}^d:\gamma^{+,p}_{s}\gamma^{-,p}_{-r}:-:\chi_{-r}\chi_s:,\, r>0,s<0;\\
 \te_{is}&:=\te_{-s,-i}:=\sum_{p=1}^d:\psi^{+,p}_{-i}\gamma^{-,p}_{s}:
 -\sum_{p=1}^d:\gamma^{+,p}_{s}\psi^{-,p}_{-i}:+:\phi_{-i}\chi_s:, \; s>0;\\
  \te_{is}&:=-\te_{-s,-i}:=\sum_{p=1}^d:\psi^{+,p}_{-i}\gamma^{-,p}_{s}:
 +\sum_{p=1}^d:\gamma^{+,p}_{s}\psi^{-,p}_{-i}:-:\phi_{-i}\chi_s:,\; s<0.
\end{align*}

\begin{rem}\label{Fhf}
It is important to observe that under the anti-linear anti-involution
$\omega$ of $\ccAhf_0$, the operators defined in the above equations
transform as follows:
$$\omega(\te_{p q}) =(-1)^{[p]+[q]}\te_{q p}, \quad \mbox{for all} \ p,
q. $$
Therefore the realization of $\hD^f$ in $\ccAhf_0$ given above defines
a $\ast$-superalgebra homomorphism from ${\mc U}(\hD^f)$ to $\ccAhf_0$.
Note that the $\ast$-structure on $\hD^f$ is the restriction of
that described in Remark \ref{F0}.
\end{rem}

Introduce the following operators on the Fock space $\Fhf_0$:
\begin{align}
\label{so++} E^{so+}_{i}&=\sum_{n\in\Z^*}:\phi_{-n}\psi^{+,i}_{n}:
-\sum_{r\in1/2+\Z_+}:\chi_{-r}\gamma^{+,i}_{r}:
+\sum_{r\in-1/2-\Z_+}:\chi_{-r}\gamma^{+,i}_{r}:,\\
\label{so--} E^{so-}_{j}&=\sum_{
n\in\Z^*}:\phi_{-n}\psi^{-,j}_{n}: +\sum_{r\in
1/2+\Z^*}:\chi_{-r}\gamma^{-,j}_{r}:,
\end{align}
where $1\le i,j\le d$. Note that \eqnref{Eij}, \eqnref{so+} and
\eqnref{so-} can be extended to actions on the Fock space
$\Fhf_0$. It is easy to see that \eqnref{Eij}, \eqnref{so+} and
\eqnref{so-} together with \eqnref{so++} and \eqnref{so--} form a
basis for the Lie algebra $\fso (2d+1)$. The action of the Lie
algebra $\fso(2d+1)$ on the Fock space $\Fhf_0$ can be lifted to
an action of the Lie group $SO(2d+1)$, which can further be extended
to an action of the Lie group $O(2d+1)$. Moreover $\Fhf_0$ is a direct sum of
finite dimensional irreducible representations of $O(2d+1)$.

Using similar arguments as before, we can show that the Lie superalgebra $\hD$
and $O(2d+1)$ form a dual pair on $\Fhf_0$ in the sense of Howe
\cite{H}. We have the following multiplicity-free decomposition of
$\Fhf_0$ with respect to the joint action of $\hD$ and $O(2d+1)$:
\begin{equation}\label{Oddduality}
\Fhf_0\cong\sum_{\la}W^\la\otimes V_{O(2d+1)}^{\la},
\end{equation}
where the summation is over a subset of all partitions of length
with $\la'_1+\la'_2\le 2d+1$, and $W^\la$ is a certain irreducible
module over $\hD$.

We define the $(2d+1)\times (2d+1)$ matrix $\widetilde\Gamma$ as
follows:
\begin{eqnarray*}
 \widetilde\Gamma :=
\begin{pmatrix}
 \gamma_{-\hf}^{+,1}&\gamma_{-\hf}^{+,2}&\cdots&\gamma_{-\hf}^{+,d}
 &\chi_{-\hf}&\gamma_{-\hf}^{-,d}&\gamma_{-\hf}^{-,d-1}&\cdots&\gamma_{-\hf}^{-,1}\\
 \psi_{-1}^{+,1}&\psi_{-1}^{+,2}&\cdots&\psi_{-1}^{+,d}&\phi_{-1}
 &\psi_{-1}^{-,d}&\psi_{-1}^{-,d-1}&\cdots&\psi_{-1}^{-,1}\\
 \psi_{-1}^{+,1}&\psi_{-1}^{+,2}&\cdots&\psi_{-1}^{+,d}&\phi_{-1}
 &\psi_{-1}^{-,d}&\psi_{-1}^{-,d-1}&\cdots&\psi_{-1}^{-,1}\\
 \vdots&\vdots&\cdots&\vdots&\vdots&\vdots&\vdots&\cdots &\vdots\\
 \psi_{-1}^{+,1}&\psi_{-1}^{+,2}&\cdots&\psi_{-1}^{+,d}&\phi_{-1}
 &\psi_{-1}^{-,d}&\psi_{-1}^{-,d-1}&\cdots&\psi_{-1}^{-,1}\\
\end{pmatrix}.
\end{eqnarray*}
For any nonnegative inetger $r$, we let $\widetilde\Gamma_r$
denote the first $r\times r$ minor of the matrix
$\widetilde\Gamma$.

Given a partition $\la=(\la_1,\la_2,\cdots,\la_{2d+1})$ of length
$2d+1$ with $\la'_1+\la'_2\le 2d+1$, we have the shifted Frobenius
notations
$\La(\la)=({\xi_{\frac{1}{2}},\xi_{1\frac{1}{2}},\cdots,\xi_{r-\frac{1}{2}}}
\mid {\xi_1,\xi_2,\cdots,\xi_r})$ for the partition $\la$ (see
\secref{Frobenius}). Let $\La^{\hD}(\la)$ be an element of the
dual space $\hD_0^*$ of the vector space $\hD_0$ defined by
\begin{equation}
\La^{\hD}(\la):=\sum_{{j\le r}\atop
j\in\hf\N}\xi_j\omega_j+\frac{2d+1}{2}\La_0.
\end{equation}

The following theorem is a consequence of the decomposition \eqnref{Oddduality}
and the explicit description of the $\hD\times\fso(2d+1)$ joint
highest weight vectors.
\begin{thm} \label{duality-hD-hf}
The Lie superalgebra $\hD$ and $O(2d+1)$ form a dual pair on
$\Fhf_0$ in the sense of Howe. In particular, we have the following
(multiplicity-free) decomposition of $\Fhf_0$ with respect to
their joint action
\begin{equation*}
\Fhf_0\cong\sum_{\la}L(\hD,\Lambda^{\hD}(\la))\otimes
V_{O(2d+1)}^{\la},
\end{equation*}
where the summation is over all partitions of length $2d+1$ with
$\la'_1+\la'_2\le 2d+1$. Furthermore, the joint highest weight
vector of the $\la$-component with respect to $\hD\times SO(2d+1)$
is given
\begin{equation*}
\Gamma_{\la_1'} \cdot  X^2_{\la_2'}\cdots
X^{\la_1}_{{\la}_{\la_1}'}\vac.
\end{equation*}
\end{thm}
\begin{rem}\label{Fhf-unitary}
By recalling Remark \ref{Fhf}, we can easily show that
all the irreducible $\hD$-modules appearing in the above theorem
are unitarizable with respect to the the anti-linear anti-involution $\omega$ described
in Remark \ref{F0}, that is, every irreducible module $L(\hD,\Lambda^{\hD}(\la))$ associated
with a partition $\lambda$ of length $2d+1$ satisfying $\la'_1+\la'_2\le 2d+1$
is unitarizable.
\end{rem}

\subsubsection{Unitarizable $\hD$-modules} Now we classify the
unitarizable irreducible quasi-finite highest weight $\hD$-modules
with respect to the the anti-linear anti-involution $\omega$ described
in Remark \ref{F0}. We have the following result.

\begin{thm}\label{unitary-hD}
Let $M$ be an irreducible quasi-finite highest weight $\hD$-module
with highest weight $\xi$. Then $M$ is unitarizable if and
only if there exists a non-negative integer or half integer $k$ such that
$\xi=\Lambda^{\hD}(\lambda)$ for some partition $\lambda$ of
length $2k$ with $\la'_1+\la'_2\le 2k$. In other words, $M$ is
unitarizable if and only if
\begin{equation*}
\La^{\hD}(\la):=\sum_{{j\le r}\atop j\in\hf\N}\xi_j\omega_j+k\La_0
\end{equation*}
such that $k\in \hf\Z_+$, $r\in\N$ and $\xi_j\in\Z$ for all $j$
satisfying the following conditions:
\begin{enumerate}
\renewcommand{\labelenumi}{(\roman{enumi})}
\item $\xi_{\hf}>\xi_{\frac{3}{2}}>\cdots>\xi_{r-\hf} \geq 0$,
$\xi_{1}>\xi_{2}>\cdots>\xi_{r} \geq 0$, and $\xi_{r-\hf}=0$ if
and only if $ r=1$ and $ \xi_\hf=\xi_1=0$, \item
$\xi_1+\xi_2+l_{1,2}(\xi_\hf)+\min\{\xi_{\frac {3}{2}},1\} \leq
2k$,
\end{enumerate}
where $l_{1,2}$ is a function from non-negative integers to itself
with $l_{1,2}(0)=0$, $l_{1,2}(1)=1$ and $l_{1,2}(x)=2$ if $x\ge
2$.
\end{thm}

\begin{proof} We have already pointed out in Remarks \ref{F0} and \ref{Fhf-unitary}
that the $L(\hD,\Lambda^{\hD}(\lambda))$ are unitarizable
irreducible quasi-finite highest weight $\hD$-modules for all
partitions $\lambda$. Now we are going to show that if $M$ is a
unitarizable irreducible quasi-finite highest weight $\hD$-module
with the highest weight $\xi$, then $\xi=\Lambda^{\hD}(\lambda)$
for some partition $\lambda$. Let $\langle \ \cdot | \cdot \
\rangle$ be a positive definite contravariant Hermitian form on
$M$ and $v_\xi$ a highest weight vector of $M$ such that $ \langle
v_\xi | v_\xi\rangle =1$. We put $\xi(\te_{i,i})=\xi_i$ for all
$i\in \hf\N$. By \thmref{criterion-hD}, there exists $r\in \N$
such that $\xi=\sum_{{j\le r}\atop
j\in\hf\N}\xi_j\omega_j+k\La_0$, where $\xi_j$, $k\in \C$. By
using similar arguments as in the proof of \thmref{unitary}, we can
show that
$\xi_i\in \Z_+$ for $i=\hf,1,\cdots,r-\hf,r$, and
$$\xi_{\hf}>\xi_{1\hf}>\cdots>\xi_{r-\hf} \geq 0,\qquad
\xi_{1}>\xi_{2}>\cdots>\xi_{r} \geq 0.$$ Moreover,
$\xi_{r-\hf}=0$ if and only if $ r=1$ and $ \xi_\hf=\xi_1=0$.

Now we choose a large positive integer $n$ such that
$\xi(\te_{n,n})=\xi(\te_{n+1,n+1})=0$. Consider the subalgebra
$sl(2,\C)$ with the standard basis $\{ 2C-\te_{n,n}-\te_{n+1,n+1}$,
$\te_{-n-1,n}$, $\te_{n,-n-1}\}$. Note that
$\omega(2C-\te_{n,n})-\te_{n+1,n+1}=2C-\te_{n,n}-\te_{n+1,n+1}$,
$\omega(\te_{n,-n-1})=\te_{-n-1,n}$ and
$\omega(-\te_{-n-1,n})=-\te_{n,-n-1}$. Unitarizability with respect to this
$sl_2$ subalgebra requires
$2k=\xi(2C)=\xi(2C-\te_{n,n}-\te_{n+1,n+1})\in\Z_+$. Hence we have
$k\in \hf\Z_+$.

Finally, we need to show
$\xi_1+\xi_2+l_{1,2}(\xi_\hf)+\min\{\xi_{\frac {3}{2}},1\} \leq
2k$. We may assume that $\xi_\hf >0$. Otherwise, we have
$\xi_1+\xi_2+l_{1,2}(\xi_\hf)+\min\{\xi_{\frac {3}{2}},1\}=0\le
2k$. Direct computations show  that $\langle \te_{1,\hf}v_\xi\mid
\te_{1,\hf}v_\xi\rangle=\xi_\hf+\xi_1> 0$ and
\begin{equation}\label{hDeqn1}
\begin{array}{ll}
0\le\| \te_{2,-1}\te_{1,\hf}v_\xi\|^2
  &=\langle\te_{1,\hf}v_\xi\mid
  \te_{-1,2}\te_{2,-1}\te_{1,\hf}v_\xi\rangle\\
  &=\langle\te_{1,\hf}v_\xi\mid
  (2C-\te_{1,1}-\te_{2,2})\te_{1,\hf}v_\xi\rangle\\
  &=(2k-\xi_1-\xi_2-1)(\xi_\hf+\xi_1).
\end{array}
\end{equation}
Thus we have $2k-\xi_1-\xi_2-1\ge 0$ since $\xi_\hf+\xi_1> 0$.
When $\xi_\hf=1$, we have $\xi_{\frac{3}{2}}=\xi_2=0$. Thus
$\xi_1+\xi_2+l_{1,2}(\xi_\hf)+\min\{\xi_{\frac {3}{2}},1\}
=\xi_1+\xi_2+1+0 \leq 2k$. Now we assume that $\xi_\hf\ge 2$. We
compute
\begin{equation}\label{hDeqn2}
\begin{array}{ll}
 \| \te_{-\hf,1}\te_{2,-1}\te_{1,\hf}v_\xi\|^2
  &=\|\te_{2,\hf}\te_{1,\hf}v_\xi\|^2\\
  &=\langle\te_{1,\hf}v_\xi\mid
  \te_{\hf,2}\te_{2,\hf}\te_{1,\hf}v_\xi\rangle\\
  &=\langle\te_{1,\hf}v_\xi\mid
  (\te_{\hf,\hf}+\te_{2,2})\te_{1,\hf}v_\xi\rangle\\
  &=(\xi_\hf+\xi_2-1)(\xi_\hf+\xi_1)\\
  &>0.
\end{array}
\end{equation}
Thus $\| \te_{2,-1}\te_{1,\hf}v_\xi\|^2>0$ and hence we have
$\xi_1+\xi_2+l_{1,2}(\xi_\hf)+\min\{\xi_{\frac
{3}{2}},1\}=\xi_1+\xi_2+2\leq 2k$ for $\xi_\hf\ge 2$ and
$\xi_{\frac {3}{2}}=0$ by using \eqnref{hDeqn1}. Eventually we
assume that $\xi_\hf\ge 2$ and $\xi_{\frac {3}{2}}>0$. Then
\begin{equation*}
\begin{array}{ll}
 \| \te_{-\hf,1}\te_{2,-1}\te_{2,\frac{3}{2}}\te_{1,\hf}v_\xi\|^2
  &=\|\te_{2,\hf}\te_{2,\frac{3}{2}}\te_{1,\hf}v_\xi\|^2\\
  &=\langle\te_{2,\frac{3}{2}}\te_{1,\hf}v_\xi\mid
  \te_{\hf,2}\te_{2,\hf}\te_{2,\frac{3}{2}}\te_{1,\hf}v_\xi\rangle\\
  &=\langle\te_{2,\frac{3}{2}}\te_{1,\hf}v_\xi\mid
  (\te_{\hf,\hf}+\te_{2,2})\te_{2,\frac{3}{2}}\te_{1,\hf}v_\xi\rangle\\
  &=(\xi_\hf+\xi_2)\|\te_{2,\frac{3}{2}}\te_{1,\hf}v_\xi\|^2\\
  &=(\xi_\hf+\xi_2)(\xi_\hf+\xi_1)(\xi_\hf+\xi_2)\\
  &>0.
\end{array}
\end{equation*}
Therefore we have
$\te_{2,-1}\te_{2,\frac{3}{2}}\te_{1,\hf}v_\xi\not=0$ and
\begin{equation*}
\begin{array}{ll}
0<\|\te_{2,-1}\te_{2,\frac{3}{2}}\te_{1,\hf}v_\xi\|^2
    &=\langle\te_{2,\frac{3}{2}}\te_{1,\hf}v_\xi\mid
  \te_{-1,2}\te_{2,-1}\te_{2,\frac{3}{2}}\te_{1,\hf}v_\xi\rangle\\
  &=\langle\te_{2,\frac{3}{2}}\te_{1,\hf}v_\xi\mid
  (2C-\te_{1,1}-\te_{2,2})\te_{2,\frac{3}{2}}\te_{1,\hf}v_\xi\rangle\\
  &=(2k-\xi_\hf-\xi_2-2)\|\te_{2,\frac{3}{2}}\te_{1,\hf}v_\xi\|^2\\
  &=(2k-\xi_\hf-\xi_2-2)(\xi_\hf+\xi_1)(\xi_\hf+\xi_2).
\end{array}
\end{equation*}
Therefore we have $\xi_\hf+\xi_2+3\le 2k$. Hence
$\xi_1+\xi_2+l_{1,2}(\xi_\hf)+\min\{\xi_{\frac
{3}{2}},1\}=\xi_\hf+\xi_2+3\leq 2k$ and the proof is completed.
\end{proof}

\section{Character formulas for unitarizable irreducible modules}
\label{sect character}
In this section we derive explicit formulae for the formal
characters of the unitarizable quasi-finite irreducible highest
weight modules over the infinite rank Lie algebras
$\hC$ and $\hD$. The method employed here is a generalization of
that developed in \cite{CL, CZ, CLZ}, which makes essential use of Howe
dualities. We mention that the character formulae for the unitarizable
irreducible modules over $\hgltwo$ (and hence $\hA$) were obtained
in \cite{CL}.

\subsection{Character formula for $\hC$} \label{character}
The central result of this subsection is \thmref{characterthm}, which
gives the character formula for the unitarizable quasi-finite irreducible highest
weight $\hC$-modules.  In order to establish the result, we need some
basic facts on characters of the symplectic group (see \cite{FH},
\cite{H1}, \cite{H2}), which we recall here.

When we deal with characters of modules over $\fsp(2m)$, we put ${\tilde
h}_{i}:=-h_{m-i+1}$ and $x_i=e^{-\epsilon_{m-i+1}}$ for
$i=1,2,\cdots,m$. That is $x_i=z^{-1}_{m-i+1}$. Recall that
$\epsilon_i( h_j)=\delta_{ij}$ and $z_i=e^{\epsilon_i}$ where
$\epsilon_i\in\h^*$ such that $\epsilon_i(h_j)=\delta_{ij}$ (see
\secref{unitarizable-hC}). The definitions of $x_i$'s and ${\tilde
h}_i$'s are somewhat nonstandard,  but they allow us to
deal with only polynomials instead of Laurent polynomials when
considering characters of certain representations of $\fsp(2m)$.
For each finite sequence of
complex numbers $\la=(\la_1,\cdots,\la_m )$, we let
$W^\la_{\fsp(2m)}:=V^{\la^*}_{\fsp(2m)}$ where
$\la^*=(-\la_m,\cdots,-\la_1)$. Note that $\la^*({\tilde
h}_{i})=\la_i$ for $i=1,\cdots,m$. $W^\la_{\fsp(2m)}$ is a
finite-dimensional irreducible representation if and only if
$\la_1\ge\la_2\ge\cdots\ge\la_m$ and $\la_i\in -\Z_+$ for
$i=1,\cdots,m$.

When we deal with characters of $Sp(2d)$-modules,
we put $z_i=e^{\epsilon_i}$ where $\epsilon_i\in\h^*$ such that
$\epsilon_i(h_{j})=\epsilon_i(e_{jj}-e_{d+j,d+j})=\delta_{ij}$
(see \secref{unitarizable-hC}). For each partition $\la$ of length
$d$ and each decreasing sequence of non-positive integers $\nu$ of
length $m$, we write
$\chi^{\la}_{Sp(2d)}(\z)=\chi^{\la}_{Sp(2d)}(z_1,\cdots,z_d)$ and
$\tilde{\chi}^{\nu}_{\fsp(2m)}(\x)=\tilde{\chi}^{\nu}_{\fsp(2m)}(x_1,\cdots,x_m)$
for the characters of the $Sp(2d)$-module $V^\la_{Sp(2d)}$ and
the $\fsp(2m)$-module $W^\nu_{\fsp(2m)}$,  to stress their dependence on
the variables $z_1,\cdots,z_d$ and $x_1,\cdots,x_m$, respectively.
It is clear that
$\chi^{\la}_{Sp(2m)}(\z)=\tilde{\chi}^{\la^*}_{\fsp(2m)}(z_1^{-1},\cdots,z_m^{-1})$
for any partition $\la$ of length $m$ since
$\chi^{\la}_{\fsp(2m)}(z_1,\cdots,z_m)
=\chi^{\la^*}_{Sp(2m)}(z_1^{-1},\cdots,z_m^{-1})$.

By the $(Sp(2d), \fsp(2m))$-duality on the
exterior algebra $\La(\C^{2d}\otimes \C^{m})$ with $m\ge d$
(\cite{H1}, \cite{H2}, also see \cite{CZ}), we have the following
identity:
\begin{align}
&(x_1\cdots x_m)^{-d}
\prod_{i=1}^{d}\prod_{j=1}^m(1+x_jz_i)(1+x_jz_i^{-1})
=\sum_{\la}\chi^{\la}_{Sp(2d)}(\z) {\tilde\chi}^{{\la'}-d (1^m)}_{\fsp(2m)}(\x).\label{spchar}
\end{align}
Here $\la$ is summed over all partitions of length $d$ with
$\la_1\le m$, and for any $k\in \N$,  $ (1^k)$ stands for the $k$-tuple
$(1,1,\cdots,1)$. Note that the partition $\la'$ is considered as
a partition of length $m$ and ${\la'}-d (1^m)$ is a decreasing
sequence of non-positive integers of length $m$.

We let
\begin{align*}
E'_r(x_1,\cdots,x_m,x^{-1}_1,\cdots,x_m^{-1}) =
&E_r(x_1,\cdots,x_m,x^{-1}_1,\cdots,x_m^{-1})\\
& - E_{r-2}(x_1,\cdots,x_m,x^{-1}_1,\cdots,x_m^{-1}),
\end{align*}
where $E_r$ is the $r$-th elementary symmetric polynomial for
$r\ge 0$ and $E_r=0$ for $r<0$. For any partition $\mu$ of length $l$
(which is an arbitrary positive integer), we define
$$|E'_\mu|=|E'_\mu(x_1,\cdots,x_m,x^{-1}_1,\cdots,x_m^{-1})|$$
by the the determinant of the $l\times l$ matrix
\begin{equation}\label{E'}
\begin{pmatrix}
E'_{\mu_1}&E'_{\mu_1+1}+E'_{\mu_1-1}&E'_{\mu_1+2}+E'_{\mu_1-2}
&\cdots&E'_{\mu_1+l-1}+E'_{\mu_1-l+1} \\
E'_{\mu_2-1}&E'_{\mu_2}+E'_{\mu_2-2}&E'_{\mu_2+1}+E'_{\mu_2-3}
&\cdots&E'_{\mu_2+l-2}+E'_{\mu_2-l} \\
\vdots&\vdots&\vdots&\cdots&\vdots \\
E'_{\mu_i-i+1}&E'_{\mu_i-i+2}+E'_{\mu_i-i}&E'_{\mu_i-i+3}+E'_{\mu_i-i-1}
&\cdots&E'_{\mu_i-i+l}+E'_{\mu_i-i-l+2} \\
\vdots&\vdots&\vdots&\cdots&\vdots \\
E'_{\mu_l-l+1}&E'_{\mu_l-l+2}+E'_{\mu_l-l}&E'_{\mu_l-l+3}+E'_{\mu_l-l-1}
&\cdots&E'_{\mu_l}+E'_{\mu_l-2l+2} \\
\end{pmatrix}.
\end{equation}
Then the character of the irreducible
$\fsp(2m)$-module $V^\la_{\fsp(2m)}$ is given \cite{FH} by
$|E'_{\la'}|=|E'_{\la'}(x_1,\cdots,x_m,x^{-1}_1,\cdots,x_m^{-1})|$,
where  $\la$ is any partition  of length $m$.

From the above result it is not very difficult to work out that
the character $\tilde{\chi}^{{\la'}-d (1^m)}_{\fsp(2m)}(\x)$ of the
finite-dimensional irreducible $\fsp(2m)$-module
$W^{{\la'}-d (1^m)}_{\fsp(2m)}$ for every partition $\la$ of length $d$ with $\la_1\le m$
is given by the determinant of the $d\times d$ matrix whose $i$-th row is
\begin{align*}
(E'_{m-\la_{d-i+1}-i+1}\quad
E'_{m-\la_{d-i+1}-i+2}+E'_{m-\la_{d-i+1}-i} \quad
E'_{m-\la_{d-i+1}-i+3}+E'_{m-\la_{d-i+1}-i-1} \allowdisplaybreaks\\
 \quad \cdots\quad
E'_{m-\la_{d-i+1}-i+d}+E'_{m-\la_{d-i+1}-i-d+2}).
\end{align*}

On the other hand, for each $r\in \Z$,
 \begin{align*}
 &E_{m-r}(x_1,\cdots,x_m,x^{-1}_1,\cdots,x_m^{-1})\\
 =&\sum_{i=0}^{m-r}E_{i}(x_1,\cdots,x_m)E_{m-r-i}(x^{-1}_1,\cdots,x_m^{-1})\\
 =&(x_1\cdots x_m)^{-1}\sum_{i=0}^{m-r}E_{i}(x_1,\cdots,x_m)
 E_{r+i}(x_1,\cdots,x_m)\\
 =&(x_1\cdots x_m)^{-1}\sum_{i=0}^{\infty}E_{i}(x_1,\cdots,x_m)
 E_{r+i}(x_1,\cdots,x_m).
 \end{align*}
 For $r\in \Z$, we define
 \begin{align*}
 &\tE_r:=\sum_{i=0}^{\infty}E_{i}(x_1,\cdots,x_m)
 E_{r+i}(x_1,\cdots,x_m), \\
&\tE'_r:=\tE_r-\tE_{r-2}.
 \end{align*}
Then we have
 \begin{equation*}
 E'_{m-r}(x_1,\cdots,x_m,x^{-1}_1,\cdots,x_m^{-1})
 = (x_1\cdots x_m)^{-1}\tE'_r.
 \end{equation*}
Hence for each partition $\la$ of length $d$, $(x_1\cdots
x_m)^{d}\tilde{\chi}^{{\la'}-d (1^m)}_{\fsp(2m)}(\x)$ equals the
determinant of the $d\times d$ matrix whose $i$-th row is
\begin{align*}
(\tE'_{\la_{d-i+1}+i-1}\quad
\tE'_{\la_{d-i+1}+i-2}+\tE'_{\la_{d-i+1}+i} \quad
\tE'_{\la_{d-i+1}+i-3}+\tE'_{\la_{d-i+1}+i+1} \allowdisplaybreaks\\
 \quad \cdots\quad
\tE'_{\la_{d-i+1}+i-d}+\tE'_{\la_{d-i+1}+i+d-2}).
\end{align*}

For each partition $\la$ of length $d$ satisfying $\la_1\le m$ and
$d\le m$, the {\em symplectic Schur polynomial of weight $d$} of
$m$ variables, denoted by
$$S_\la^{\fsp,d}(\x)=S_\la^{\fsp,d}(x_1,\cdots,x_m),$$
is defined by the determinant of the $d\times d$
matrix whose $i$-th row is
\begin{align*}
(\tE'_{\la_{d-i+1}+i-1}\quad
\tE'_{\la_{d-i+1}+i-2}+\tE'_{\la_{d-i+1}+i} \quad
\tE'_{\la_{d-i+1}+i-3}+\tE'_{\la_{d-i+1}+i+1} \allowdisplaybreaks\\
 \quad \cdots\quad
\tE'_{\la_{d-i+1}+i-d}+\tE'_{\la_{d-i+1}+i+d-2}).
\end{align*}
Thus we have
 \begin{equation*}
 {\tilde\chi}^{{\la'}-d (1^m)}_{\fsp(2m)}(x_1,\cdots,x_m) \\
 =(x_1\cdots x_m)^{-d}S_\la^{\fsp,d}(x_1,\cdots,x_m).
\end{equation*}
Hence we can rewrite the combinatorial formula \eqnref{spchar} as
follows:
\begin{align}\label{combin-Sp}
 \prod_{i=1}^{d}\prod_{j=1}^m(1+x_jz_i)(1+x_jz_i^{-1})
=\sum_{\la}\chi^{\la}_{Sp(2d)}(\z)S_\la^{\fsp,d}(\x).
\end{align}
Here $\la$ is summed over all partitions of length $d$ with
$\la_1\le m$.

Now, for every partition $\la$ of length $d$, we still use
$$S_\la^{\fsp,d}(\x)=S_\la^{\fsp,d}(x_1,\cdots,x_m,\cdots)$$
to denote the {\em symplectic Schur function of weight $d$}
of infinitely many variables, which is the determinant of the
$d\times d$ matrix whose $i$-th row is
\begin{align*}
(\te'_{\la_{d-i+1}+i-1}\quad
\te'_{\la_{d-i+1}+i-2}+\te'_{\la_{d-i+1}+i} \quad
\te'_{\la_{d-i+1}+i-3}+\te'_{\la_{d-i+1}+i+1} \allowdisplaybreaks\\
 \quad \cdots\quad
\te'_{\la_{d-i+1}+i-d}+\te'_{\la_{d-i+1}+i+d-2}),
\end{align*}
where $ \te'_r=\te_r-\te_{r-2}$ and
$\te_r=\sum_{i=0}^{\infty}e_{i}(x_1,x_2,\cdots) e_{r+i}(x_1,x_2,
\cdots)$. Hereafter $e_{i}(x_1,x_2,\cdots)$ stands for the
i-th elementary symmetric function of infinitely many variables. Therefore,
the symplectic Schur function $S_\la^{\fsp,d}(x_1,x_2,\cdots)$ is
the inverse limit of the symplectic Schur polynomials
$S_\la^{\fsp,d}(x_1,\cdots,x_m)$,  and we thus have
$S_\la^{\fsp,d}(x_1,\cdots,x_m)$ $=$
$S_\la^{\fsp,d}(x_1,\cdots,x_m,0,0,\cdots)$.

For any given partition $\la$ of length $d$, we let
$$DS_\la^{\fsp,d}(\x) =DS_\la^{\fsp,d}(x_1,x_2,\cdots)$$
denote the {\em skew symplectic Schur function of weight $d$}
of infinitely many variables, which is
defined by the determinant of the $d\times d$ matrix
whose $i$-th row is
\begin{align*}
({\tilde H}'_{\la_{d-i+1}+i-1}\quad {\tilde
H}'_{\la_{d-i+1}+i-2}+{\tilde H}'_{\la_{d-i+1}+i} \quad
{\tilde H}'_{\la_{d-i+1}+i-3}+{\tilde H}'_{\la_{d-i+1}+i+1} \allowdisplaybreaks\\
 \quad \cdots\quad
{\tilde H}'_{\la_{d-i+1}+i-d}+{\tilde H}'_{\la_{d-i+1}+i+d-2}),
\end{align*}
where
\begin{align*}
&{\tilde H}'_r:={\tilde H}_r-{\tilde H}_{r-2},\\
&{\tilde
H}_r:=\sum_{i=0}^{\infty}H_{i}(x_1,x_2,\cdots)
H_{r+i}(x_1,x_2,\cdots),
\end{align*}
and $H_{i}(x_1,x_2,\cdots)$ are the
complete symmetric functions of infinitely many variables.
Note that $(x_1\cdots x_m)^{-d}
DS_\la^{\fsp,d}(x_1^{-1},\cdots,x_m^{-1},0,0,\cdots)$ is the
character of some infinite dimensional unitarizable module over
the Lie algebra $\mathfrak{so}(2m)$.

Analogous to hook Schur functions (see \cite{BR} and \cite{CL}),
for each partition $\la$ of length $d$, we define the {\em symplectic
hook Schur function} $$HS_\la^{\fsp,d}(\x,\y)
=S_\la^{\fsp,d}(x_1,x_2,\cdots,y_\hf,y_{\frac{3}{2}},\cdots)$$ {\em of weight $d$}
in infinitely many variables by
 $$
 HS_\la^{\fsp,d}(\x,\y):=\sigma(S_\la^{\fsp,d}(\x,\y)),
 $$
where $\sigma$ is the involution of the ring of symmetric
functions (see for example \cite{M}), which sends the elementary
symmetric functions of $y_j$'s to the complete symmetric functions
of $y_j$'s. Recall that the hook Schur function (cf. \cite{CL})
$$HS_{\la}(\x,\y)=\sigma(S_{\la}(\x,\y))$$  is a
symmetric function of the variables $\bf x$ and the variables $\bf
y$ separably where $\la$ is a partition. Hereafter $S_\la$
stands for the Schur function associated with the partition $\la$.
Then we have
 $$
 HS_{(1^d)}(\x,\y)=\sigma(e_d(\x,\y))
 =\sum_{i=0}^{d}e_{i}(\x)e_{d-i}(\y).
 $$
The symplectic hook Schur function can be written in terms of hook
Schur functions as follows.

\begin{prop} Let $\x=\{x_1,x_2,\cdots\}$ and
${\bf y}=\{y_1,y_2,\cdots\}$ be two infinite sets of  variables. For
each partition $\la$ of length $d$, the symplectic hook Schur
function $HS_\la^{\fsp,d}(\x,\y)$ of weight $d$ equals the
determinant of the
 following $d\times d$ matrix whose $i$-th row is
\begin{align*}
(\tilde{f}'_{\la_{d-i+1}+i-1}\quad
\tilde{f}'_{\la_{d-i+1}+i-2}+\tilde{f}'_{\la_{d-i+1}+i} \quad
\tilde{f}'_{\la_{d-i+1}+i-3}+\tilde{f}'_{\la_{d-i+1}+i+1} \allowdisplaybreaks\\
 \quad \cdots\quad
\tilde{f}'_{\la_{d-i+1}+i-d}+\tilde{f}'_{\la_{d-i+1}+i+d-2}),
\end{align*}
where $ \tilde{f}'_r=\tilde{f}_r-\tilde{f}_{r-2}$ and
$\tilde{f}_r=\sum_{i=0}^{\infty}HS_{(1^i)}(\x,\y)HS_{(1^{r+i})}(\x,\y)$.
\end{prop}

\begin{proof}
Let $\sigma$ denote the involution of the ring of symmetric
functions, which sends the elementary symmetric functions of
variables $\bf y$ to the complete symmetric functions of variables
$\bf y$. The proposition follows by applying the involution
$\sigma$ to the determinant of the
 following $d\times d$ matrix whose $i$-th row is
\begin{align*}
(\te'_{\la_{d-i+1}+i-1}\quad
\te'_{\la_{d-i+1}+i-2}+\te'_{\la_{d-i+1}+i} \quad
\te'_{\la_{d-i+1}+i-3}+\te'_{\la_{d-i+1}+i+1} \allowdisplaybreaks\\
 \quad \cdots\quad
\te'_{\la_{d-i+1}+i-d}+\te'_{\la_{d-i+1}+i+d-2}),
\end{align*}
where $ \te'_r=\te_r-\te_{r-2}$, $\te_r=\sum_{i=0}^{\infty}e_{i}
e_{r+i}$ and $e_{i}=e_{i}({\bf x},{\bf y})$ are the elementary
symmetric functions of infinitely many variables ${\bf x}$ and ${\bf y}$.
\end{proof}

\begin{prop}\label{combin1} Let $\x=\{x_1,x_2,\cdots\}$
be an infinite set of  variables and ${\bf
z}=\{z_1,z_2,\cdots,z_d\}$ be $d$ variables. Then
\begin{align}
 \prod_{i=1}^{d}\prod_{j=1}^{\infty}(1+x_jz_i)(1+x_jz_i^{-1})
=\sum_{\la}\chi^{\la}_{Sp(2d)}(\z)S_\la^{\fsp,d}(\x),
\end{align}
where $\la$ is summed over all partitions of length $d$,
 and
\begin{align}
 \prod_{i=1}^{d}\prod_{j=1}^{\infty}(1-x_jz_i)^{-1}(1-x_jz_i^{-1})^{-1}
=\sum_{\la}\chi^{\la}_{Sp(2d)}(\z)DS_\la^{\fsp,d}(\x),
\end{align}
 where $\la$ is summed over all partitions of length $d$.
\end{prop}

\begin{proof} The first identity follows from \eqnref{combin-Sp}
by putting $m \rightarrow \infty$. Let $\sigma$ denote the
involution of the ring of symmetric functions, which sends the
elementary symmetric functions of $x_j$'s to the complete
symmetric functions of $x_j$'s. By applying the involution
$\sigma$ to both sides of the first identity of the
proposition, we obtain the second identity.
\end{proof}

\begin{prop}\label{HS} Let $\x=\{x_1,x_2,\cdots\}$ and
$\y=\{y_\hf,y_{\frac{3}{2}},\cdots\}$
be two infinite sets of variables and ${\bf
z}=\{z_1,z_2,\cdots,z_d\}$ be $d$ variables. Then
\begin{align*}
 \prod_{i=1}^{d}\prod_{j=1}^{\infty}\prod_{k=\hf}^{\infty}
 {\frac{(1+x_jz_i)(1+x_jz_i^{-1})}{(1-y_kz_i)(1-y_kz_i^{-1})}}
=\sum_{\la}\chi^{\la}_{Sp(2d)}(\z)HS_\la^{\fsp,d}(\x,\y),
\end{align*}
 where $\la$ is summed over all partitions of length $d$.
\end{prop}

\begin{proof}By \propref{combin1}, we have
\begin{align}
 &\prod_{i=1}^{d}\prod_{j=1}^{\infty}\prod_{k=\hf}^{\infty}
{(1+x_jz_i)(1+x_jz_i^{-1})}{(1+y_kz_i)(1+y_kz_i^{-1})}\label{eqnSp1}\\
=&\sum_{\la}\chi^{\la}_{Sp(2d)}(\z)S_\la^{\fsp,d}(\x,\y),\nonumber
\end{align}
where $\la$ is summed over all partitions of length $d$.
The proposition follows by applying to both sides of equation \eqnref{eqnSp1}
the involution of the ring of symmetric functions, which sends the
elementary symmetric functions of $\y$ to the complete symmetric
functions of $\y$.
\end{proof}

We need the following lemma to prove the main result of this
subsection.

\begin{lem}\cite{CZ}\label{lem}
Suppose that $f^\la(\x)$ and $g^\la(\x)$ are power series in the
variables $\x$ and suppose that
\begin{equation*}
\sum_{\la}f^\la(\x)\chi^\la_{Sp(2d)}(\z)=\sum_\la
g^\la(\x)\chi^\la_{Sp(2d)}(\z),
\end{equation*}
where $\la$ is summed over all partitions of length $d$. Then
$f^\la(\x)=g^\la(\x)$, for all $\la$.
\end{lem}

\begin{prop} Let $\x=\{x_1,x_2,\cdots\}$ and
${\bf y}=\{y_1,y_2,\cdots\}$ be two infinite sets of variables. For
each partition $\la$ of length $d$, we have
 $$
HS_\la^{\fsp,d}(\x,\y):=\sum_{\mu\nu}c^\la_{\mu\nu}(\z)S_\mu^{\fsp,d}(\x)
DS_\nu^{\fsp,d}(\y)
 $$
where $\mu$ and $\nu$ are summed over all partitions of length
$d$. Here the non-negative integers $c^\la_{\mu\nu}$ are the
multiplicity of $V^\la_{Sp(2d)}$ in the tensor product
decomposition of $V^\mu_{Sp(2d)}\otimes V^\nu_{Sp(2d)}$.
\end{prop}

\begin{proof}
By \propref{combin1}, we have
\begin{align*}
 &\prod_{i=1}^{d}\prod_{j=1}^{\infty}\prod_{k=\hf}^{\infty}
 {\frac{(1+x_jz_i)(1+x_jz_i^{-1})}{(1-y_kz_i)(1-y_kz_i^{-1})}}\\
=&\sum_{\mu}\chi^{\mu}_{Sp(2d)}(\z)S_\mu^{\fsp,d}(\x)
\sum_{\nu}\chi^{\nu}_{Sp(2d)}(\z)D_\nu^{\fsp,d}(\y),
\end{align*}
 where $\mu$ and $\nu$ are summed over all partitions of length $d$.
On the other hand,
\begin{equation*}
\chi^{\mu}_{Sp(2d)}(\z)\chi^{\nu}_{Sp(2d)}(\z)
=\sum_{\la}c^\la_{\mu\nu}\chi^{\la}_{Sp(2d)}(\z),
\end{equation*}
where the summation is over all partitions of length $d$. Thus
\begin{align*}
&\sum_{\mu}\chi^{\mu}_{Sp(2d)}(\z)S_\mu^{\fsp,d}(\x)
\sum_{\nu}\chi^{\nu}_{Sp(2d)}(\z)D_\nu^{\fsp,d}(\y)\\
=&\sum_{\la}\chi^{\la}_{Sp(2d)}\big(\sum_{\mu\nu}c^\la_{\mu\nu}(\z)S_\mu^{\fsp,d}(\x)
D_\nu^{\fsp,d}(\y)\big),
\end{align*}
 where $\la$, $\mu$ and $\nu$ are summed over all partitions of length
 $d$. Therefore we have
 \begin{align}\label{eqnSp2}
  &\prod_{i=1}^{d}\prod_{j=1}^{\infty}\prod_{k=\hf}^{\infty}
 {\frac{(1+x_jz_i)(1+x_jz_i^{-1})}{(1-y_kz_i)(1-y_kz_i^{-1})}}\\
=&\sum_{\la}\chi^{\la}_{Sp(2d)}\big(\sum_{\mu\nu}c^\la_{\mu\nu}(\z)S_\mu^{\fsp,d}(\x)
D_\nu^{\fsp,d}(\y)\big),\nonumber
\end{align}
 where $\la$, $\mu$ and $\nu$ are summed over all partitions of length
 $d$.
Now the proposition follows from \eqnref{eqnSp2}, \propref{HS} and
\lemref{lem}.
\end{proof}

Now we turn to the computation of the formal character of $\F_0$ with
respect to the abelian algebra $\sum_{s\in\half\N}\C
\te_{ss}\oplus\sum_{i=1}^d\C E_{ii}$. We need the following commutation
relations:  for $i\in\N$, $r\in\frac{1}{2}+\Z_+$,
\begin{align*}
&[\te_{ii},\psi^{+,p}_{-n}]=\delta_{in}\psi^{+,p}_{-n},
&[\te_{ii},\psi^{-,p}_{-n}]=\delta_{-in}\psi^{-,p}_{-n},\\
&[\te_{rr},\gamma^{+,p}_{-s}]=\delta_{rs}\gamma^{+,p}_{-s},
&[\te_{rr},\gamma^{-,p}_{-s}]=\delta_{-rs}\gamma^{-,p}_{-s},\\
&[\te_{rr},\psi^{\pm,p}_{-n}]=[\te_{ii},\gamma^{\pm,p}_{-r}]=0.
\end{align*}
Furthermore for $i=1,\cdots,d$, we have
\begin{align*}
&[E_{ii},\psi^{+,p}_{-n}]=\delta_{ip}\psi^{+,p}_{-n},
&[E_{ii},\psi^{-,p}_{-n}]=-\delta_{ip}\psi^{-,p}_{-n},\\
&[E_{ii},\gamma^{+,p}_{-r}]=\delta_{ip}\gamma^{+,p}_{-r},
&[E_{ii},\gamma^{-,p}_{-r}]=-\delta_{ip}\gamma^{-,p}_{-r}.\\
\end{align*}

Let $e$ be a formal indeterminate. For $j\in\N$, $r\in\hf
+\Z_+$, $i=1,\cdots,d$, set
\begin{equation*}
z_i=e^{\epsilon_i},\quad x_j=e^{\omega_j},\quad y_r=e^{\omega_r},
\end{equation*}
where $\epsilon_1,\cdots,\epsilon_d$ and $\omega_s$ are the
respective fundamental weights of $Sp(2d)$ and $\hC$ introduced
earlier. By using the commutation relations established above, we can easily show that
the formal character of $\F_0$, with
respect to the abelian algebra $\sum_{s\in\half\N}\C
\te_{ss}\oplus\sum_{i=1}^d\C E_{ii}$, is given by
\begin{equation}\label{character1}
{\rm ch}\F_0=\prod_{i=1}^{d}\prod_{j=1}^{\infty}
\prod_{k=\hf}^{\infty}
 {\frac{(1+x_jz_i)(1+x_jz_i^{-1})}{(1-y_kz_i)(1-y_kz_i^{-1})}}.
\end{equation}
By \propref{HS} we can rewrite \eqnref{character1} as
\begin{equation}\label{character2}
{\rm
ch}\F_0=\sum_{\la}\chi^{\la}_{Sp(2d)}(\z)HS_\la^{\fsp,d}(\x,\y),
\end{equation}
where $\la$ is summed over all partitions of length $d$. On the
other hand, \thmref{duality-hC} implies that
\begin{equation}\label{character3}
{\rm ch}\F_0\cong\sum_{\la}{\rm
ch}L(\hC,\Lambda^{\hC}(\la))\chi^{\la}_{Sp(2d)}(\z),
\end{equation}
where $\la$ is summed over all partitions of length $d$. Using
\eqnref{character3} and \eqnref{character2} together with
\lemref{lem} we have the following character formula.

\begin{thm}\label{characterthm} For each partition $\la$ of length $d$,
the formal character of the irreducible $\hC$-module $L(\hC,\Lambda^{\hC}(\la))$
is given by
\begin{equation*}
{\rm ch}L(\hC,\Lambda^{\hC}(\la))=HS_\la^{\fsp,d}(\x,\y).
\end{equation*}
\end{thm}

\subsection{Character formula for $\hD$} \label{character-hD}
We now construct a character formula for the unitarizable irreducible
quasi-finite highest weight $\hD$-modules.
Let us start by recalling some results on formal characters of
finite dimensional representations of the orthogonal group (see
\cite{FH}, \cite{ H1}, \cite{H2}).

When we deal with characters of modules over $\fso(2m)$, we will put ${\tilde
h}_{i}:=-h_{m-i+1}$ and $x_i=e^{-\epsilon_{m-i+1}}$ for
$i=1,2,\cdots,m$. That is $x_i=z^{-1}_{m-i+1}$. Recall that
$\epsilon_i( h_j)=\delta_{ij}$ and $z_i=e^{\epsilon_i}$ where
$\epsilon_i\in\h^*$ such that $\epsilon_i(h_j)=\delta_{ij}$ (see
\secref{unitarizable-hD}). For each finite sequence of complex
numbers $\la=(\la_1,\cdots,\la_m )$, we let
$W^\la_{\fso(2m)}:=V^{\la^*}_{\fso(2m)}$ where
$\la^*=(-\la_m,\cdots,-\la_1)$. Note that $\la({\tilde
h}_{i})=\la_i$ for $i=1,\cdots,m$. $W^\la_{\fso(2m)}$ is a
finite-dimensional irreducible representation if and only if
$-|\la_1|\ge\la_2\ge\cdots\ge\la_m$ with either $\la_i\in \Z$ or
$\la_i\in \hf+\Z$ for $i=2,\cdots,m$.

When considering characters of $O(n)$-modules,
we put $z_i=e^{\epsilon_i}$ with $\epsilon_i\in\h^*$ such that
$\epsilon_i(h_{j})=\epsilon_i(e_{jj}-e_{d+j,d+j})=\delta_{ij}$ for $n=2d$,
and $\epsilon_i(h_{j})=\epsilon_i(e_{jj}-e_{d+j+1,d+j+1})=\delta_{ij}$
for $n=2d+1$ (see \secref{unitarizable-hC}). For each partition $\la$ of length
$n$ and each sequence of complex numbers $\nu$ of length $m$, we
write $\chi^{\la}_{O(n)}(\z)=\chi^{\la}_{O(n)}(z_1,\cdots,z_d)$
and
$\tilde{\chi}^{\nu}_{\fso(2m)}(\x)=\tilde{\chi}^{\nu}_{\fso(2m)}(x_1,\cdots,x_m)$
for the character of $O(n)$-module $V^\la_{O(n)}$ and
$\fso(2m)$-module $W^\nu_{\fso(2m)}$ to stress their dependence on
the variables $z_1,\cdots,z_d$ and $x_1,\cdots,x_m$, respectively.

Recall that when $n=2 d+1$, the irreducible $O(n)$-modules $V^\la_{O(n)}$
and $V^{\bar\la}_{O(n)}$ restrict to isomorphic $SO(n)$-modules. To distinguish these
$O(n)$-representations at the level of characters, we let $\epsilon$ be
the eigenvalue of $-I_{n}\in O(n)$ so that $\epsilon^2=1$.
Denote by $\chi_{O(n)}^\la(\epsilon, {\bf z})$ the character of $V^\la_{O(n)}$
(with $\la'_1+\la'_2\le n$) with respect to the Cartan subalgebra
$\sum \C h_i$ together with $-I_n$. It is easy to see that
$$\chi_{O(n)}^\la(\epsilon, {\bf z})=  \epsilon^{|\la|} \chi_{SO(n)}^\la({\bf z})$$
where $|\la|$ is the size of $\la$.

By the classical $(\fso(2m), O(n))$-duality on the
exterior algebra $\La(\C^{n}\otimes \C^{m})$ with $m\ge n$
(\cite{H1}, \cite{H2}, also see \cite{CZ}), we have the following
identities for $n=2d$ and $n=2d+1$ respectively :
\begin{align}
&(x_1\cdots x_m)^{-\frac{n}{2}}
\prod_{i=1}^{d}\prod_{j=1}^m(1+x_jz_i)(1+x_jz_i^{-1})
=\sum_{\la}\chi^{\la}_{O(n)}(\z)
{\tilde\chi}^{{\la'}-\frac{n}{2}(1^m)}_{\fso(2m)}(\x),\label{even-char}\\
&(x_1\cdots x_m)^{-\frac{n}{2}}
\prod_{i=1}^{d}\prod_{j=1}^m(1+\epsilon x_j z_i)(1+\epsilon  x_j z_i^{-1})
(1+\epsilon x_j)=\sum_{\la}\chi^{\la}_{O(n)}(\epsilon, \z)
{\tilde\chi}^{{\la'}-\frac{n}{2}(1^m)}_{\fso(2m)}(\x).\label{odd-char}
\end{align}
The summations on the right hand sides of both equations range over
all partitions of length $n$ with $\la'_1+\la_2'\le n$ and $\la_1\le m$. Note that the
partition $\la'$ is considered as a partition of length $m$ and
$\frac{n}{2}\ge\la_2'\ge\cdots\ge\la_m'$ together with $n\ge
\la_1'$.

Recall that
\begin{align*}
E_r:=E_r(x_1,\cdots,x_m,x^{-1}_1,\cdots,x_m^{-1}),
\end{align*}
where $E_r$ is the $r$-th elementary symmetric polynomial for
$r\ge 0$ and $E_r=0$ for $r<0$ (see \secref{character}). For each
partition $\mu$ of length $l$, we let $|E_\mu|$ denote the
determinant of the $l\times l$ matrix with $i$-th row
\begin{equation*}
\begin{pmatrix}
E_{\mu_i-i+1}&E_{\mu_i-i+2}+E_{\mu_i-i}&E_{\mu_i-i+3}+E_{\mu_i-i-1}
&\cdots&E_{\mu_i-i+l}+E_{\mu_i-i-l+2}
\end{pmatrix}.
\end{equation*}
For any partition $\nu$ of length $m$ with $\nu_m= 0$, the
formal character of the finite-dimensional irreducible $\fso(2m)$-module
$V^\nu_{\fso(2m)}$ equals $|E_{\nu'}|$. On the other hand, if the partition
$\nu=(\nu_1,\cdots,\nu_m)$ is of length $m$ with $\nu_m\not= 0$, the
characters of the finite-dimensional irreducible representations
$V^\nu_{\fso(2m)}$ and
$V^{(\nu_1,\cdots,\nu_{m-1},-\nu_m)}_{\fso(2m)}$ are respectively equal
to (see \cite{FH})
\begin{align*}
\hf|E_{\nu'}|+\hf\left(\prod_{i=1}^m(x_i-x_i^{-1})\right)|E'_{(\nu-(1^m)
)'}|, \\
 \hf|E_{\nu'}|-\hf\left(\prod_{i=1}^m(x_i-x_i^{-1})\right)|E'_{(\nu-(1^m))' }|,
\end{align*}
where $|E'_{(\nu-(1^m))'}|$ is the determinant of the $(\nu_1-1)\times (\nu_1-1)$
matrix defined by \eqnref{E'}.  Note that $\nu-(1^m)$ is a partition,
and its transpose partition has length $\nu_1-1$.
Also for any partition $\nu=(\nu_1,\cdots,\nu_m)$,
the characters of the finite-dimensional irreducible
modules $V^{\nu+\hf(1^m)}_{\fso(2m)}$ and
$V^{(\nu_1+\hf,\cdots,\nu_{m-1}+\hf,-\nu_m-\hf)}_{\fso(2m)}$ respectively
equal to
\begin{align*}
 \hf\big(\prod_{i=1}^m(x_i^\hf+x_i^{-\hf})\big)|M^+_{\nu'}|
 +\hf\big(\prod_{i=1}^m(x_i^\hf-x_i^{-\hf})\big)|M^-_{\nu'}|,\\
 \hf\big(\prod_{i=1}^m(x_i^\hf+x_i^{-\hf})\big)|M^+_{\nu'}|
 -\hf\big(\prod_{i=1}^m(x_i^\hf-x_i^{-\hf})\big)|M^-_{\nu'}|.
\end{align*}
Hereafter, we let
$$|M^-_\mu|:=|M^-_\mu (x_1,\cdots,x_m,x^{-1}_1,\cdots,x_m^{-1})|$$
and
$$|M^+_\mu|:=|M^+_\mu (x_1,\cdots,x_m,x^{-1}_1,\cdots,x_m^{-1})|$$
denote the determinants of the $l\times l$ matrices respectively having the
$i$-th rows
\begin{equation*}
\begin{pmatrix}
E_{\mu_i-i+1}-E_{\mu_i-i}&E_{\mu_i-i+2}-E_{\mu_i-i-1}&
\cdots&E_{\mu_i-i+n}-E_{\mu_i-i-n+1}
\end{pmatrix}
\end{equation*}
and
\begin{equation*}
\begin{pmatrix}
E_{\mu_i-i+1}+E_{\mu_i-i}&E_{\mu_i-i+2}+E_{\mu_i-i-1}&
\cdots&E_{\mu_i-i+n}+E_{\mu_i-i-n+1}
\end{pmatrix},
\end{equation*}
where  $\mu=(\mu_1,\cdots,\mu_l)$ is any partition of length $l$.

Let $\x=\{x_1,x_2,\cdots\}$ be an infinite set of  variables.
Analogous to the symplectic Schur polynomials, we also have the
{\em orthogonal Schur functionof weight $\frac{n}{2}$},
$$S_\la^{\fso,\frac{n}{2}}(\x)=S_\la^{\fso,\frac{n}{2}}(x_1,x_2,\cdots),$$
defined for each partition $\la$
of length $n=2d$ with $\la_1'+\la_2'\le 2d$  by
 \begin{align*}
 S_\la^{\fso,\frac{n}{2}}(\x)
 :=
 \left \{ \begin{array}{ll} |\te_\la|,
\qquad & \mbox{if \ \  $\la'_1={\frac{n}{2}}$;}\\
\\
 \hf|\te_\la|+\hf\big(\sum_{i=0}^\infty e_i\big)
 \big(\sum_{i=0}^\infty
 (-1)^{i}e_i\big)|\te^\diamond_{\la-(1^m) }| & \mbox{if \ \
$\la'_1<{\frac{n}{2}}$;} \\
\\
\hf|\te_{\bar\la}|-\hf\big(\sum_{i=0}^\infty e_i\big)
 \big(\sum_{i=0}^\infty
 (-1)^{i}e_i\big)|\te^\diamond_{\bar\la-(1^m) }| & \mbox{if \ \
$\la'_1>{\frac{n}{2}}$,}\\
\end{array}
\right.
\end{align*}
where $|\te_\la|$ denotes the determinant of the $d\times d$
matrix whose $i$-th row is
\begin{align*}
(\te_{\la_{d-i+1}+i-1}\quad
\te_{\la_{d-i+1}+i-2}+\te_{\la_{d-i+1}+i} \quad
\te_{\la_{d-i+1}+i-3}+\te_{\la_{d-i+1}+i+1} \allowdisplaybreaks\\
 \quad \cdots\quad
\te_{\la_{d-i+1}+i-d}+\te_{\la_{d-i+1}+i+d-2}),
\end{align*}
$|\te^\diamond_\la|$ denotes the determinant of the $(d-1)\times
(d-1)$ matrix whose $i$-th row is
\begin{align*}
(\te'_{\la_{d-i}+i-1}\quad \te'_{\la_{d-i}+i-2}+\te'_{\la_{d-i}+i}
\quad
\te'_{\la_{d-i}+i-3}+\te'_{\la_{d-i}+i+1} \allowdisplaybreaks\\
 \quad \cdots\quad
\te'_{\la_{d-i}+i-d}+\te'_{\la_{d-i}+i+d-2}).
\end{align*}
Here $\te_r=\sum_{i=0}^{\infty}e_{i}(\x) e_{r+i}(\x)$,
 $\te_r'=\te_r-\te_{r-2}$ and $e_{i}(\x)$ is the $i$-th
elementary symmetric function in the infinite set of variables
$\x$. Recall that $\overline{\la}$ is a partition of length $n$
obtained from the Young diagram of $\la$ by replacing its first
column by a column of length $n-\la_1'$ (see
\secref{unitarizable-hD}).

Similarly, for each partition $\la$ of
length $n=2d+1$ with $\la_1'+\la_2'\le 2d+1$, the {\em orthogonal
Schur function} {\em of weight $\frac{n}{2}$}
$$S_\la^{\fso,\frac{n}{2}}(\x)=S_\la^{\fso,\frac{n}{2}}(x_1,x_2,\cdots)$$
 in infinitely many variables is defined by
 \begin{align*}
 S_\la^{\fso,\frac{n}{2}}(\x)
  :=
 \left \{ \begin{array}{ll}
 \hf\big(\sum_{i=0}^\infty e_i\big)|\tm^+_\la|
 +\hf\big(\sum_{i=0}^\infty
 (-1)^{i}e_i\big)|\tm^-_\la|
 & \mbox{if \ \
$\la'_1\le {\frac{n}{2}}$;} \\
 \\
 \hf\big(\sum_{i=0}^\infty e_i\big)|\tm^+_{\bar\la}|
 -\hf\big(\sum_{i=0}^\infty
 (-1)^{i}e_i\big)|\tm^-_{\bar\la}| & \mbox{if \ \
$\la'_1>{\frac{n}{2}}$,}
\end{array}
\right.
\end{align*}
where $|\tm_\la^+|$ denotes the determinant of the $d\times d$
matrix whose $i$-th row is
\begin{align*}
(\te_{\la_{d-i+1}+i-1}+\te_{\la_{d-i+1}+i}\quad
\te_{\la_{d-i+1}+i-2}+\te_{\la_{d-i+1}+i+1} \quad
  \allowdisplaybreaks\\
   \cdots\quad
\te_{\la_{d-i+1}+i-d}+\te_{\la_{d-i+1}+i+d-1}),
\end{align*}
and $|\tm_\la^-|$ denotes the determinant of the $d\times d$
matrix whose $i$-th row is
\begin{align*}
(\te_{\la_{d-i+1}+i-1}-\te_{\la_{d-i+1}+i}\quad
\te_{\la_{d-i+1}+i-2}-\te_{\la_{d-i+1}+i+1} \quad
  \allowdisplaybreaks\\
   \cdots\quad
\te_{\la_{d-i+1}+i-d}-\te_{\la_{d-i+1}+i+d-1}).
\end{align*}
We put
$S_\la^{\fso,d}(x_1,\cdots,x_m)=S_\la^{\fso,d}(x_1,\cdots,x_m,0,0,\cdots)$.
Using similar arguments as for $\fsp(2m)$ in
\secref{character}, we have
 \begin{align*}
 &(x_1\cdots x_m)^{\frac{n}{2}}\tilde{\chi}^{{\la'}-d(1^m)}_{\fso(2m)}(x_1,\cdots,x_m)
=S_\la^{\fso,d}(x_1,\cdots,x_m)
 \end{align*}
for any partition $\la$ of length $n$ with
$\la_1'+\la_2'\le n$ and $\la_1\le m$.
Now we can rewrite the combinatorial formulae \eqnref{even-char}
and \eqnref{odd-char} respectively as follows (for $m\ge n$):
\begin{align}\label{combin-even}
 \prod_{i=1}^{d}\prod_{j=1}^m(1+x_jz_i)(1+x_jz_i^{-1})
=\sum_{\la}\chi^{\la}_{O(n)}(\z)S_\la^{\fso,\frac{n}{2}}(x_1,\cdots,x_m),
\end{align}
for even integer $n=2d$ and
\begin{align}\label{combin-odd}
 \prod_{i=1}^{d}\prod_{j=1}^m(1+\epsilon x_jz_i)(1+\epsilon x_jz_i^{-1})(1+\epsilon x_j)
=\sum_{\la}\chi^{\la}_{O(n)}(\epsilon, \z)S_\la^{\fso,\frac{n}{2}}(x_1,\cdots,x_m),
\end{align}
for odd integer $n=2d+1$. In both equations the summations over $\la$ range
over all partitions of length $n$
satisfying $\la'_1+\la'_2\le n$ and $\la_1\le m$.

Let $\x=\{x_1,x_2,\cdots\}$ and
$\y=\{y_\hf,y_{1\hf},\cdots\}$ be two infinite sets of  variables.
For each partition $\la$ of length $n$, the {\em skew orthogonal
Schur function of weight $\frac{n}{2}$} of infinitely many variables
denoted by
$D_\la^{\fso,\frac{n}{2}}(\x)=D_\la^{\fso,\frac{n}{2}}(x_1,x_2,\cdots)$
is defined by $\sigma(D_\la^{\fso,\frac{n}{2}}(\x))$ where
$\sigma$ is the involution of the ring of symmetric functions
sending the elementary symmetric functions of $x_j$'s to the
complete symmetric functions of $x_j$'s. Also, for each partition
$\la$ of length $n$, the {\em hook orthogonal Schur function of
weight $\frac{n}{2}$} of infinitely many variables denoted by
$HS_\la^{\fso,\frac{n}{2}}(\x,\y)$ is defined by
$\sigma(S_\la^{\fso,\frac{n}{2}}(\x,\y))$ where $\sigma$ is the
involution of the ring of symmetric functions sending the
elementary symmetric functions of $y_j$'s to the complete
symmetric functions of $y_j$'s.

By using \eqnref{combin-even}, \eqnref{combin-odd} and similar
arguments as in \secref{character}, we can prove the following two
propositions.

\begin{prop}\label{combin1-evenodd} Let $\x=\{x_1,x_2,\cdots\}$
be an infinite set of variables and ${\bf
z}=\{z_1,z_2,\cdots,z_d\}$ be $d$ variables.
\begin{itemize}
\item[(i)] When $n=2d$,   we have
\begin{align}
&\prod_{i=1}^{d}\prod_{j=1}^\infty(1+x_jz_i)(1+x_jz_i^{-1})
=\sum_{\la}\chi^{\la}_{O(n)}(\z)S_\la^{\fso,\frac{n}{2}}(\x),\\
&\prod_{i=1}^{d}\prod_{j=1}^\infty(1-x_jz_i)^{-1}(1-x_jz_i^{-1})^{-1}
=\sum_{\la}\chi^{\la}_{O(n)}(\z)D_\la^{\fso,\frac{n}{2}}(\x),
\end{align}
where the summations on the right hand sides of both equations range over all
partitions of length $n$ satisfying $\la'_1+\la'_2\le n$.
\item[(ii)] When $n=2d+1$,  we have
\begin{align}
&\prod_{i=1}^{d}\prod_{j=1}^\infty(1+\epsilon x_jz_i)(1+\epsilon x_jz_i^{-1})(1+\epsilon x_j)
=\sum_{\la}\chi^{\la}_{O(n)}(\epsilon, \z)S_\la^{\fso,\frac{n}{2}}(\x),\\
&\prod_{i=1}^{d}\prod_{j=1}^\infty(1-\epsilon x_jz_i)^{-1}(1-\epsilon x_jz_i^{-1})^{-1}(1-\epsilon x_j)^{-1}
=\sum_{\la}\chi^{\la}_{O(n)}(\epsilon, \z)D_\la^{\fso,\frac{n}{2}}(\x),
\end{align}
where the summations on the right hand sides of both equations range over all
partitions of length $n$ satisfying $\la'_1+\la'_2\le n$.
\end{itemize}\end{prop}

\begin{prop}\label{HS-O} Let $\x=\{x_1,x_2,\cdots\}$ and $\y=\{y_\hf,y_{1\hf},\cdots\}$
be two infinite sets of variables and ${\bf
z}=\{z_1,z_2,\cdots,z_d\}$ be $d$ variables.
\begin{itemize}
\item[(i)] When $n=2d$, we have
\begin{align*}
 \prod_{i=1}^{d}\prod_{j=1}^{\infty}\prod_{k=\hf}^{\infty}
 {\frac{(1+x_jz_i)(1+x_jz_i^{-1})}{(1-y_kz_i)(1-y_kz_i^{-1})}}
=\sum_{\la}\chi^{\la}_{O(n)}(\z)HS_\la^{\fso,\frac{n}{2}}(\x,\y),
\end{align*}
 where $\la$ is summed over all partitions of length $n$.
\item[(ii)]  When $n=2d+1$, we have
\begin{align*}
 \prod_{i=1}^{d}\prod_{j=1}^{\infty}\prod_{k=\hf}^{\infty}
 {\frac{(1+\epsilon x_jz_i)(1+\epsilon x_jz_i^{-1})(1+\epsilon x_j)}{(1-\epsilon y_kz_i)
 (1-\epsilon y_kz_i^{-1})(1-\epsilon y_k)}}
=\sum_{\la}\chi^{\la}_{O(n)}(\epsilon, \z)HS_\la^{\fso,\frac{n}{2}}(\x,\y),
\end{align*}
 where $\la$ is summed over all partitions of length $n$.
\end{itemize}\end{prop}

We need the following results to prove \thmref{D-character}.

\begin{lem}\cite{CZ}\label{lem2}
Let $f^\la(\y)$ and $g^\la(\y)$ be power series in the
variables $\y$.
\begin{itemize}
\item[(i)] Suppose that $n$ is odd and
\begin{equation*}
\sum_{\la}f^\la(\y)\chi^\la_{O(n)}(\epsilon,\x)=\sum_\la
g^\la(\y)\chi^\la_{O(n)}(\epsilon,\x),
\end{equation*}
where the summation is over all partitions of length $n$. Then
$f^\la(\y)=g^\la(\y)$, for all $\la$.
\item[(ii)] Suppose that $n$
is even and
\begin{equation*}
\sum_{\la}f^\la(\y)\chi^\la_{O(n)}(\x)=\sum_\la
g^\la(\y)\chi^\la_{O(n)}(\x),
\end{equation*}
where the summation is over all partitions of length $n$.  Then
$f^\la(\y)+f^{\bar{\la}}(\y)=g^\la(\y)+g^{\bar{\la}}(\y)$.
\end{itemize}
\end{lem}

\begin{prop} Let $\x=\{x_1,x_2,\cdots\}$ and
${\bf y}=\{y_1,y_2,\cdots\}$ be two set of infinitely many variables. Let
$n$ be a fixed non-negative integer and $b^\la_{\mu\nu}$ denote
the multiplicity of $V^\la_{O(n)}$ in the tensor product
decomposition of $V^\mu_{O(n)}\otimes V^\nu_{O(n)}$.
\begin{itemize} \item[(i)] Suppose that
$n$ is odd, for any partition $\la$ of length $n$  with
$\la'_1+\la'_2\le n$, we have
 $$
HS_\la^{\fso,{\frac{n}{2}}}(\x,\y)=\sum_{\mu, \nu}b^\la_{\mu \nu}(\z)S_\mu^{\fso,{\frac{n}{2}}}(\x)
DS_\nu^{\fso,{\frac{n}{2}}}(\y)
 $$
where $\mu$ and $\nu$ are summed over all partitions of length
$n$. \item[(ii)] Suppose that $n$ is even, for any partition of
length $n=$ with $\la'_1+\la'_2\le n$, we have
 \begin{align*}
 &HS_\la^{\fso,\frac{n}{2}}(\x,\y)+HS_{\bar\la}^{\fso,\frac{n}{2}}(\x,\y)\\
 =&\big(\sum_{\mu, \nu}b^\la_{\mu \nu}(\z)S_\mu^{\fso,\frac{n}{2}}(\x)
D_\nu^{\fso,\frac{n}{2}}(\y)\big)+\big(\sum_{\mu, \nu}b^{\bar\la}_{\mu \nu}
(\z)S_\mu^{\fso,\frac{n}{2}}(\x)
D_\nu^{\fso,\frac{n}{2}}(\y)\big),
  \end{align*}
where $\mu$ and $\nu$ are summed over all partitions of length
 $n$ satisfying $\mu'_1+\mu'_2\le n$
 and $\nu'_1+\nu'_2\le n$, respectively.
 \end{itemize}
\end{prop}

\begin{proof}
We shall only prove the case with $n$ being even, as the odd case is
analogous. Let $n=2d$. By \propref{combin1-evenodd}, we have
\begin{align*}
 &\prod_{i=1}^{d}\prod_{j=1}^{\infty}\prod_{k=\hf}^{\infty}
 {\frac{(1+x_jz_i)(1+x_jz_i^{-1})}{(1-y_kz_i)(1-y_kz_i^{-1})}}\\
=&\sum_{\mu}\chi^{\mu}_{O(n)}(\z)S_\mu^{\fso,\frac{n}{2}}(\x)
\sum_{\nu}\chi^{\nu}_{O(n)}(\z)D_\nu^{\fso,\frac{n}{2}}(\y),
\end{align*}
 where $\mu$ and $\nu$ are summed over all partitions of length $n$
 satisfying $\la'_1+\la'_2\le n$ and
 $\mu'_1+\mu'_2\le n$, respectively.
 On the other hand,
\begin{equation*}
\chi^{\mu}_{O(n)}(\z)\chi^{\nu}_{O(n)}(\z)
=\sum_{\la}b^\la_{\mu\nu}\chi^{\la}_{O(n)}(\z),
\end{equation*}
where the summation is over all partitions of length $n$ with
$\la'_1+\la'_2\le n$. Thus
\begin{align*}
&\sum_{\mu}\chi^{\mu}_{O(n)}(\z)S_\mu^{\fso,\frac{n}{2}}(\x)
\sum_{\nu}\chi^{\nu}_{O(n)}(\z)D_\nu^{\fso,\frac{n}{2}}(\y)\\
=&\sum_{\la}\chi^{\la}_{O(n)}\big(\sum_{\mu, \nu}b^\la_{\mu\nu}(\z)S_\mu^{\fso,\frac{n}{2}}(\x)
D_\nu^{\fso,\frac{n}{2}}(\y)\big),
\end{align*}
 where $\la$, $\mu$ and $\nu$ are summed over all partitions of length
 $n$ satisfying $\la'_1+\la'_2\le n$, $\mu'_1+\mu'_2\le n$
 and $\nu'_1+\nu'_2\le n$, respectively. Therefore we have
 \begin{align}\label{eqnO2}
  &\prod_{i=1}^{d}\prod_{j=1}^{\infty}\prod_{k=\hf}^{\infty}
 {\frac{(1+x_jz_i)(1+x_jz_i^{-1})}{(1-y_kz_i)(1-y_kz_i^{-1})}}\\
=&\sum_{\la}\chi^{\la}_{O(n)}\big(\sum_{\mu\nu}b^\la_{\mu\nu}(\z)S_\mu^{\fso,\frac{n}{2}}(\x)
D_\nu^{\fso,\frac{n}{2}}(\y)\big),\nonumber
\end{align}
 where $\la$, $\mu$ and $\nu$ are summed over all partitions of length
 $n$ satisfying $\la'_1+\la'_2\le n$, $\mu'_1+\mu'_2\le n$
 and $\nu'_1+\nu'_2\le n$, respectively.
Now by using \propref{HS-O}, equation \eqnref{eqnO2} and \lemref{lem2},
for any partition of length $n$ with $\la'_1+\la'_2\le n$, we have
 \begin{align*}
 &HS_\la^{\fso,\frac{n}{2}}(\x,\y)+HS_{\bar\la}^{\fso,\frac{n}{2}}(\x,\y)\\
 =&\big(\sum_{\mu, \nu}b^\la_{\mu\nu}(\z)S_\mu^{\fso,\frac{n}{2}}(\x)
D_\nu^{\fso,\frac{n}{2}}(\y)\big)+\big(\sum_{\mu, \nu}b^{\bar\la}_{\mu\nu}(\z)
S_\mu^{\fso,\frac{n}{2}}(\x) D_\nu^{\fso,\frac{n}{2}}(\y)\big),
  \end{align*}
where $\mu$ and $\nu$ are summed over all partitions of length
 $n$ satisfying $\mu'_1+\mu'_2\le n$
 and $\nu'_1+\nu'_2\le n$, respectively.
\end{proof}

Now we turn to the computation of the formal characters of the Fock spaces.
Let $d=\left[\frac{n}{2}\right]$, that is,
$n=2d$ if $n$ is even and $n=2d+1$ if $n$ is odd.
Let $e$ be a formal indeterminate and set for $j\in\N$, $r\in\hf
+\Z_+$, $i=1,\cdots,d$
\begin{equation*}
z_i=e^{\epsilon_i},\quad x_j=e^{\omega_j},\quad y_r=e^{\omega_r},
\end{equation*}
where $\epsilon_1,\cdots,\epsilon_d$ and $\omega_s$ are the
respective fundamental weights of $O(n)$ and $\hD$ introduced
earlier. By using similar arguments as in \secref{character}, we can
easily show that the character of $\F_0$, with respect to the
abelian algebra $\sum_{s\in\half\N}\C \te_{ss}\oplus\sum_{i=1}^d\C
E_{ii}$, is given by
\begin{equation}\label{F-even}
{\rm ch}\F_0=\prod_{i=1}^{d}\prod_{j=1}^{\infty}
\prod_{k=\hf}^{\infty}
 {\frac{(1+x_jz_i)(1+x_jz_i^{-1})}{(1-y_kz_i)(1-y_kz_i^{-1})}}.
\end{equation}
Similarly, the character of $\Fhf_0$, with respect to the abelian algebra
$\sum_{s\in\half\N}\C \te_{ss}\oplus\sum_{i=1}^l\C E_{ii}$, is
given by
\begin{equation}\label{F-odd}
{\rm ch}\Fhf_0=\prod_{i=1}^{d}\prod_{j=1}^{\infty}
\prod_{k=\hf}^{\infty}
 {\frac{(1+\epsilon x_jz_i)(1+\epsilon x_jz_i^{-1})(1+\epsilon x_j)}
 {(1-\epsilon y_kz_i)(1-\epsilon y_kz_i^{-1})(1-\epsilon y_k)}}.
\end{equation}
By \propref{HS-O}, we can rewrite \eqnref{F-even} as
\begin{equation}\label{character-O-2} {\rm
ch}{\mf F}^{\frac{n}{2}}_0=
\sum_{\la}\chi^{\la}_{O(n)}(\z)HS_\la^{\fso,\frac{n}{2}}(\x,\y),
\end{equation}
where $\la$ is summed over all partitions of length $n$ with
$\la'_1+\la'_2\le n$. On the other hand, \thmref{duality-hD} implies
that
\begin{equation}\label{character-O-3}
{\rm ch}{\mf F}^{\frac{n}{2}}_0=\sum_{\la}{\rm
ch}L(\hD,\Lambda^{\hD}(\la))\chi^{\la}_{O(n)}(\z),
\end{equation}
where $\la$ is summed over all partitions of length $n$ with
$\la'_1+\la'_2\le n$ . Combining \eqnref{character-O-3} with
\eqnref{character-O-2}, we arrive at
\begin{equation} \sum_{\la}\chi^{\la}_{O(n)}(\z)HS_\la^{\fso,\frac{n}{2}}(\x,\y)
=\sum_{\la}{\rm
ch}L(\hD,\Lambda^{\hD}(\la))\chi^{\la}_{O(n)}(\z).
\end{equation}
In a similar way, we can also derive the following equation for $n$ odd:
\begin{equation} \sum_{\la}\chi^{\la}_{O(n)}(\epsilon, \z)HS_\la^{\fso,\frac{n}{2}}(\x,\y)
=\sum_{\la}{\rm
ch}L(\hD,\Lambda^{\hD}(\la))\chi^{\la}_{O(n)}(\epsilon, \z).
\end{equation}

Applying \lemref{lem2} to these equations,  we obtain the following character
formulae:

\begin{thm} \label{D-character} For each partition $\la$ of length $n=2d+1$ with
$\la'_1+\la'_2\le n$, we have
\begin{equation*}
{\rm ch}L(\hD,\Lambda^{\hD}(\la))
=HS_\la^{\fso,\frac{n}{2}}(\x,\y).
\end{equation*}
For each partition $\la$ of length $n=2d$ with $\la'_1+\la'_2\le
n$, we have
\begin{equation*}
{\rm ch}L(\hD,\Lambda^{\hD}(\la))+{\rm
ch}L(\hD,\Lambda^{\hD}(\bar\la))
=HS_\la^{\fso,\frac{n}{2}}(\x,\y)+HS_{\bar\la}^{\fso,\frac{n}{2}}(\x,\y).
\end{equation*}
\end{thm}

\begin{rem}
From \thmref{characterthm} and \thmref{D-character} we can easily
extract character formulae for the unitarizable quasi-finite irreducible
highest weight modules over the ${\mathfrak{so}}$ and ${\mathfrak{sp}}$
type subalgebras of $\hglone$.
\end{rem}
\bigskip

\frenchspacing

\bigskip

\noindent{\bf Acknowledgements.} We thank Shun-Jen Cheng for
discussions.
Financial support from the Australian Research Council and
the National Science Council of the Republic of
China is gratefully acknowledged.

\end{document}